\documentclass[preprint]{elsarticle}

\usepackage[margin=2.5cm]{geometry}

\usepackage{booktabs}
\usepackage{amsmath}
\usepackage{amssymb}
\usepackage{xcolor}
\usepackage{caption}
\usepackage{subcaption}
\usepackage{rotating}
\usepackage{listings}
\usepackage{algorithm}
\usepackage[noend]{algpseudocode}
\usepackage{appendix}

\usepackage{float}

\usepackage{soul} 

\newcommand{\commentout}[1]{}

\graphicspath{{pics/}}

\newcommand{\qb}{\mathbf q}

\newcommand{\Fm}{\mathcal F}
\newcommand{\Sm}{\mathcal S}
\newcommand{\Cm}{\mathcal C}

\journal{Atmospheric Research}

\begin{document}

\begin{frontmatter}

\title{Comparison of adaptive mesh refinement techniques for numerical weather prediction.}
\author{Daniel S. Abdi\corref{cor1}}
\author{Ann Almgren\corref{cor2}}
\author{Francis X. Giraldo\corref{cor3}}
\author{Isidora Jankov\corref{cor4}}
\cortext[cor1]{CIRES, University of Colorado Boulder, CO, USA Email: daniel.abdi@noaa.gov}
\cortext[cor2]{Lawrence Berkeley National Laboratory, Berkeley, CA, USA}
\cortext[cor3]{Naval Postgraduate School, Monterey, CA, USA}
\cortext[cor4]{NOAA Global Systems Laboratory, Boulder, CO, USA}

\begin{abstract}
This paper examines the application of adaptive mesh refinement (AMR) in the field of numerical weather prediction (NWP). We implement and assess two distinct AMR approaches and evaluate their performance through standard NWP benchmarks. In both cases, we solve the fully compressible Euler equations, fundamental to many non-hydrostatic weather models.

The first approach utilizes oct-tree cell-based mesh refinement coupled with a high-order discontinuous Galerkin method for spatial discretization. In the second approach, we employ level-based AMR with the finite difference method. Our study provides insights into the accuracy and benefits of employing these AMR methodologies for the multi-scale problem of NWP. Additionally, we explore  essential properties including their impact on mass and energy conservation. Moreover, we present and evaluate an AMR solution transfer strategy for the tree-based AMR approach that is simple to implement, memory-efficient, and ensures conservation for both flow in the box and sphere.

Furthermore, we discuss scalability, performance portability, and the practical utility of the AMR methodology within an NWP framework -- crucial considerations in selecting an AMR approach. The current de facto standard for mesh refinement in NWP employs a relatively simplistic approach of static nested grids, either within a general circulation model or a separately operated regional model with loose one-way synchronization. It is our hope that this study will stimulate further interest in the adoption of AMR frameworks like AMReX in NWP. These frameworks offer a triple advantage: a robust dynamic AMR for tracking localized and consequential features such as tropical cyclones, extreme scalability, and performance portability.
\end{abstract}


\end{frontmatter}



\section{Introduction}
The study of atmospheric phenomena plays a fundamental role in enhancing our understanding of the Earth's weather and climate system and its consequential impacts on human society and the environment. Atmospheric scientists strive to develop precise models capable of simulating intricate atmospheric processes with high fidelity, empowering us to make informed decisions regarding weather forecasting, climate change, and air quality management. Nevertheless, modeling the atmosphere presents distinctive challenges owing to the extensive range of spatial and temporal scales involved, as well as the imperative need to capture interactions among various physical and chemical processes \citep{wallace2006}. 

The employment of uniform grids to address multiscale problems, which encompass small flow scales over expansive domains, becomes infeasible, even when using high-performance supercomputers. NWP models have traditionally relied on uniform grids, with resolution chosen a priori to strike a balance between computational efficiency and accuracy. However, this approach has limitations when dealing with phenomena characterized by sharp gradients, localized features, or widely varying scales. In such cases, fixed grids may lead to exorbitant computational costs or a failure to resolve essential atmospheric details. A typical compromise in NWP is to operate Global Circulation Models (GCMs) with coarse resolution in the range of 10km to 30km (e.g. the Global Forecast System has 13km resolution), and regional models with a  much finer resolution (e.g. High Resolution Rapid Refresh (HRRR) at 3km). However, regional models, also known as Limited Area Models (LAMs), often employ static grids, which may hinder their ability to effectively track phenomena such  as approaching tropical cyclones. Moreover, LAMs rely on GCMs for lateral boundary conditions, often employing one-way synchronization, which may not be able to adequately capture large scale feedback triggered by localized phenomena such as tropical cyclones \citep{ferguson2016}.

Adaptive Mesh Refinement (AMR) has emerged as a powerful computational technique to address these challenges in atmospheric modeling. AMR enables dynamic grid refinement and coarsening in response to the evolving solution, optimizing computational resources where they are most needed. This adaptability makes AMR well-suited for modeling atmospheric processes that exhibit significant spatial and temporal variability, such as tropical cyclones, convective storms, turbulence, and the interaction between land surfaces and the atmosphere.  We note that static mesh refinement is a subset of AMR; if the grid hierarchy does not change between time steps then there is no additional cost associated with the adaptivity (other than the cost of assessing whether any change was needed, and this step can be eliminated as well).

In the context of high-order spectral element methods, mesh refinement can be classified into three categories a) h-refinement, where elements in regions of interest (ROI) are split into smaller child elements. b) p-refinement where the order of polynomial interpolation in an element within the ROI is increased while the grid remains the same, and  c) r-refinement where grid points are re-located to the ROI. The drawback of r-refinement is that it may result in a mesh of poor quality with high skewness and aspect ratio. 

Our primary focus in this work is h-refinement as applied to the finite volume (FV) / finite difference (FD)  and the high-order discontinuous Galerkin  spectral element methods (dGSEM) used in NWP. 

For the purposes of this paper, we introduce a taxonomy of AMR strategies for structured grids that refer only to the creation and structure of the grid hierarchy (independent of the governing equations, choice of methodology and coupling strategy).  The first is cell-based refinement; with this strategy, each cell that meets the specified refinement criteria is individually split into four (in 2d) or eight (in 3d) elements, and the resulting relationship is stored in a quad- or oct-tree structure.  Highly scalable libraries exist for this approach such as p4est \citep{p4est} and libMesh \citep{libMesh}. The second approach is a generalization of the cell-based approach that still uses a tree structure but requires the individual regions of the mesh to be logically rectangular; we typically refer to these regions as ``grids" or ``patches" or ``boxes".  This strategy is also known as ``patch-based refinement" since fine patches of some minimum size in each dimension are constructed.   Taking a broader perspective, one could in fact view cell-based refinement as a subset of the quad-tree or oct-tree patch-based refinement, wherein the minimum grid width is one. Notably, with the use of high-order spectral element methods, this distinction becomes less clear, as the element itself can be considered as a block of nodes, while typical low order cell-centered finite-volume methods have just one node per element.  In software such as FLASH \cite{fryxell2000}, if a single cell of a patch is tagged for refinement then the whole patch is refined, ensuring that any coarse grid is either entirely covered or entirely uncovered, where a cell at one level of refinement is referred to as ``covered" if there are fine cells covering the same physical location.  The grids at all levels are typically required to be the same size.  If the minimum grid width is too large then this can generate overly large regions of refinement; if it is too small then the area-to-volume ratio of the patches results in potentially overly high ratio of communication to computation costs.  The third type of AMR we consider here is also patch-based, in that data is always contained in boxes of at least a minimum width.  This third type of refinement does not use a tree-structure, rather the data is considered to live at ``levels of refinement" (the level is defined by the mesh resolution at that level). The patches are constructed to be the smallest boxes containing the tagged cells while maintaining a logically rectangular shape, and need not all have the same size.  This approach is especially convenient when using subcycling in time, since the data structures are optimized to consider communication between grids of the same resolution rather than coarse/fine communication.  For the purposes of this paper, we will refer to the three types as cell-based,  quad-tree/oct-tree patch-based, and level-based patch-based refinement.  These last two types also fall in the category of ``block-structured refinement."

\begin{figure}
\centering
\begin{subfigure}{0.49\textwidth}
\centering
\includegraphics[width=0.6\linewidth]{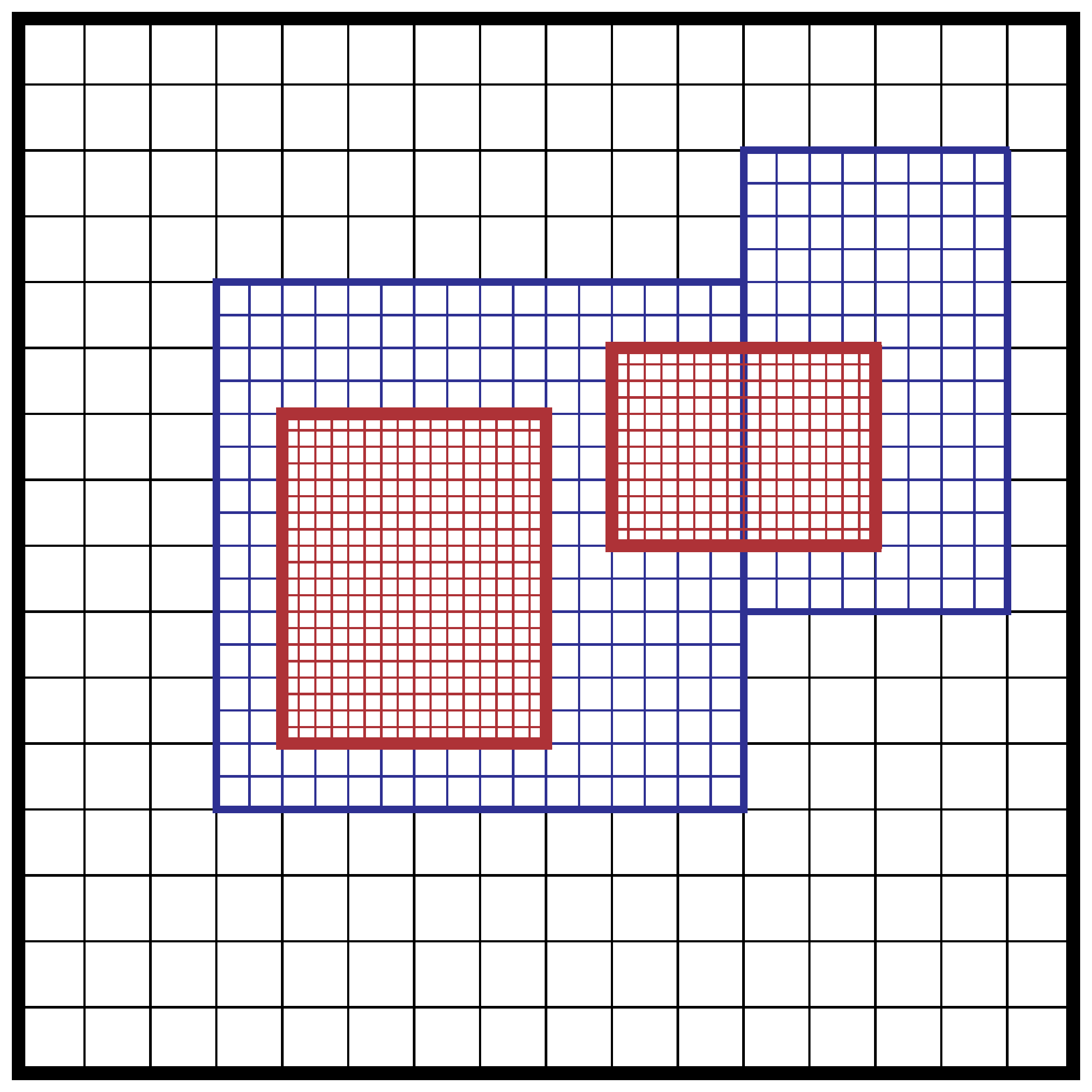}
\end{subfigure}
\begin{subfigure}{0.49\textwidth}
\centering
\includegraphics[width=\linewidth]{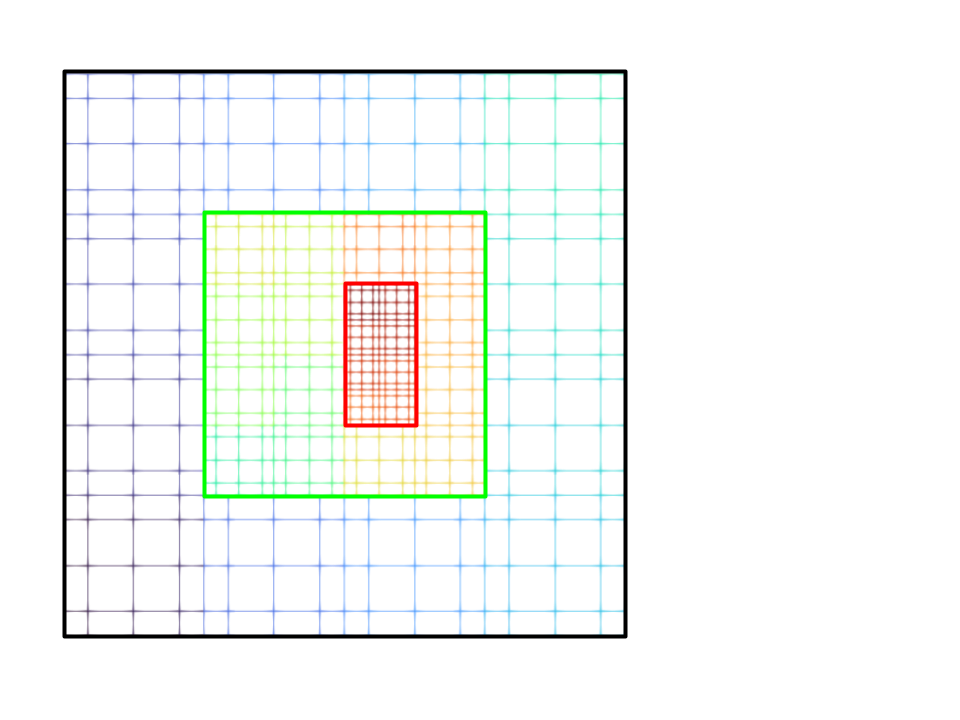} 
\end{subfigure}
\caption{Illustration of level-based and tree-based AMR. The left figure depicts level-based AMR as used in AMReX (Courtesy of Almgren et al.) \citep{amrex2021}. Three overlapping levels are depicted: level-0 (gray), level-1 (blue) and level-2 (red). Note that the boxes within a given level are non-overlapping even though the different levels themselves overlap. The right figure depicts tree-based AMR  as used in NebulaSEM \citep{NebulaSEM}. Here, the blue, green and red boxes all lie within the same level and are also non-overlapping.}
\label{fig:amr-grids}
\end{figure}

In this study, we exclusively explore the cell-based and level-based patch-based AMR methods, ignoring the tree-based patch-based method. For conciseness, we will henceforth denote the cell-based method as the tree-based method and the level-based patch-based method as the level-based method throughout the remainder of the paper. An illustration of these two methods is presented in Figure~\ref{fig:amr-grids}.

Level-based AMR was first introduced in \citet{berger1984} in which a system of conservation laws was solved on a multi-level grid hierarchy with a finite volume approach. Richardson extrapolation was used to estimate the local truncation error and define the refinement criteria 
which determined the regions of refinement, and subcycling in time was used to maintain the ratio of timestep to mesh spacing across levels. 
The first application of this method for an NWP application was by \citet{skamarock1989} in which they solved the primitive hydrostatic equations using the finite difference method and the method of \citet{berger1984} for adaptive mesh refinement.
Level-based AMR is popular because it permits the use of existing solvers for locally uniform meshes. 
A number of open-source libraries for level-based AMR (such as AMReX, {AMRClaw}, Chombo, FLASH, Uintah, et al.) exist. In this work, we use AMReX \citep{amrex2021}, a software framework for massively parallel block-structured AMR. 

The remainder of the paper is organized as follows. In Section \ref{sec-models}, we provide a brief description of the atmospheric models used in this study. Section \ref{sec-amr} provides a description of the tree-based and level-based AMR approaches used in this work. Section \ref{sec-results} discusses various numerical tests to validate the implementation and accuracy of AMR.  Section \ref{sec-choice} is dedicated to examining the suitability of the two AMR approaches for NWP applications. Finally, our conclusions are presented in Section \ref{sec-conclusions}.

\section{Model description}
\label{sec-models}
In this study, we employ two dynamical cores with distinct AMR implementations. The first one, developed in-house and called NebulaSEM \citep{NebulaSEM}, is a versatile Computational Fluid Dynamics (CFD) code that supports both the finite volume  and discontinuous Galerkin spectral element discretizations. It provides support for tree-based AMR on a generic polyhedral mesh for finite-volume discretization, and quad-/oct-tree based AMR on hexahedral mesh for dGSEM. It offers both implicit and explicit time discretizations options. The second code,  Energy Research and Forecasting (ERF) \citep{almgren2023}, utilizes the AMReX \citep{amrex2021} framework that is both performance portable and highly scalable. ERF uses standard finite volume / finite difference discretizations on the classic Arakawa C-grid. The time discretization in ERF is a third-order Runge-Kutta scheme with substepping of perturbation quantities at the acoustic time scale. Both codes have the capability to solve the fully compressible Euler equations for dry dynamics, which form the foundation of many non-hydrostatic dynamical cores. While ERF is primarily a mesoscale model written following the principles of the well-known Weather Research and Forecasting (WRF) model, NebulaSEM can also simulate global atmospheric circulation on a cubed-sphere grid.

\subsection{Governing equations}
The prognostic variables for dry dynamics, denoted as  $\qb = {(\rho, \mathbf U^\top, \Theta, \mathbf \Cm)}^\top$, encompass critical parameters for atmospheric modeling. In this context, $\rho$ represents density, $\mathbf{U}=\rho \mathbf u$ is momentum,  $\Theta=\rho \theta$, where $\theta$ is the potential temperature, and $\mathbf \Cm = \rho \mathbf c$, where $\mathbf c$ is a vector of tracers, and the superscript $(\cdot)^\top$ denotes the transpose operator. The velocity components are $\mathbf u=(u,v,w)^\top$. The governing equations for the hyperbolic Euler equations in conservation form are expressed as follows:

\begin{equation}
\begin{aligned}
& \frac{\partial\rho}{\partial t} + \nabla \cdot \mathbf{U}  = 0\\
& \frac{\partial \mathbf U}{\partial t} + \nabla \cdot \left( \frac{ \mathbf U \otimes \mathbf U}{\rho} + P\mathbf I_3 \right)= \mathbf F\\
& \frac{\partial \Theta}{\partial t} + \nabla \cdot \left( \frac{\Theta \mathbf U}{\rho} \right) = 0\\
& \frac{\partial \Cm}{\partial t} + \nabla \cdot \left( \frac{\Cm \mathbf U}{\rho} \right) = 0.
\end{aligned}
\end{equation}
These equations can be succinctly expressed in compact vector notation as:

\begin{equation}
\label{eulervector}
\frac{\partial\qb}{\partial t} + \nabla \cdot \Fm(\qb) = \Sm(\qb).
\end{equation}
Here, $\qb$ is the solution vector, $\Fm = {(\mathbf U, \frac{\mathbf U \otimes \mathbf U}{\rho} + P \mathbf I_3, \frac{\Theta \mathbf U}{\rho},  \frac{\mathbf \Cm \mathbf U}{\rho})}^\top$ is the flux vector, $\Sm(\qb)={(0,\mathbf F,0,0)}^\top $ is the source vector where $\mathbf F$ encompasses various body forces including buoyancy and the Coriolis force, and $\mathbf I_3$ is the rank-3 identity matrix. In order to close the above system of equations, the equation of state

\begin{equation}
P=P_0 \left( \frac{R \Theta}{P_0} \right)^{\gamma}
\end{equation}
is incorporated where $P_0$ represents the baseline sea-level pressure, $R = c_p - c_v$ denotes the air's gas constant defined as the difference between specific heat at constant pressure ($c_p$) and specific heat at constant volume ($c_v$), and $\gamma = \frac{c_p}{c_v}$ represents the specific heat ratio.

To enhance numerical stability, the density and pressure variables are defined as perturbations from a hydrostatically stratified background state

\[
\begin{aligned}
\rho(\mathbf x, t)=\overline\rho(z) + \rho'(\mathbf x, t)\\
P(\mathbf x, t)=\overline P(z) + P'(\mathbf x, t)
\end{aligned}
\]
where  ($\mathbf x$,t) are the space-time coordinates. 
This formulation allows discrete mass conservation due to the flux form of the continuity equation.  However, total energy 
is not discretely conserved due to the choice of potential temperature rather than the total energy as a prognostic variable.

These equations essentially govern the dry dynamics employed by both codes, with some minor differences, 
such as the optional use of terrain-fitted coordinates and map projection factors in ERF. Additionally, we are neglecting
turbulence, viscous diffusivity and other forcing terms for the sake of simplicity.

The aforementioned set of equations is frequently employed in atmospheric modeling. While the continuity and momentum equations adopt standard conservative formulations, the thermodynamic equation is expressed in terms of potential temperature, which presents advantages and drawbacks (see \citet{reddy2023} for details). Consequently, the equation set fails to conserve total energy in a closed system. However, it is possible to conserve energy using a modified thermodynamic equation that utilizes total energy $E$ as the prognostic variable.

\subsection{The dGSEM formulation}

The discontinuous Galerkin spectral element method in NebulaSEM is briefly introduced here.  Consider the decomposition of a domain $\Omega \subset \mathbb{R}^3$ into non-overlapping elements $\Omega = \bigcup_e \Omega_e$. Starting from Equation \ref{eulervector}, the weak form is derived by applying the standard Galerkin method. This method involves multiplying the equation by a test function $v$ and integrating by parts leading to the weak form

\begin{equation}
\left(v, \frac{\partial \qb}{\partial t}\right)_{\Omega_e} + \left\langle v \mathbf n, \Fm^* (\qb)\right\rangle_{\partial\Omega_e} - \left(\nabla v, \Fm(\qb)\right)_{\Omega_e} = \left(v, \Sm(\qb)\right)_{\Omega_e}
\end{equation}
where $\mathbf n$ is the outward pointing normal on the surface of the element, $()$ and $\langle\rangle$ represent the inner products on the volume and surface, respectively. Applying a second integration by parts on the flux integral, we obtain the strong form

\begin{equation}
\left(v, \frac{\partial \qb}{\partial t}\right)_{\Omega_e} + \left\langle v \mathbf n, \Fm^* (\qb) - \Fm(\qb)\right\rangle_{\partial\Omega_e} + \left(v, \nabla \cdot \Fm(\qb)\right)_{\Omega_e} = \left(v, \Sm(\qb)\right)_{\Omega_e}
\end{equation}
The numerical flux $\Fm^*$ is computed using the generalized upwind method of Rusanov, which is appropriate for hyperbolic systems

\begin{equation}
\Fm^*(\qb) = \{\Fm(\qb)\} + \mathbf{\hat{n}}\frac{|\widehat{c}|}{2} [\![\qb ]\!]
\end{equation}
where $|\widehat{c}|$ represents the speed of sound, $\{\}$ denotes an average and $[\![ ]\!]$ signifies a jump across a face.

The element-based dGSEM method possesses several desirable characteristics that make it a compelling choice for NWP. These attributes include arbitrarily high-order accuracy, greater geometric flexibility compared to global spectral methods, an easily invertible diagonal mass matrix that facilitates solution by explicit methods, high scalability due the high arithmetic intensity per element and minimal ghost layer communications, amenability to GPU acceleration (see f.i. \citet{abdi8,abdi18} for GPU accelerated dGSEM weather model), retaining conservation property of any physical conservation law (e.g. mass and total energy) in discrete form, as well as providing a straightforward avenue for the application of hp-refinement.

\subsection{The FD formulation on the Arakawa C-grid}

ERF uses the Arakawa C-grid, a popular staggered grid arrangement in NWP. In this grid, the normal components of the wind  $(u,v,w)$ are calculated on faces, while other quantities such as $\rho,\Theta, P$  are calculated at the centers of cells. 

Advection in ERF allows users the choice of standard spatial discretization techniques 
ranging from $2^{nd}$ to $6^{th}$ order, as well as blended $3^{rd}/4^{th}$, 
blended $5^{th}/6^{th}$, and $3^{rd}$- and $5^{th}$-order 
Weighted Essentially Non-Oscillatory (WENO) methods, for the primitive variables (e.g. velocity and potential temperature)
in order to form the fluxes.
Since the fluxes in the continuity equation (which is used to update density) are the momenta $\mathbf{U}$, 
which are themselves solution variables, there is no need for additional interpolation.

The second-order discretization, for example to create the flux $\mathbf{U} q$ at the face $m-1/2$, has the form 
\[
\begin{aligned}
\qb^{2nd}_{m-1/2}  = \frac{1}{2} (\qb_{m} +  \qb_{m-1}).
\end{aligned}
\]
while the standard higher-order schemes take the form
\[
\begin{aligned}
\qb^{4th}_{m-1/2}  = \frac{7}{12} (\qb_{m} +  \qb_{m-1}) -  \frac{1}{12} ((\qb_{m+1} - \qb_{m-2})) \\
\qb^{3rd}_{m-1/2}  = \underbrace{\qb^{4th}_{m-1/2}}_{\text{centered flux}} + \underbrace{ sign(\mathbf u)\frac{1}{12} (\qb_{m+1} -  \qb_{m-2}) -  3 * (\qb_{m} +  \qb_{m-2})}_{\text{upwind bias}}. \\
\end{aligned}
\]

Compared to high-order element-based dGSEM methods \citep{giraldoBook,hesthavenBook}, high-order FD methods present several challenges related to boundary treatment. High-order interpolation schemes need to be modified  to ensure accuracy and stability of numerical simulations at boundaries comparable to that in the interior of the domain. 

\section{Adaptive Mesh Refinement}
\label{sec-amr}

In this section, we provide some insights into the adaptive mesh refinement techniques used in this study, specifically, the tree-based and level-based AMR strategies.

\subsection{Tree-based AMR}
\label{sec-tree-amr}
Tree-based AMR approaches utilize a single multi-resolution locally structured, globally unstructured grid. Quad- or oct-trees partition the domain into non-overlapping patches. Purely cell-based approaches refine cells individually, whereas block-based approaches refine a block of cells (e.g. 8x8 cells) together. In the context of high-order dGSEM, this distinction is blurred because an element contains many internal nodes.
 
\begin{figure}
\centering
\includegraphics[width=0.6\linewidth]{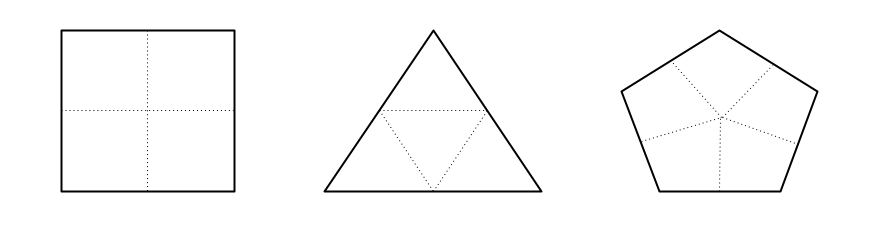}
\caption{Supported face refinement patterns. Quadrilateral and triangular faces are sub-divided into 4 elements, while $n$-sided polygon is sub-divided to $n$ quadrilaterals.}
\label{faces}
\end{figure}

The AMR implementation within NebulaSEM, originally designed for a finite-volume solver, is tailored to handle polyhedral elements with arbitrary polygonal face shapes. The process of refining polyhedral cells is executed by initially refining the faces enclosing the volume. Subsequently, pyramids are constructed, with each refined face serving as the base and the centroid of the polyhedra as the apex. The patterns for face refinement used in this work are depicted in Figure \ref{faces}. Quadrilateral and triangular faces are subdivided in the traditional manner to generate four smaller elements, while faces with $n$-sided convex polygonal shapes are split into $n$ quadrilaterals. Afterward, the pyramids that share faces are merged. Following this rather convoluted procedure, it is possible to arrive at the conventional splitting of a hexahedron into 8 child elements, while at the same time providing a generic but basic algorithm for polyhedral cell refinement.

A tree structure is employed to manage relationships between refined cells and their children. At the beginning of simulations, each cell in the potentially unstructured coarse mesh is treated as the root of a tree -- known as the ``forest-of-trees" approach in \citet{p4est}. 
In our implementation, cells designated for refinement or coarsening are promptly removed once the new finer or coarser cells, respectively, are created, retaining only the active cells to minimize storage requirements. 
It is worth noting that the current implementation is not highly scalable, as grid refinement and coarsening are performed on a single node. Subsequently, the grid is partitioned using the METIS library and distributed to slave compute nodes, each storing a portion of the refined grid. While highly scalable implementations do exist for quad-/oct-tree refinement (e.g. \citet{p4est}), the absence of a generic implementation for polyhedral refinement  was what motivated our undertaking. Admittedly outsourcing complex AMR infrastructure to external libraries is often the right approach. In the future, we intend to support the recently released extension of p4est for hybrid meshes, known as $t8code$ \citep{holke2018}. Operations such as parallel mesh refinement, coarsening, and re-partitioning pose significant computational demands. Highly scalable implementations often utilize Space Filling Curves (SFCs)\citep{bader2012} to store and partition the domain in an efficient but slightly less optimal way than algorithms in Par/METIS.  The saving in memory and run-time from not using unstructured grid partitioning method such as METIS can be significant for large scale simulations \citep{p4est}.

The existing implementation of the tree-based AMR lacks the incorporation of subcycling in time. Instead, a uniform time step is employed throughout the entire grid, with the smallest time step being enforced for every grid cell. This stands as a significant challenge in contrast to level-based AMR methods, where the integration of subcycling is more straightforward and intuitive. The absence of subcycling in time  leads to inefficiencies in terms of computational resources and hinder the method's ability to fully exploit the benefits of adaptive mesh refinement.

\subsection{Level-based AMR}
Level-based AMR is the more prevalent choice for meteorological applications (see f.i. \citet{skamarock2019,harris2013,zangl2022} ). This approach enables the reuse of existing solvers for structured, locally uniform grids. Limited area models (LAMs) can be viewed as a straightforward implementation of static (non-adaptive) mesh refinement where information flows one-way from GCM to LAM via lateral boundary conditions. In fact, LAMs often utilize an entirely different codebase with distinct physics options. One-way nested simulations may also be executed lockstep with the coarse grid simulation if both share the same codebase. In a two-way nesting approach, as implemented in the Weather Research and Forecasting (WRF) \citep{skamarock2019} model, additional feedback happens from the child to the parent for two-way communication. 
The coarse grid solution of the parent is replaced with the suitably averaged fine grid solution of the child.

\subsubsection{Grid Creation and Regridding}

AMReX, the block-structured AMR library used to evaluate level-based AMR, operates with a multi-level grid hierarchy structure  as shown in Figure~\ref{fig:amr-grids}. 
The shape of the nested region in AMReX-based codes such as ERF is more general than that allowed in WRF \citep{skamarock2019}. 
AMReX requires that the union of grids at a specific level $\ell > 0$ 
must be entirely contained in the coarser level $\ell-1$; however, 
unlike WRF there is no strict parent-child relationship between the grids on different levels,  
and no one grid at level $\ell$ must be contained within a single grid at level $\ell-1.$ {\it Proper nesting} enforces that
there is sufficient distance between the level $\ell / (\ell+1)$ interface and the level $(\ell-1) / \ell$ interface that 
interpolation from level $\ell$ to level $\ell+1$ is sufficient to fill the level $\ell+1$ ghost cells; we never need to reach
into level $\ell-1$ data for this operation.

New logically rectangular grids are dynamically created and destroyed from a list of cells tagged for refinement and coarsening using the Berger-Rigoutsos clustering algorithm \citep{berger1991}. The algorithm's main goal is to minimize inter-grid communication. AMReX imposes additional constraints on the optimization process through the grid size, including maximum length in each direction and divisibility by a blocking factor. The grid efficiency parameter allows control over how tightly the new grid boxes should fit the cells tagged for refinement. This combination of constraints and parameters ensures that grid management is both effective and efficient.

In using AMR it is important to assess the solution for the need to regrid (i.e. potentially create a new set of refined patches) 
often enough that the important dynamic features of the flow do not leave the refined region.  One could do this 
every time step, but in practice we typically regrid at regular intervals, e.g. every $T$ timesteps. In order to ensure 
that the features don't leave the fine region during that time, when creating the fine grids we ensure that there are 
at least $B$ cells at the coarse level between the feature and the coarse-fine interface.  
For compressible flows for which the timestep is determined by the acoustic rather than advective time scale, $B$ need not
be particularly large.  For example, for a flow with Mach number 0.1 running with an acoustic CFL constraint of 1, it would
take 10 timesteps for a feature to move one grid cell, thus setting $T = 10$ and $B = 2$ would be very reasonable. We note
that both $T$ and $B$ are runtime parameters specified by the user. These two parameters are also utilized in the tree-based AMR approach.

\begin{figure}[h]
\centering
\includegraphics[width=0.4\linewidth]{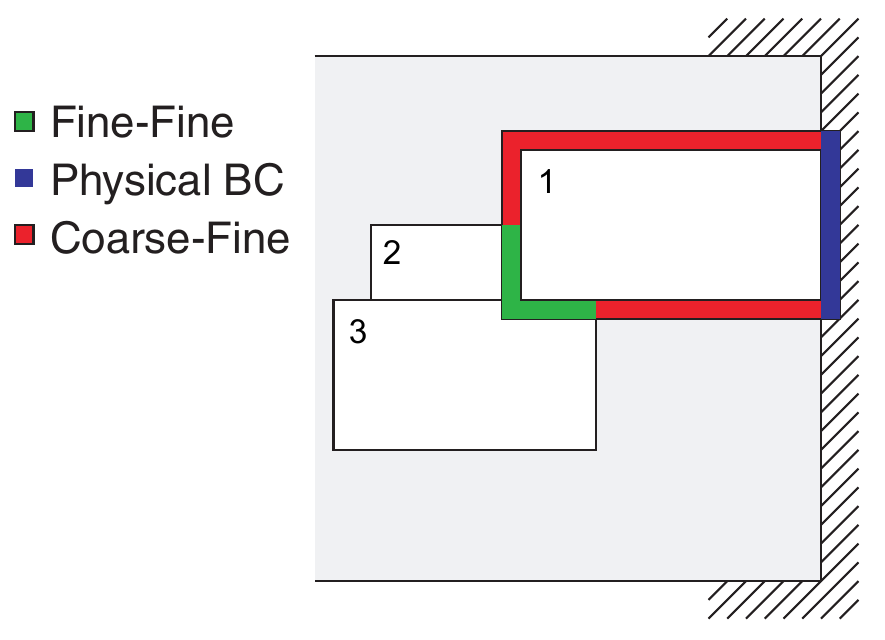}
\caption{Illustration of how ghost cells are filled in AMReX's block-structured multi-level grid hierarchy (Courtesy of Almgren et al). This figure shows two levels, with three grids at the finer level. The colored regions here show the ghost cells of grid 1 and the 
different ways in which the data are filled in those regions.  Cells in the green region are ``valid'' cells of grids 2 and 3; 
values in this region are copied from grids 2 and 3 to level 1.  Cells in the red region do not live in any other fine grid; data here 
is interpolated from the next coarser level. Finally, grids in the blue region are outside the domain and are filled according to the 
domain boundary conditions.}
\label{fig:amrex-ghost-cells}
\end{figure}

\subsubsection{Filling ghost cells}

AMReX provides general functionality for filling ghost cells (halo cells) in any given grid; the application code must supply any
problem-specific information.   The paradigm for filling ghost cells is shown in Figure~\ref{fig:amrex-ghost-cells} and is as follows: 
(1) if the ghost cell to be filled is contained in the domain and does not overlap another grid at the same level, the data is filled by interpolation from the next coarser level, else 
(2) if the ghost cell is contained in the domain and does overlap another grid at the same level, the data is filled by a copy from the other grid, else 
(3) if the ghost cell is not contained in the domain it is filled using problem-specific boundary conditions as specified in the application code such as ERF, which may be general such as ``symmetry" or may require (possibly time-interpolated) data from external files, analogous to how WRF treats real boundary conditions.  
Periodic boundary conditions are considered to be of type (2) above.
We also note that while AMReX generically views boundary conditions as only modifying values on faces on the domain boundary or cells outside the domain,
AMReX does not prohibit the type of nudging used in WRF and other atmospheric modeling codes, in which cells interior to the grid are modified using external data.  (Specifically, ERF can read wrfbdy files and set domain boundary conditions identically to the way WRF does.)


The spatial resolution ratio between consecutive level grids is often set to a small integer, typically 2, 3 or 4. (Refinement ratios of 2 or 4 are typically used in applications in which all of the data is at cell centers or corners; refinement ratio of 3 appears mostly in codes with staggered meshes.)  More general integer refinement ratios are allowed by AMReX if an interpolation stencil is provided.

\subsubsection{Subcycling in time}
Subcycling in time, i.e. keeping the ratio of the time step to mesh spacing the same across levels, is utilized in many applications 
in order to keep the Courant Friedrichs Levy (CFL) number $CFL = \frac{c\Delta t}{\Delta x}$ the same.  (Here $c$ is the sound speed
and this is the acoustic rather than advective CFL number.)
After each coarse grid time step, and possibly multiple subcycled fine grid steps, a synchronization process is performed to allow feedback from data at the fine level to data at the coarse level. The synchronization consists of two steps. For systems of equations written in flux form, a
refluxing operation is performed that effectively over-writes the faced-based fluxes used in the evolution of cell-centered quantities 
on the coarse level by the time- and space-averaged face-based fluxes from the fine level.  
(In practice this takes the form of updating the data in all coarse cells immediately adjacent to, but not covered by, 
grids at the fine level.)  Refluxing guarantees that if a cell-centered quantity is conserved by the numerical 
algorithm on a single level (due to the flux form of the equations), then this quantity is
conserved in the hierarchical solution in a simulation with AMR.
However, the fine and coarse solution can still be inconsistent with each other; to address this
we typically average the fine solution onto the coarse regions underneath.

Moreover, AMReX is purposefully designed to leverage the capabilities of modern high-performance computing (HPC) architectures. It utilizes the Message Passing Interface (MPI) for coarse-grained parallelization making it suitable for large-scale simulations distributed across multiple nodes. Additionally, AMReX provides several backends to achieve fine-grained parallelization, including OpenMP for CPU, and CUDA/HIP/SYCL for GPU acceleration. Each grid resides entirely in an MPI rank, which can also host other patches. Fine-grained parallelization for computations inside a grid is enabled in such a way that threads of execution process specific regions or logical ``tiles" of the grid. This combination of coarse-grained and fine-grained parallelization strategies ensure that AMReX is well-suited for exascale simulations.

Load balancing in AMReX uses one of two methods \citep{amrex2021}: space filling curves, similar to what's used in tree-based methods, or solving the knapsack problem \citep{hochbaum1987,amrex2021} where the weights 
can be determined by the number of cells in a patch or by a user-defined function. 
Load balancing is applied at each level independently since the levels are advanced sequentially.

\subsection{Solution transfer and non-conformal face flux computations}
A common concern with the non-overlapping tree-based AMR methods and high-order dGSEM method is whether the refined grid should be conformal or non-conformal. In a conformal grid, the faces align perfectly, which means the dG  spectral element solver does not require specific adjustments to handle AMR. On the other hand, maintaining a conformal grid after refinement is a challenging task and often leads to poor quality meshes. It is important to note that this is not a concern for a Finite Volume (FV) solver because the method operates the same way regardless of the number of faces shared between elements. 

Adaptive mesh refinement necessitates the transfer of solution between different levels of refinement. In NWP applications, it is important for the interpolation method to conserve essential quantities such as the total mass within the system. Therefore, the transfer of solutions from parent to child elements and across faces between elements should be performed in a manner that preserves the global conservation properties of dGSEM. One such method that ensures conservation for AMR with the dGSEM discretization is the ``mortar method" discussed in \citet{kopriva1996,kopera2014,kozdon2018}.

While the mortar method is unrivaled in terms of accuracy, yielding results identical to that of a conformal grid, it poses implementation difficulties associated with generating appropriate mortar elements, especially in three dimensions. An alternative approach that may be easier to implement is the simple point-to-point interpolation method discussed in \citet{laughton2021}. The point-to-point interpolation is attractive in this regard, however, it lacks mass conservation and requires dealiasing to counter aliasing and oscillatory effects introduced by discontinuities in interpolation at element interfaces. Under-resolved simulations are particularly prone to admitting discontinuities and introducing spurious noise into elements.

AMR introduces two types of error: 1) \emph{solution transfer error}, which is the error introduced during the transfer of fields from the parent grid/elements to the refined grid. Addressing this type of error can be challenging  on curved surfaces, as will be discussed later. 2) \emph{flux computation error}, which is the error introduced when computing fluxes at non-conformal faces.  It's worth noting that the finite volume (FV) method is not affected by this error, as mentioned earlier. These errors are also present for level-based methods, with the only distinction being they are addressed for an entire ``grid", rather than on a cell-by-cell basis, as is the case here.

First, we consider the case of solution transfer in the finite-volume method for an element that is split into four hexahedra of different volumes. In this case, a refinement operator that assigns the value contained in the parent to all children  i.e. $\qb_c = \qb_p$ is adequate because of piecewise constant representation of fields in each element. On the other hand, during the coarsening process when the children are merged back into the parent element, assigning the average of children values to the parent i.e. $\qb_p = \sum_{c=1}^{4} \qb_c / 4$ results in a loss of conservation. To preserve conservation, the average should be volume-weighted instead as follows:

\[
\qb_p = \frac{\sum_{c=1}^{4} \qb_c * V_c}{\sum_{c=1}^4 V_c}.
\]

Similarly, for the high-order dGSEM, a conservative coarsening operator can be obtained by integrating over the Legendre-Gauss-Lobatto (LGL) points. Leveraging the tensor-product nature of dGSEM, we first construct projection matrices in 1D and then apply the interpolation in a tensor-product fashion. These matrices are also used for projections of fluxes on non-conformal faces, which is discussed below.

\begin{figure}
\centering
\includegraphics[width=0.4\linewidth]{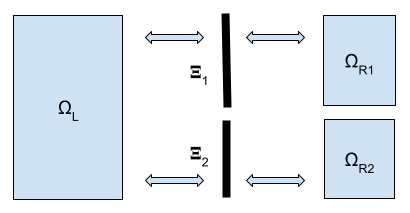}
\caption{Mortar configuration for h-refinement of the right element.}
\label{mortar}
\end{figure}

The grid after refinement/coarsening can become non-conforming, i.e. a case where a single face of an element is shared by two elements on the opposing side. While this does not pose a challenge for the finite volume method, the dG spectral element discretization requires special interpolation techniques at non-conformal faces so as to preserve conservation. We implemented and compared two approaches: a) simple pointwise matching and b) the mortar method for hyperbolic systems introduced by \citet{kopriva1996}. The essence of the second method is that the interpolation should maintain the global conservation property of dGSEM, and at the same time satisfy the ``upwind" nature of hyperbolic systems. The latter is satisfied as a result of a specific mortar configuration shown in Figure \ref{mortar}. The use of two mortars $\Xi_1$ and $\Xi_2$ instead of one allows for discontinuity of solutions at points of intersection of the mortars (see \citep{kopriva1996,giraldoBook,kozdon2018,kopera2014} for details). 

A three-step procedure is involved in the mortar method: 1) the solution is projected from element faces to corresponding mortars, 2) the Riemann problem is solved on the mortars, and 3) the fluxes are projected back onto elements. The projection matrices from subdomain to mortar and vice-versa are constructed by minimizing the $L^2$ norm of the difference between solutions at the face of the element $\Omega_L$ and the mortars. Given local coordinate systems for the mortar $x=[-1,1]$ and element $\xi=[-1,1]$ , the subdomain to mortar projection is formulated minimizing the integral in an $L^2$ sense

\[
\int_{\Xi_i} (\qb^{\Xi_i}(x) - \qb^{\Omega_L}(\xi)) \psi(x) dx = 0\\
\]
which yields the projection

\[
\qb^{\Xi_i} = P^{\Omega_L->\Xi_i} \qb^{\Omega_L} = {M^{\Xi_i}}^{-1} S^{\Omega_L->\Xi_i}\qb^{\Omega_L}
\]
where $M^{\Xi_i}$ is the mass matrix on the mortar, $S^{\Omega_L->\Xi_i}$ is the mixed mass matrix of $\Omega_L$ onto the mortar. Note that the projection matrix from the right elements to the mortars is the identity matrix. We construct the projection matrix in 1-D and apply tensor-products for 2-D and 3-D projections. For the reverse projection from mortars to element face, the following integral is minimized instead

\[
\sum_{i=1}^2 \left( \int_{\Xi_i} (\qb^{\Xi_i}(x) - \qb^{\Omega_L}(\xi)) \psi(\xi) d\xi \right) = 0\\
\]
which yields the projection

\[
\qb^{\Omega_L} = \sum_{i=1}^2 \left( P^{\Xi_i->\Omega_L} \qb^{\Xi_i} \right ) = \sum_{i=1}^2 s_{\Xi_i}{M^{\Omega_L}}^{-1} S^{\Xi_i->\Omega_L}\qb^{\Xi_i}.
\]
where $M^{\Omega_L}$ is the mass matrix of the element face, $S^{\Xi_i->\Omega_L}$ is the mixed mass matrix of the mortar onto $\Omega_L$ and is the transpose of $S^{\Omega_L->\Xi}$, $s_{\Xi_i}$ is a scale factor for the mortar size relative to element size, e.g.  $s_{\Xi_i}$ is typically  $0.5$ for a 2:1 balanced tree. 


These projection matrices are constructed once for the reference element and serve both flux computation at non-conforming faces and solution transfer during refinement and coarsening. While this procedure is adequate for straight-faced affine elements, curvilinear elements utilized for simulation on the sphere require that the projection matrices be constructed for each element separately using the actual mass matrices of the element, instead of the reference element's. 

\subsection{Considerations for spherical geometry}
\label{spherical-geom}
In contrast to limited area models (LAMs), which simulate atmospheric flows in a finite box, general circulation models (GCMs) operate on spherical grids without lateral boundaries. These spherical grids come in various forms, such as latitude-longitude, cubed-sphere, icosahedral, hexagonal, and more. However, not all grid types are equally compatible with adaptive mesh refinement (AMR) techniques.

Among these grid types, cubed-sphere grids are particularly well-suited for both tree-based and level-based AMR libraries. The choice of grid can have a significant impact on the model's performance and ability to capture flow features accurately. For instance, latitude-longitude grids can suffer from a polar singularity issue where grid points are clustered at the poles. Icosahedral grids offer good uniformity and are suitable for dynamic adaptive mesh refinement (see references \citep{giraldo2000,giraldoHesthavenWarburton2002}), however, it is an unstructured grid and unsuitable for many block-structured AMR libraries. It is, however, possible to apply efficient AMR on icosahedral-triangular grid as demonstrated in \citet{zangl2022} for the ICON model.  On the other hand, to support an icosahedral-quadrilateral refinement, each of  the 20 triangles of the icosahedron needs to be subdivided into 3 quadrilaterals, or the 20 triangles combined to form 10 quadrilateral panels. This can potentially introduce several weak singularities at corners of the quadrilaterals compared to just eight for a cubed-sphere grid.
Hexagonal grids, as used in the MPAS  weather model \citep{skamarock2012}, pose unique challenges for dynamic mesh refinement. Unlike quad- or hex-based grids, they lack a direct and straightforward approach to dynamically refine the mesh. While it is possible to generate a static variable-resolution hexagonal grid based on a specified cell density function, dynamically refining the grid to capture moving features may not be feasible or could require regridding the entire domain, even for localized flow features. Moreover, the adaptive mesh generator must generate a conformal grid, thereby adding to the overall complexity of the process.

Both LAMs and GCMs take advantage of the stratified nature of the atmosphere by treating vertical (convective) processes differently from horizontal (advective) processes. In hydrostatic models, vertical processes are essentially ignored. When it comes to AMR, the considerable difference in spatial scales, nearly three orders of magnitude, between the vertical and horizontal directions highlights that AMR is typically more critical for the horizontal aspect than the vertical. Therefore, AMR often involves anisotropic refinement, where the vertical direction is either excluded from refinement or is refined to a lesser extent than the horizontal. Libraries like p4est \citep{p4est} and AMReX\citep{amrex2021} have built-in support for this type of anisotropic refinement.

Spherical geometry introduces additional complexity, particularly regarding interpolation and restriction operators. To ensure conservation of important physical quantities, these operators must be adapted for use in curvilinear coordinates on the sphere 
(see, e.g., \citep{vanderholst2007} and \citep{jablonowski2006}). Moreover, the projection matrices described in the previous section should be computed individually for each element. However, we do not use this approach and propose a simple and efficient alternative described below.

Instead of computing projection matrices individually for each element on a curved surface, we propose a novel approach that is simpler to implement and efficient because only one set of projection matrices are needed.  We utilize the projection matrices computed using the reference element for all elements and then address the solution transfer errors introduced by this assumption. The specific details of solution transfer are as follows. a) For coarsening, the contribution of each child element is weighed by its volume, similar to the finite-volume method, ensuring conservation. b) For refinement, the total integrated field value after refinement is adjusted to be equal to that before refinement. A single scale factor is computed and applied uniformly to all fields in the child elements. This approach ensures conservation, however, we believe a solution transfer strategy that utilizes projection matrices per element should yield a more accurate result that takes less time to adjust. 

Additionally, the areas and volumes of the curvilinear elements are computed using the spherical triangle area formula, ensuring that  the total surface area before and after refinement remains unchanged. We use the spherical excess

\[
A = r_e^2(\theta_1+\theta_2+\theta_2 - \pi)
\]
where $(\theta_1,\theta_2,\theta_3)$ are the internal angles, and $r_e$ is the radius of the earth. Alongside metrics for mass and energy loss, we include a ``volume loss" metric to monitor potential loss in volume due to refinement. If straight-faced affine elements were used for the sphere, the total surface area would increase after refinement, analogous to the concept of approximating $\pi$ using an inscribed polygon, resulting in a net area/volume gain. 

\subsection{Scalability and performance portability}
NWP  involves computations that demand immense processing power. Exascale supercomputers, capable of performing one quintillion calculations per second, hold the potential to significantly enhance the resolution and accuracy of NWP. Nevertheless, the development of a massively scalable NWP model that effectively utilizes current and emerging HPC technologies in a performance-portable manner poses a considerable challenge. This challenge is further exacerbated by AMR because a) managing several levels of grids in parallel require intricate algorithms and communication strategies and b) the mesh manipulation code may not be suitable to GPU acceleration. The approach adopted by the AMR codes in this study is briefly outlined below.

ERF reaps immediate benefit from employing the AMReX framework, intentionally designed with scalability and performance portability in mind. Because level-based methods advance the solution at different levels with potentially different timesteps, the parallel communication strategy is notably more complex than that of tree-based methods. The solution procedure involves several steps, including a refinement operator for coarse-to-fine level communication, the execution of the solver on each mesh level, a coarsening operator for fine-to-coarse level communication, and a flux-correction step to ensure conservation of quantities. Many level-based AMR frameworks are implemented using this communication protocol, which can pose challenges for finite-element method implementations \citep{kirby2018}. 
The grids can have a complex layout making inter-level and inter-processor exchange non-trivial. For example, AMReX employs a hash-based algorithm for performing intersection between nested grids \citep{amrex2021}. AMReX adopts hierarchical parallelism with MPI for coarse-grained parallelism, and utilizes OpenMP/CUDA/HIP/SYCL to leverage many-core compute architectures inside nodes. A software layer, akin to that of the Kokkos library~\citep{edwards2013}, conceals implementation details and provide performance-portable code across various devices.

AMReX has been employed in PeleC \citep{frahan2023}, a solver designed for compressible reacting flows, to prepare it for exascale computing. PeleC has showcased scalability on the entire Summit supercomputer, equipped with 27,648 V100 GPUs, achieving a weak scaling parallel efficiency of 65\%. It is noteworthy that the initial PeleC code utilized OpenACC for GPU acceleration, demonstrating similar level of scalability on the supercomputer. Nevertheless, AMReX introduces an additional advantage by enabling compatibility with Intel and AMD GPUs through their native language implementations.

In contrast, NebulaSEM does not scale well on many nodes primarily  because the regridding process is carried out on a single node. However, scalable tree-based AMR libraries, such as \citet{p4est}, do exist.  Unlike level-based methods, tree-based AMR utilize non-overlapping patches, with communication occurring solely at mesh-partition boundaries between two neighboring elements. Although NebulaSEM technically uses a cell-based AMR, its utilization of high-order elements aids in achieving performance portability to many-core architectures. Acceleration is attained using OpenMP and OpenACC directives on vector computations, with managed memory facilitating automatic data transfers. Nevertheless, the mesh adaptation code is still executed on the CPU due to the complexity of managing an unstructured grid on the GPU.

\section{Results}
\label{sec-results}
Here, we present the results of a comparison between tree-based and level-based AMR methods applied to a set of Numerical Weather Prediction (NWP) benchmark problems. While both codes support 3D AMR, for the purpose of a straightforward comparison, the test cases we analyze involve only 2D AMR. Moreover, to assess the accuracy of AMR simulations, we focus on problems with known exact solutions, such as the isentropic vortex problem, and those featuring a reversing flow field where the solution at the end of the simulation period is precisely known. In instances where obtaining exact solutions is not feasible, we conduct two uniform grid simulations that encapsulate the AMR simulation, representing the coarsest and finest resolutions used in the AMR simulation.

\subsection{Isentropic vortex problem}
The first test case we study is  the classic isentropic vortex problem (see f.i \citet{spiegel2015}) -- a unique type of flow where the prognostic variables are linked in a constant entropy manner.  This test case serves as a litmus test for assessing the accuracy of numerical methods and the impact of grid refinement, as it offers a rare advantage:  an exact analytical solution to the Euler equations. The isentropic vortex is characterized by a circular region of rotating fluid, with the vortex's boundary representing a zone of sharp velocity gradients. It stands as a pivotal challenge for examining a numerical method's capacity to maintain and preserve flow features, particularly vortices, over extended periods. An essential aspect of this examination is the assessment of numerical dissipation, a critical consideration in the realm of NWP models.

The free stream conditions for the problem are defined as follows

\[
\rho = 1, u = U_\infty, v = V_\infty, \theta = \theta_\infty, p= p_\infty.
\]
The initial conditions are obtained by introducing perturbations to the mean velocity and potential temperature, which are expressed as:

\[
(u',v') = \frac{\beta}{2\pi} \exp\left(\frac{1-r^2}{2}\right) (-y+y_c,x-x_c)
\]
\[
\theta' = \frac{(\gamma - 1)\beta^2}{8\gamma\pi^2}\exp{(1-r^2)}
\]
where $\gamma$ is the gas constant of air as previously defined, and
\[
r= \sqrt{{(x-x_c)}^2+{(y-y_c)}^2}.
\]
The pressure and potential temperature are linked with the isentropic condition
\[
p=\rho^\gamma = (\theta + \theta')^{\frac{\gamma}{\gamma - 1}}.
\]

For the simulation of the isentropic vortex problem, a computational domain of dimensions  $[-6m,6m]$ x $[-6m,6m]$ is considered. Other parameter values are set as follows: $(x_c,y_c)=(0,0)$, $\beta=5$, $U_\infty = 1$ m/s,  $V_\infty=1$ m/s , $p_\infty = 1$ Pa and $\theta_\infty$=1 K. The domain is discretized into $32\times32$  elements with polynomial order of $N=3$ for the tree-based method, and into $96\times96$ elements for the level-based method. 
For both AMR methods, only one level of refinement is allowed. The simulation runs for 12 seconds with a fixed time step of $\Delta t = 0.00025$ seconds using the explicit forward Euler method for the tree-based AMR. For the level-based method, the coarse grid time step size is $\Delta t = 0.00015$ s.  We run the level-based method with grid efficiencies of 0.7 and 0.9.

\begin{figure}
\begin{subfigure}{\textwidth}
\includegraphics[width=0.24\linewidth]{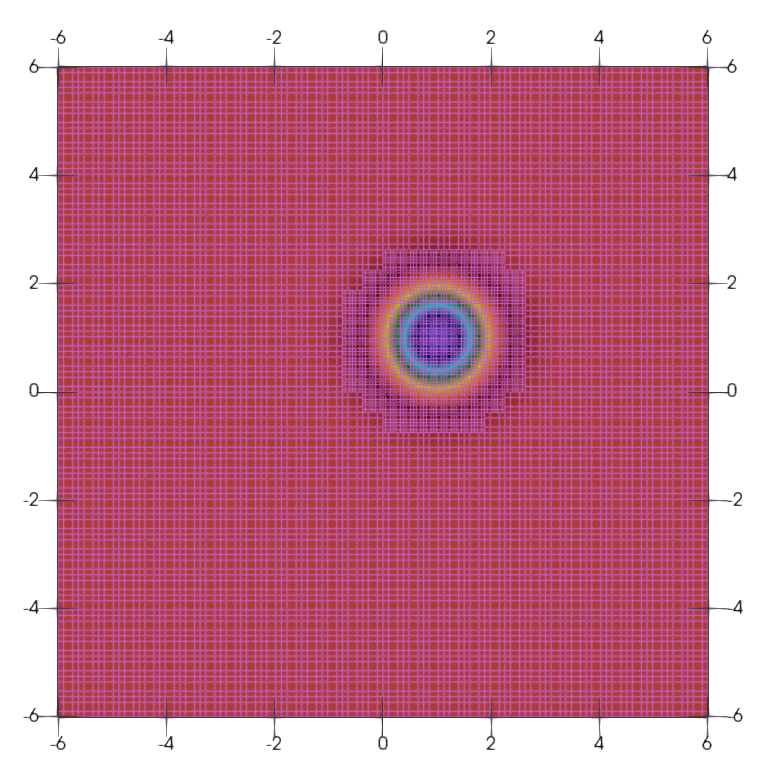}
\includegraphics[width=0.24\linewidth]{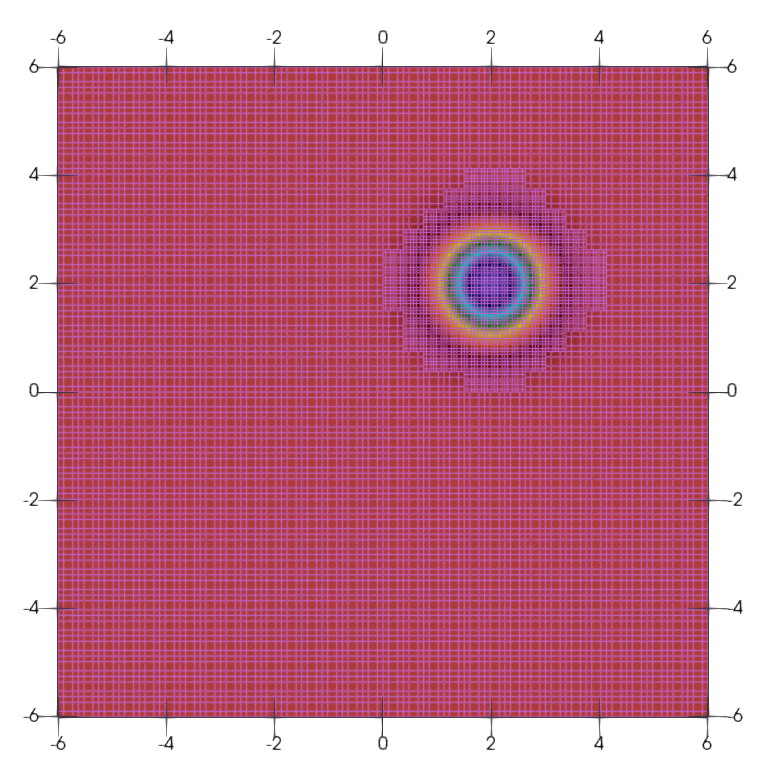}
\includegraphics[width=0.24\linewidth]{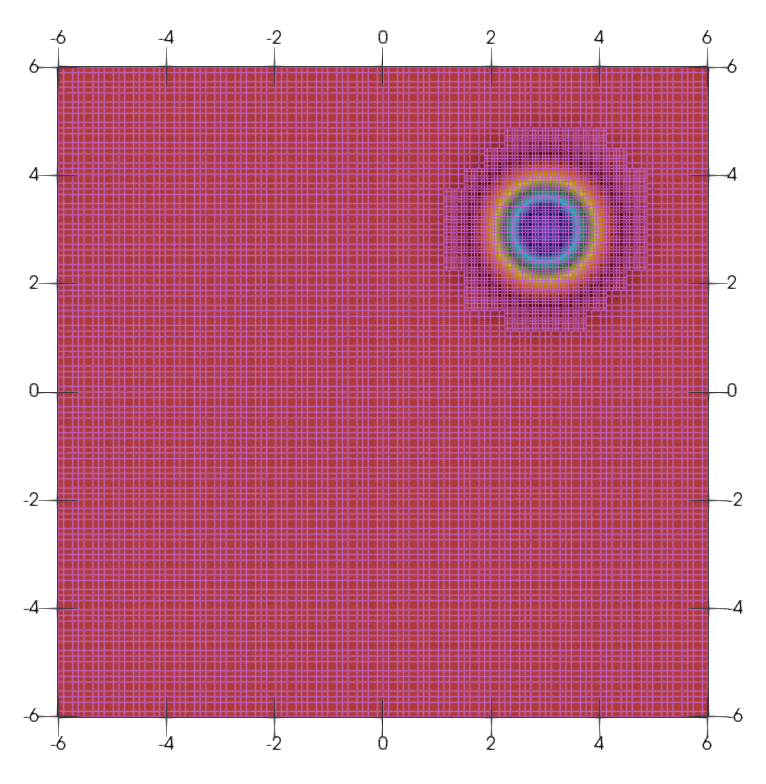}
\includegraphics[width=0.24\linewidth]{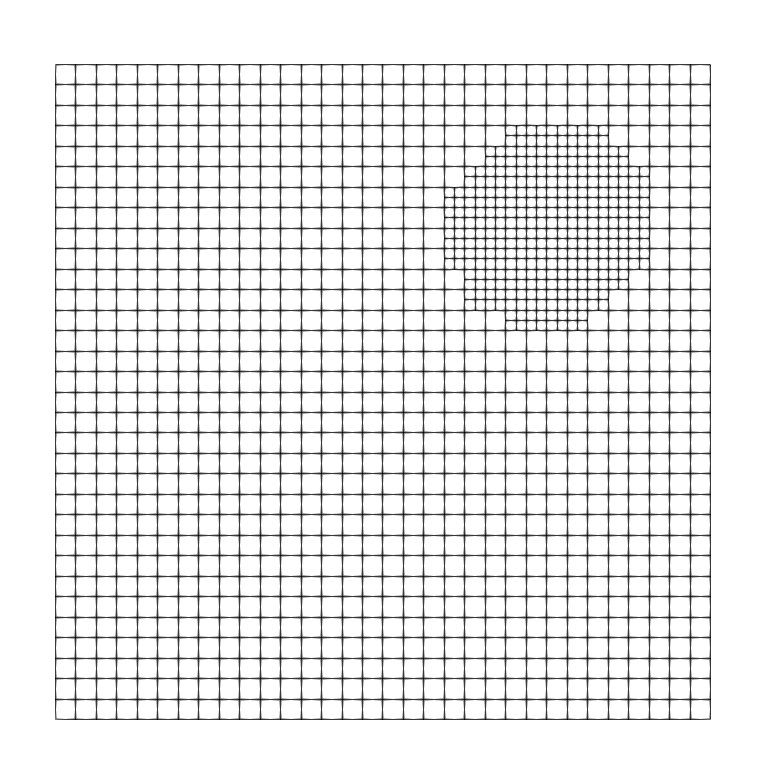}
\caption{Tree-based AMR}
\end{subfigure}
\begin{subfigure}{\textwidth}
\includegraphics[width=0.24\linewidth]{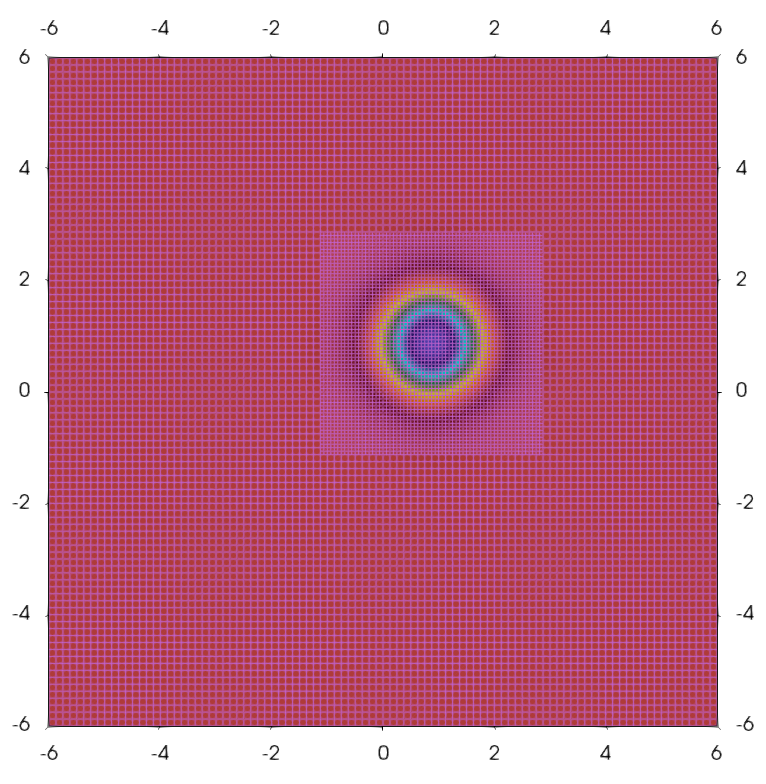}
\includegraphics[width=0.24\linewidth]{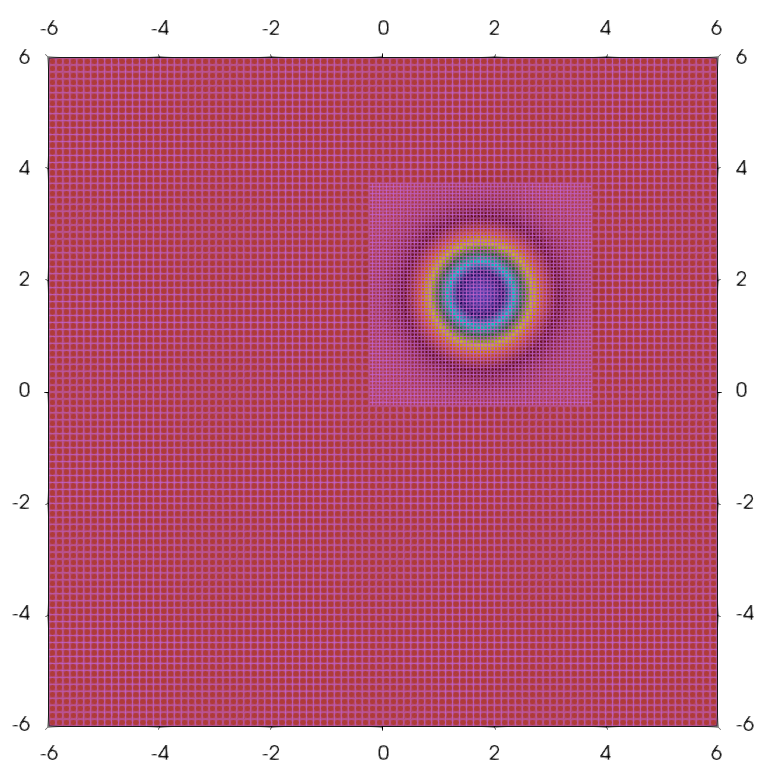}
\includegraphics[width=0.24\linewidth]{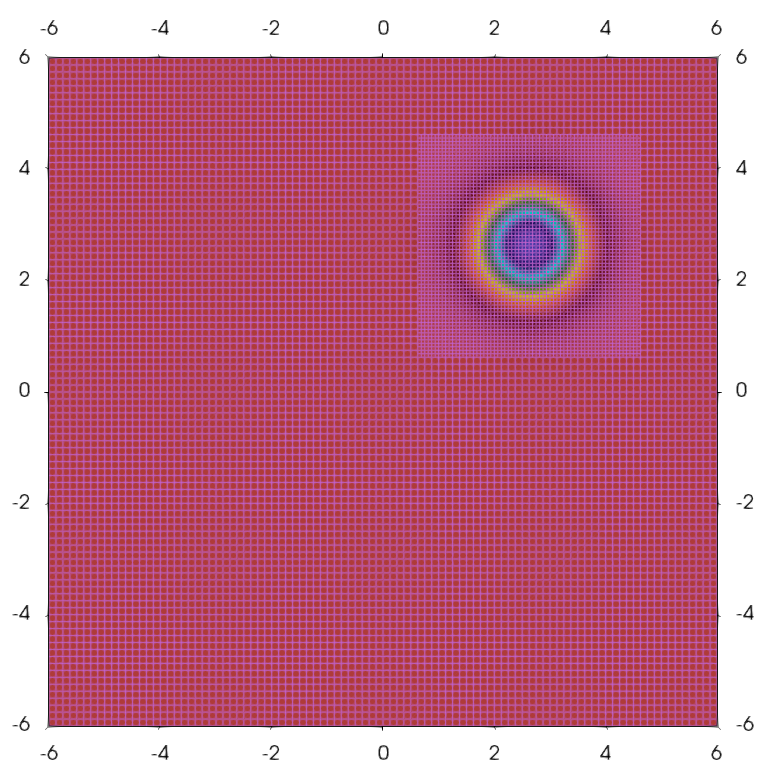}
\includegraphics[width=0.24\linewidth]{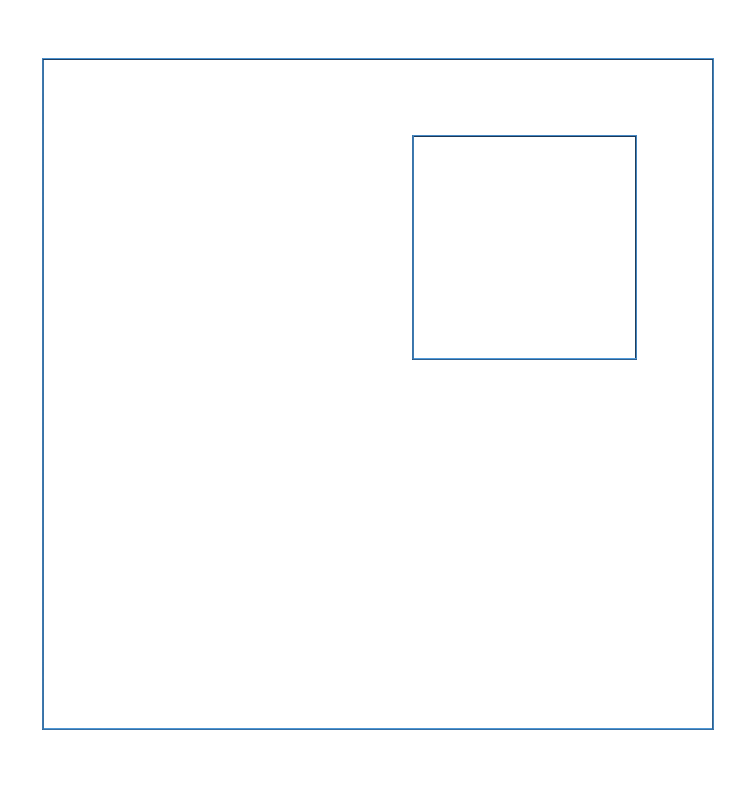}
\caption{Level-based AMR with grid efficiency of 0.7}
\end{subfigure}
\begin{subfigure}{\textwidth}
\includegraphics[width=0.24\linewidth]{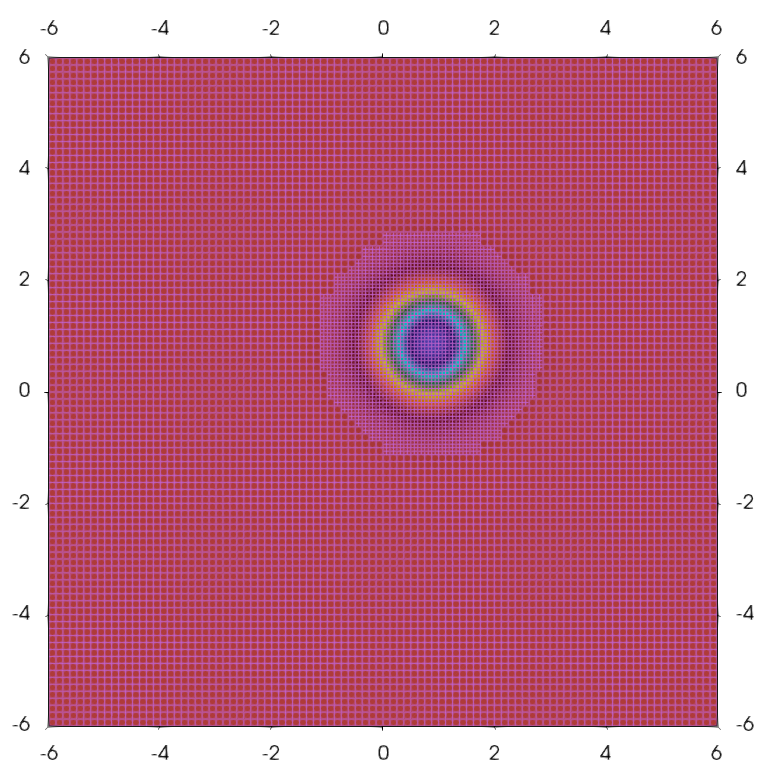}
\includegraphics[width=0.24\linewidth]{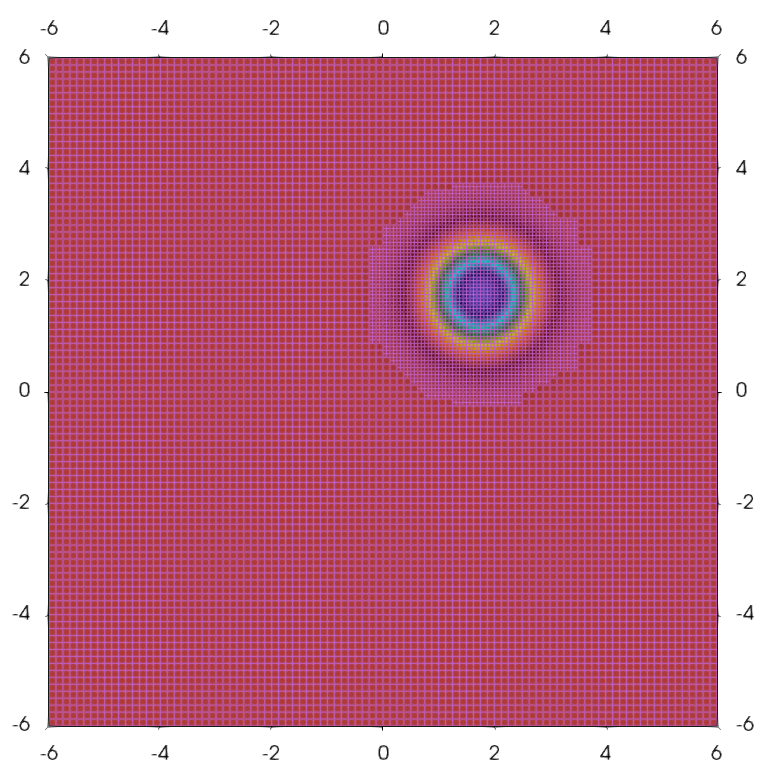}
\includegraphics[width=0.24\linewidth]{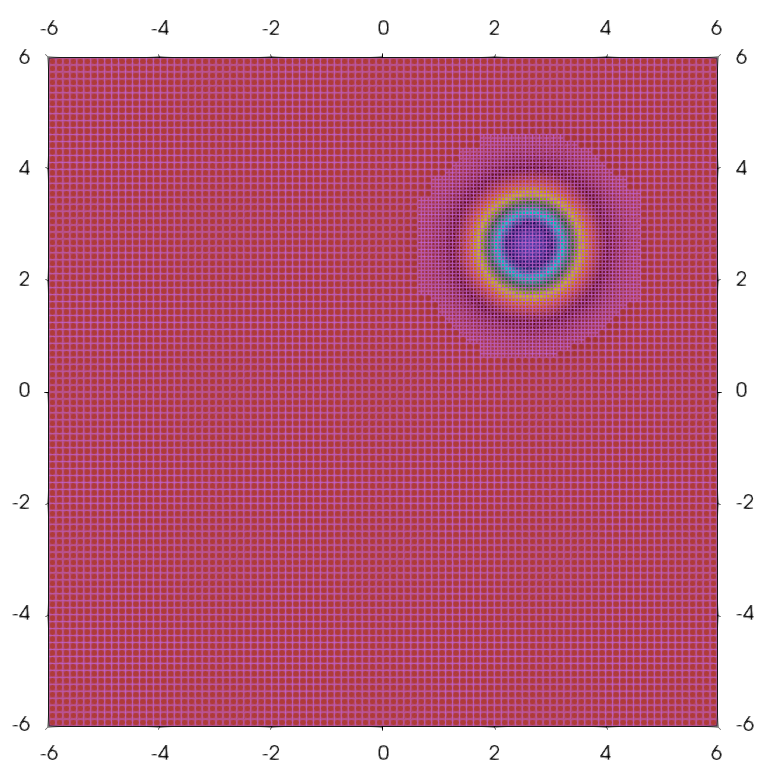}
\includegraphics[width=0.24\linewidth]{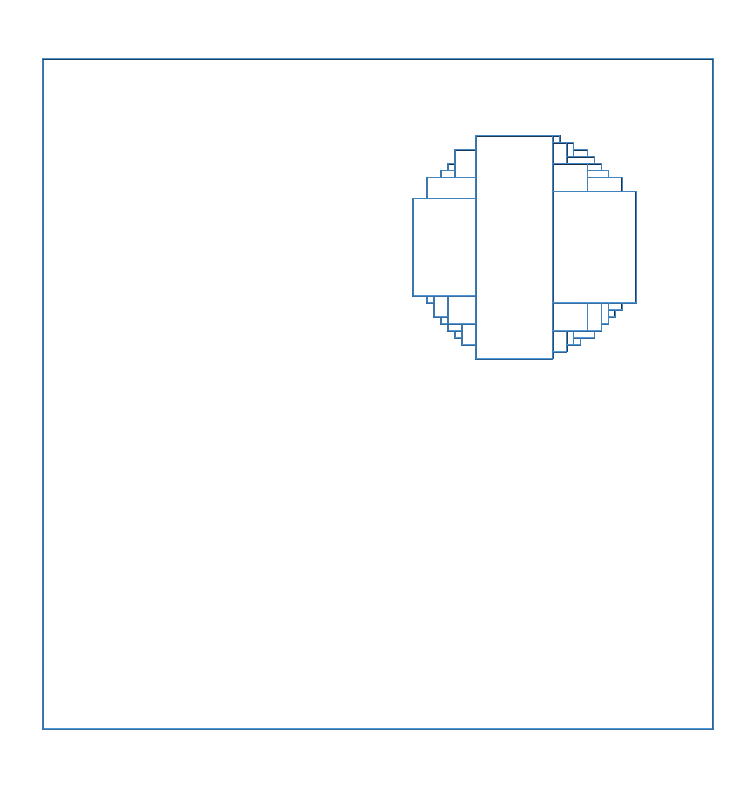}
\caption{Level-based AMR with grid efficiency of 0.9}
\end{subfigure}
\caption{Evolution of the isentropic vortex for 3 seconds (snapshots at 1, 2 and 3 seconds). The potential temperature is dipicted with three different configurations: a) tree-based AMR with a buffer zone of 2 cells b) level-based AMR with grid efficiency of 0.7. Notably, only one rectangular patch is created in this case. c) level-based AMR with grid efficiency of 0.9, resulting in the creation of multiple patches.}
\label{isentropic-vortex}
\end{figure}

Throughout the simulation, the isentropic vortex maintains its shape while moving along the diagonal at a constant velocity, as illustrated in Fig.~\ref{isentropic-vortex}. Both tree-based and level-based AMR approaches demonstrate an effective tracking of the vortex's evolution. Notably, when using the level-based AMR with a grid efficiency of 0.7, it requires only one rectangular subgrid, akin to a moving nest simulation. However, increasing the grid efficiency to 0.9 results in AMReX generating up to 20 rectangular subgrids of variable sizes on the same level to better fit the vortex's shape. Hence, it's imperative to carefully select the AMR parameters to avoid the creation of an excessively large number of grids. As the grid efficiency approaches 1, the level-based AMR effectively becomes a cell-based AMR. We should note that the AMReX infrastructure is not designed to operate efficiently in this limit.

\begin{figure}
\centering
\includegraphics[width=0.32\linewidth]{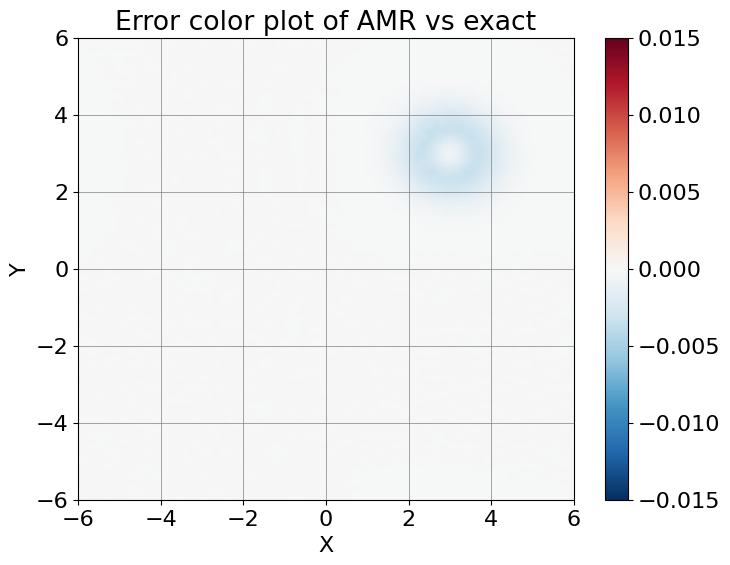}
\includegraphics[width=0.32\linewidth]{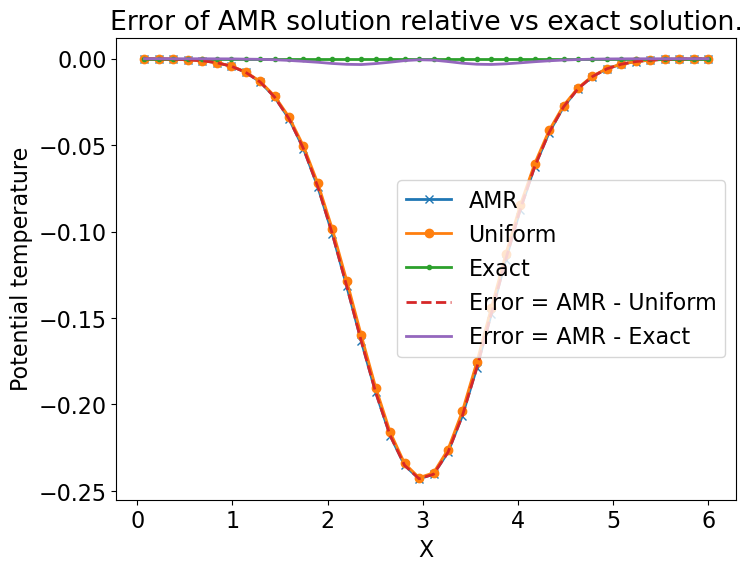}
\includegraphics[width=0.32\linewidth]{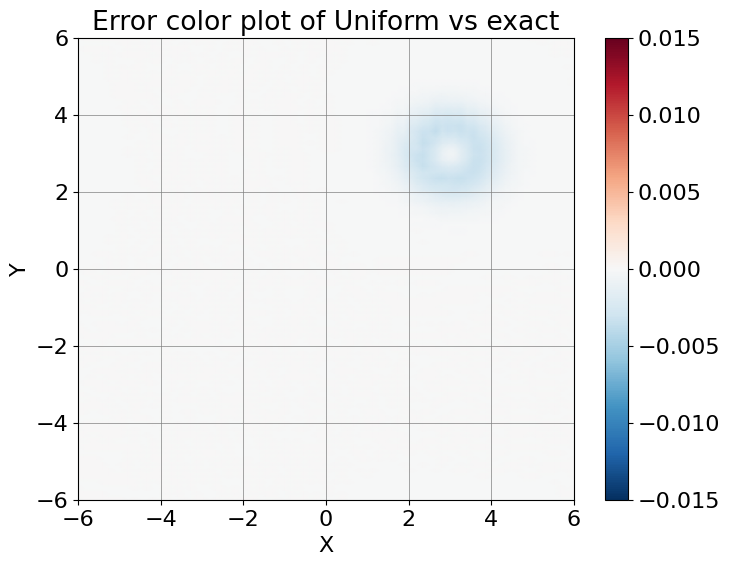}
\caption{Potential temperature error for the isentropic vortex problem at t = 3 sec using the tree-based AMR method. Left figure shows that the AMR simulation error is mostly zero except within the bubble where small differences can be observed. The middle figure depicts profile of potential temperature at y = 3m. The right figure shows the error plot for the uniform grid simulation of $32\times32$ elements.}
\label{isentropic-vortex-error}
\end{figure}

To assess the accuracy of the tree-based AMR simulation, we compute the $L^2$ norm of the error with respect to the exact solution, which is a simple translation of the vortex at $1$ m/s.

\[
L_{\Omega}^2 = \frac{1}{ N_{dof}} \sqrt{\sum_{i=1}^{N_{dof}}|\qb - \qb_{exact}|^2}
\]
The computed $L^2$ norm of the error after 3 seconds is approximately $8.75 \times 10^{-4}$.  The plot depicting the deviation of the solution from the exact solution is presented in Figure \ref{isentropic-vortex-error}. It illustrates that, outside of the bubble (localized feature), the error is predominantly zero. 

\subsection{Linear advection }
Atmospheric dynamics are primarily dominated by advective processes. Therefore, AMR which enhances the solution of the advection problem plays an important role in the accuracy of NWP models \citep{jablonowski2006}. We investigate a standard test case, swirling deformation flow introduced in \citet{leveque1996}, with a flow field that induces deformation as opposed to widely used solid-body advection tests. A time-reversing flow field is used in such a way that the exact solution after one full period (T) becomes the same as the initial condition.  At half the period (T/2), the deformation should peak forming a thin filament.

\begin{figure}
\centering
\begin{subfigure}{0.32\textwidth}
\centering
\includegraphics[width=0.8\linewidth]{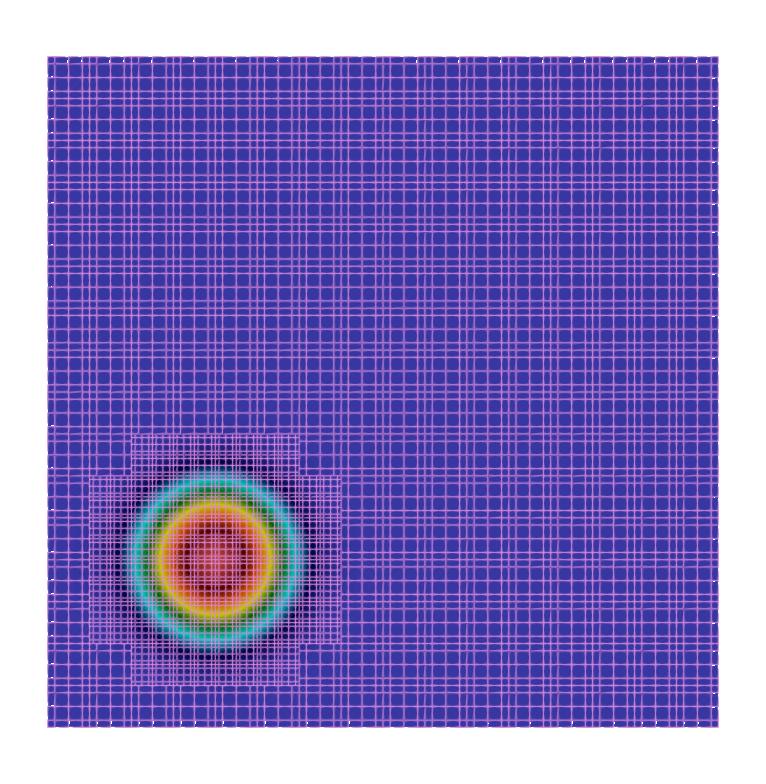}
\caption{t=0}
\end{subfigure}
\begin{subfigure}{0.32\textwidth}
\centering
\includegraphics[width=0.8\linewidth]{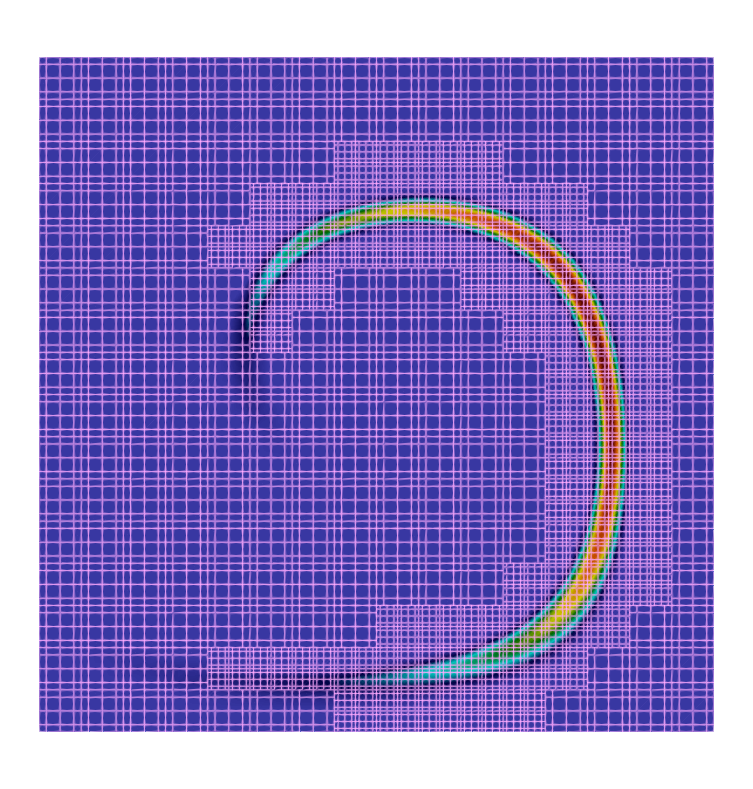}
\caption{t=T/4}
\end{subfigure}
\begin{subfigure}{0.32\textwidth}
\centering
\includegraphics[width=0.8\linewidth]{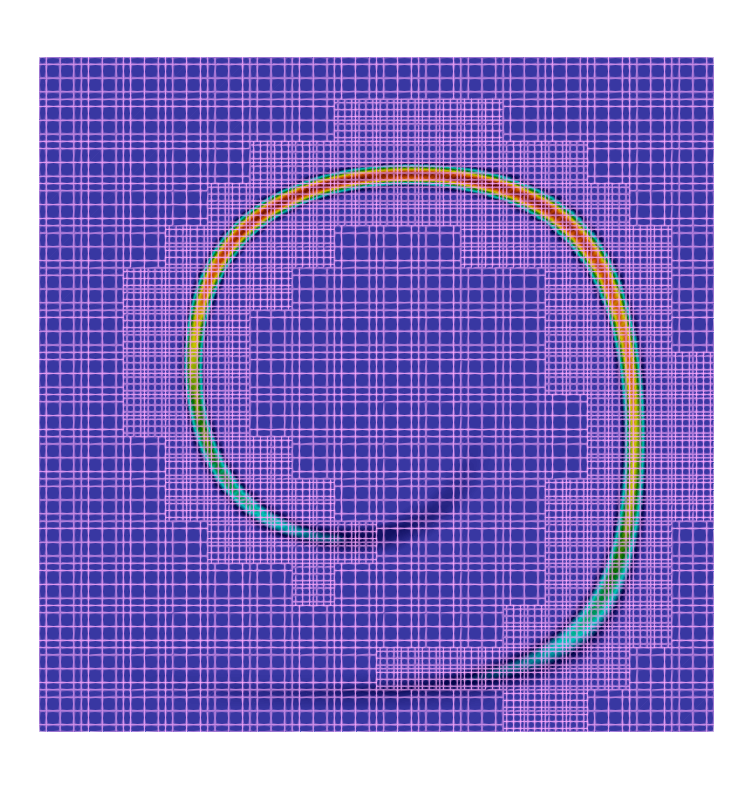}
\caption{t=T/2}
\end{subfigure}
\begin{subfigure}{0.32\textwidth}
\centering
\includegraphics[width=0.8\linewidth]{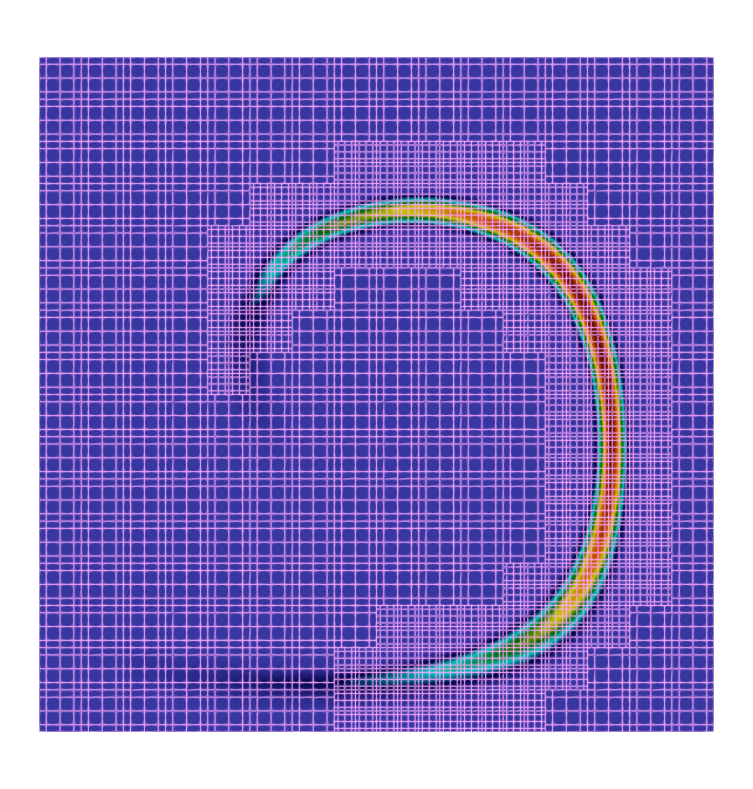}
\caption{t=3T/4}
\end{subfigure}
\begin{subfigure}{0.32\textwidth}
\centering
\includegraphics[width=\linewidth]{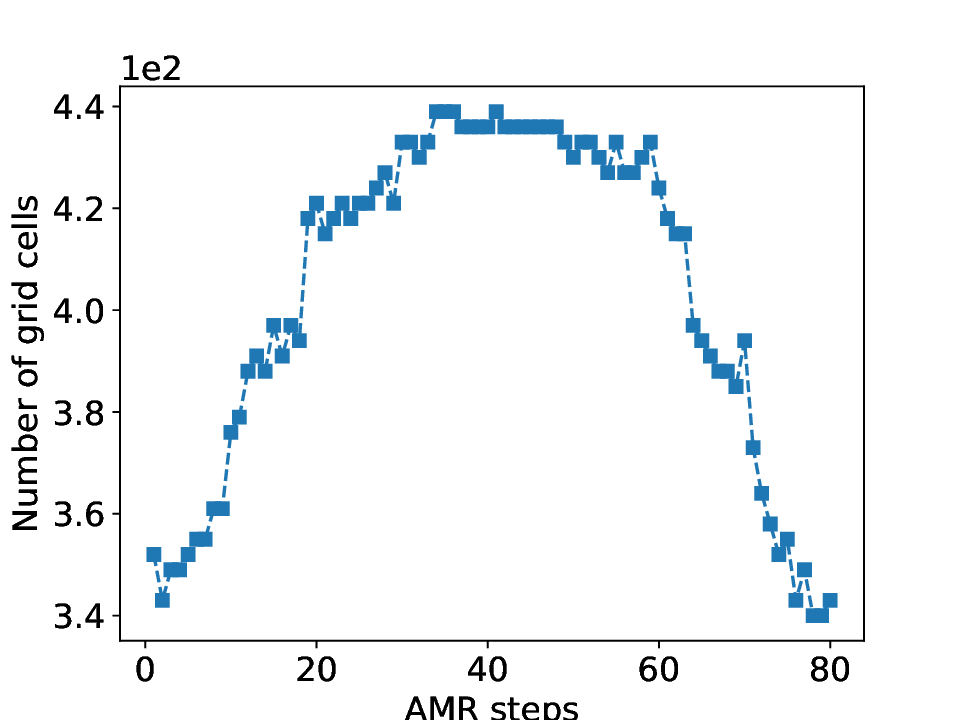}
\caption{Number of grid cells}
\end{subfigure}
\begin{subfigure}{0.32\textwidth}
\centering
\includegraphics[width=0.8\linewidth]{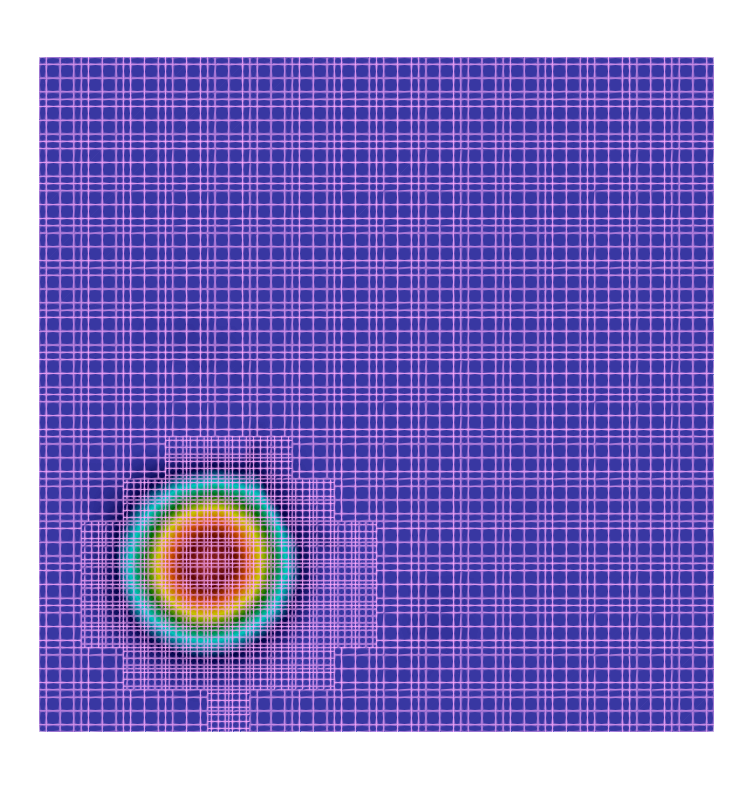}
\caption{t=T}
\end{subfigure}
\begin{subfigure}{0.32\textwidth}
\centering
\scalebox{1}[-1]{\rotatebox{-90}{\includegraphics[width=0.8\linewidth]{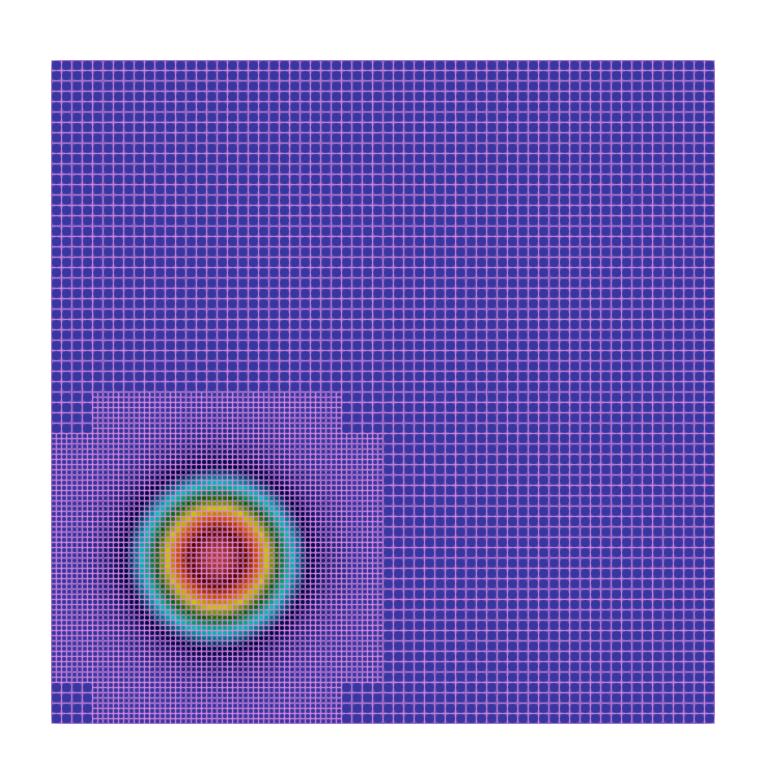}}}
\caption{t=0}
\end{subfigure}
\begin{subfigure}{0.32\textwidth}
\centering
\scalebox{1}[-1]{\rotatebox{-90}{\includegraphics[width=0.8\linewidth]{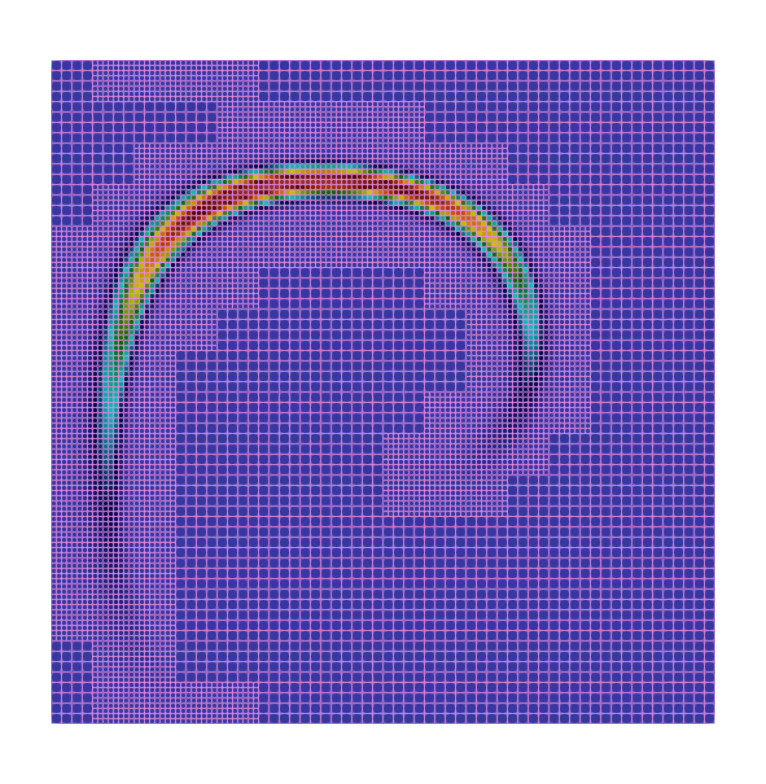}}}
\caption{t=T/4}
\end{subfigure}
\begin{subfigure}{0.32\textwidth}
\centering
\scalebox{1}[-1]{\rotatebox{-90}{\includegraphics[width=0.8\linewidth]{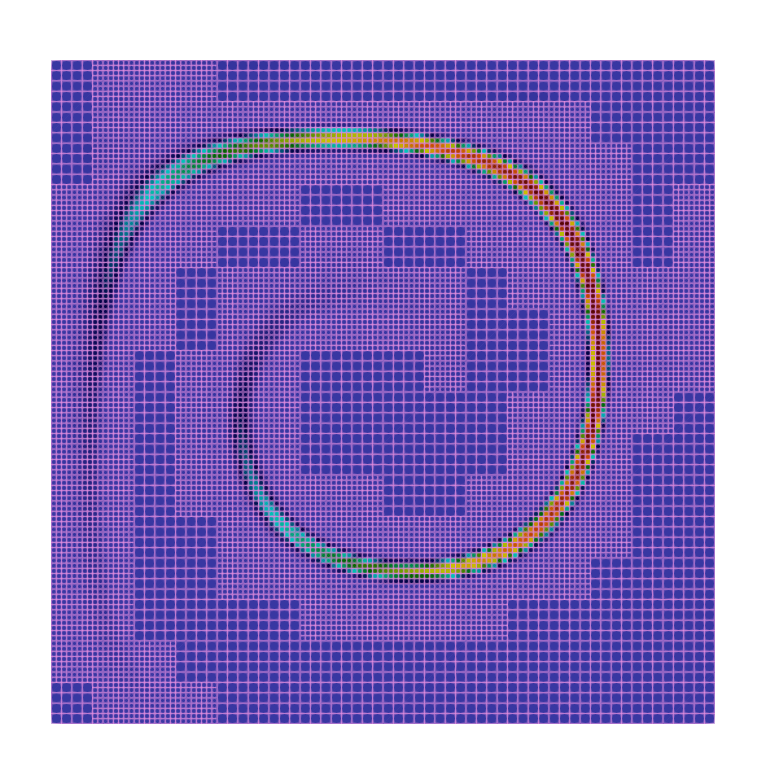}}}
\caption{t=T/2}
\end{subfigure}
\begin{subfigure}{0.32\textwidth}
\centering
\scalebox{1}[-1]{\rotatebox{-90}{\includegraphics[width=0.8\linewidth]{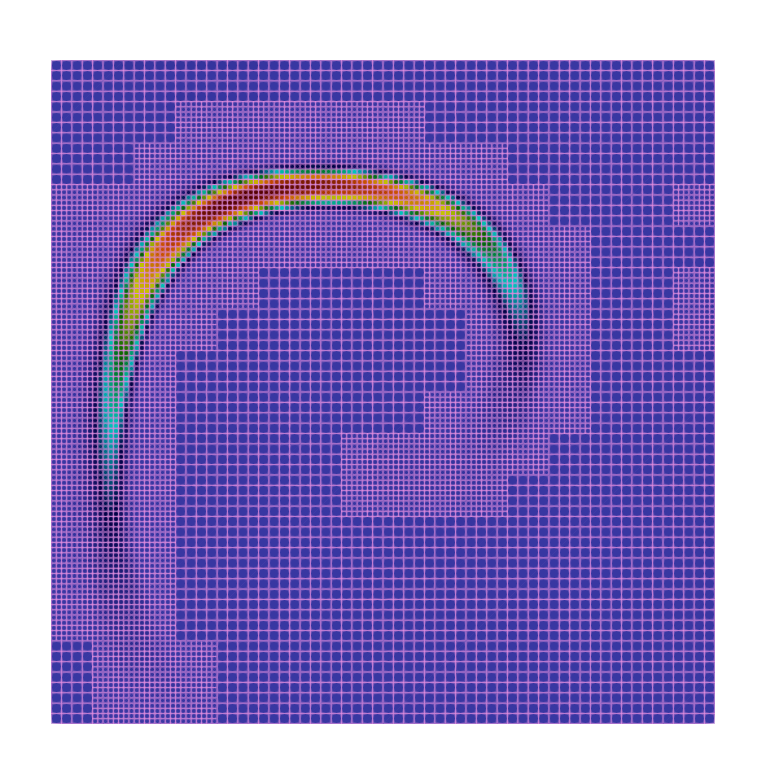}}}
\caption{t=3T/4}
\end{subfigure}
\begin{subfigure}{0.32\textwidth}
\centering
\includegraphics[width=\linewidth]{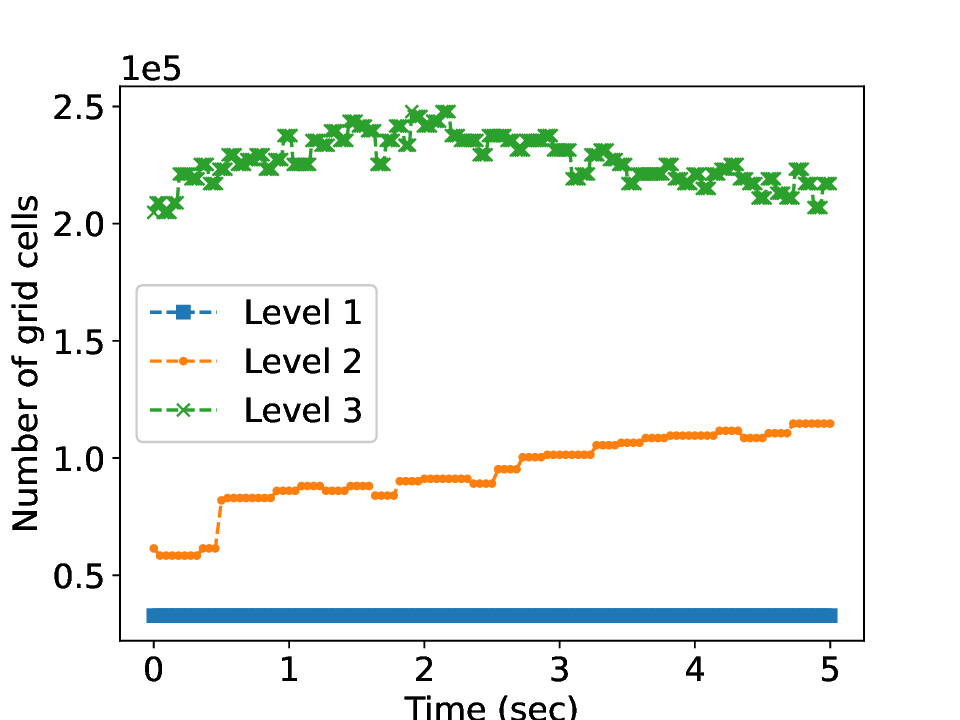}
\caption{Number of grid cells}
\end{subfigure}
\begin{subfigure}{0.32\textwidth}
\centering
\scalebox{1}[-1]{\rotatebox{-90}{\includegraphics[width=0.8\linewidth]{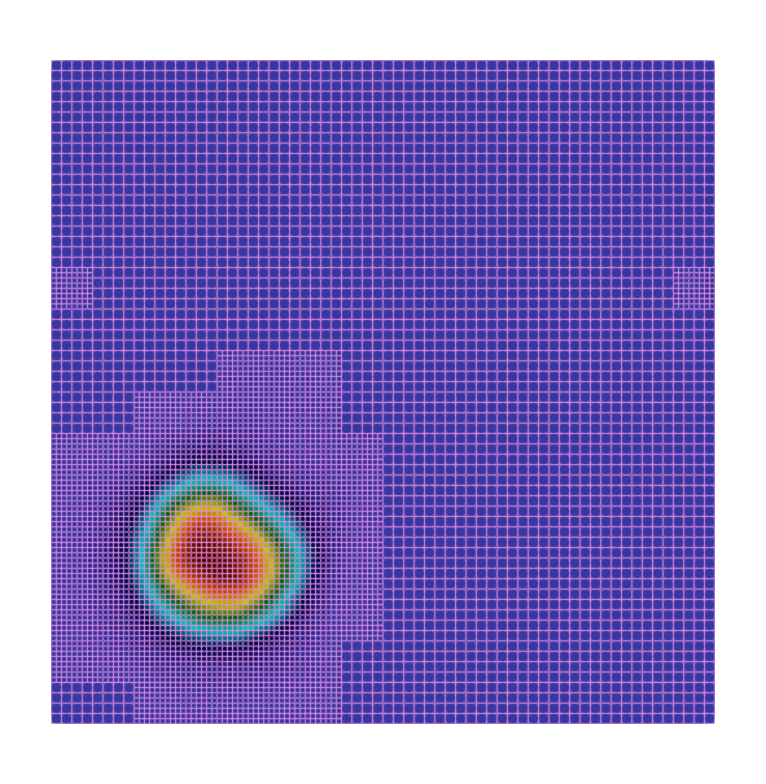}}}
\caption{t=T}
\end{subfigure}
\caption{Evolution of an initially cosine-bell shaped tracer of the swirling flow test case introduce in \citet{leveque1996} for one full period using the two AMR approaches. Figures (a)-(f) show results for the tree-based AMR with a buffer zone of 2 cells and one level of refinement. Figures (g)-(l) show results for the level-based AMR with grid efficiency of 0.95 and three levels of refinement. In both cases, a thin filament forms at half the period, at which the wind flow reverses and the filament deforms back to a shape nearly identical to the initial state.  Both approaches adeptly follow the shape of the filament. }
\label{advection-leveque}
\end{figure}

This test case involves the advection of a cosine-bell distribution of a tracer located in the bottom-left of the domain, as depicted in Fig. \ref{advection-leveque}. The initial shape of the tracer is defined by the following equation for a given cutoff radius $r_c$ and maximum concentration $h_{max}$

\begin{equation}
h' = 
\begin{cases}
        0  &  \text{for } r > r_c \\
        \frac{h_{max}}{2} (1 + \cos(\frac{\pi r}{r_c})) &  \text{for } r \le r_c.
\end{cases}
\label{cosine-bell}
\end{equation}
The driving wind field is given by

\[
\begin{aligned}
u = sin^2(\pi x) sin(2 \pi y) g(t)\\
v = -sin^2(\pi y) sin(2 \pi x) g(t)
\end{aligned}
\]
where $g(t)$ is used to introduce time dependence in the flow field for $0 \le t \le T$

\begin{equation}
g(t) = cos(\pi t / T).
\label{cosine-time}
\end{equation}

For the simulation, we use a computational domain of dimensions $[0,1]$ x $[0,1]$, with the initial tracer bubble centered at $[0.25,0.25]$ with radius $0.25$. We employ both the tree-based AMR and level-based AMR methods initially subdividing the domain into a grid of $16\times16$ elements with polynomial order $N = 4$, and into $64\times 64$ elements for the latter. The simulation is run for a time period $T=5$ , and regridding is carried out 80 times and 250 times for the tree-based and level-based methods, respectively. 

As depicted in Fig. \ref{advection-leveque}, the tracer bubble undergoes significant deformation due to the background wind field, reaching its  peak distortion half way through the simulation before returning to its original position at the end of one full period. Both methods effectively track the evolving shape of the tracer bubble, with the number of grid cells peaking at $T/2$ for the tree-base method. For the level-based method, the number of cells in level 0 (the background) is the same throughout the simulation, unlike the other levels. We used a grid efficiency of 0.95 just for comparisons sake, but the default value of 0.7 is more efficient with only a few boxes created to track the filament.

\begin{figure}
\centering
\begin{subfigure}{0.32\textwidth}
\centering
\includegraphics[width=\linewidth]{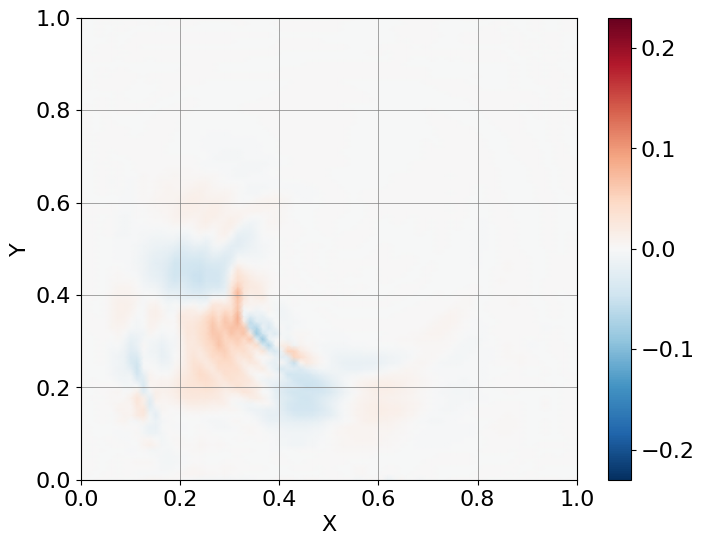}
\caption{AMR initial $16\times16$ elements}
\end{subfigure}
\begin{subfigure}{0.32\textwidth}
\centering
\includegraphics[width=\linewidth]{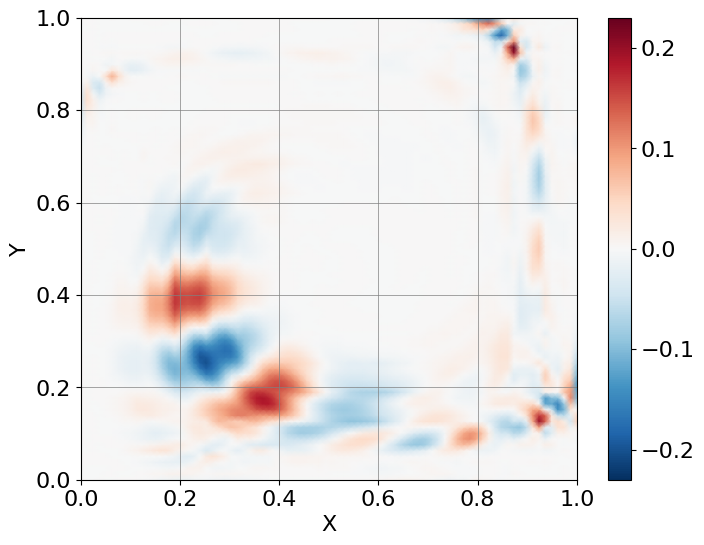}
\caption{Conformal $16\times16$ elements}
\end{subfigure}
\begin{subfigure}{0.32\textwidth}
\centering
\includegraphics[width=\linewidth]{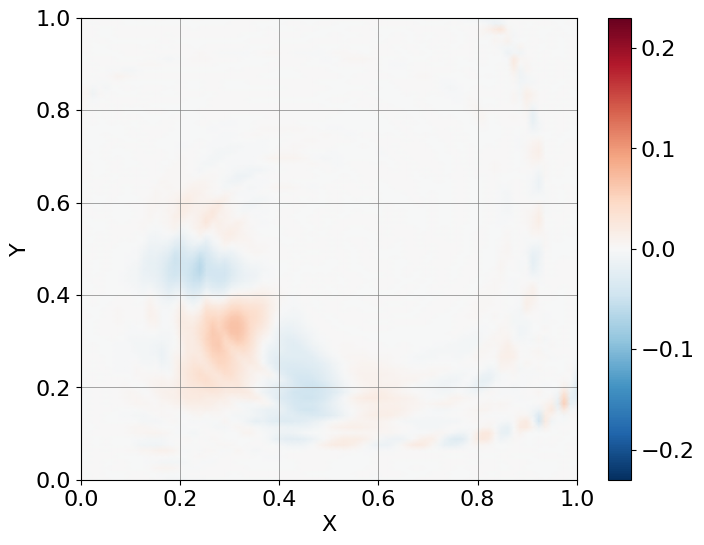}
\caption{Conformal $32\times32$ elements}
\end{subfigure}
\caption{Comparison of solution errors of the test case by \citet{leveque1996} after a full period (t=T), at which point the tracer shape should match its initial configuration. The solution error with the tree-based AMR (a) is considerably lower than the solution obtained with a conformal grid (b) using the initial coarse resolution applied in the AMR simulation. The AMR result is comparable to the fine uniform grid resolution that uses $32\times32$ elements instead.}
\label{leveque-error}
\end{figure}

The comparison of the tracer distribution error at $t=T$ relative to $t=0$ is depicted in Fig. \ref{leveque-error}. Adaptive mesh refinement (AMR) has significantly enhanced the solution compared to a conformal grid, utilizing the initial coarse resolution of $16\times16$ elements in the AMR simulation. Achieving a comparable level of accuracy as the AMR simulation would require doubling the resolution of the uniform grid to $32\times32$ elements. Additionally, the time required for the AMR simulation, the coarse uniform grid, and the fine uniform grid simulations are 11 seconds, 18 seconds, and 51 seconds, respectively. This underscores the advantage of AMR in improving accuracy to the level of the fine uniform grid simulation while adding minimal running costs.

\begin{figure}
\centering
\begin{subfigure}{0.49\textwidth}
\includegraphics[width=\linewidth]{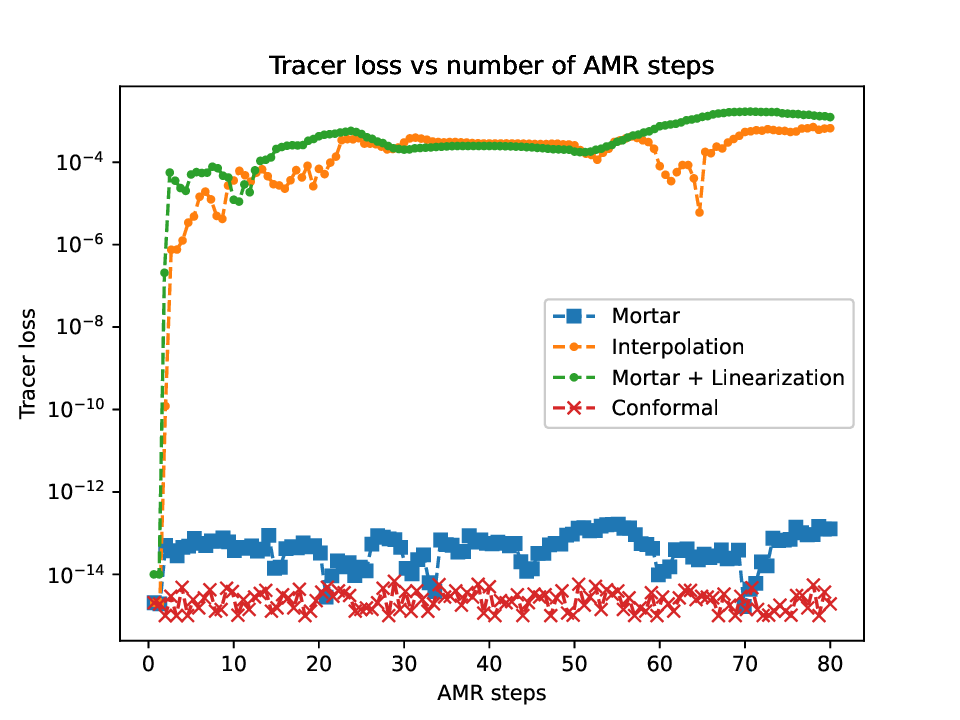}
\end{subfigure}
\begin{subfigure}{0.49\textwidth}
\includegraphics[width=\linewidth]{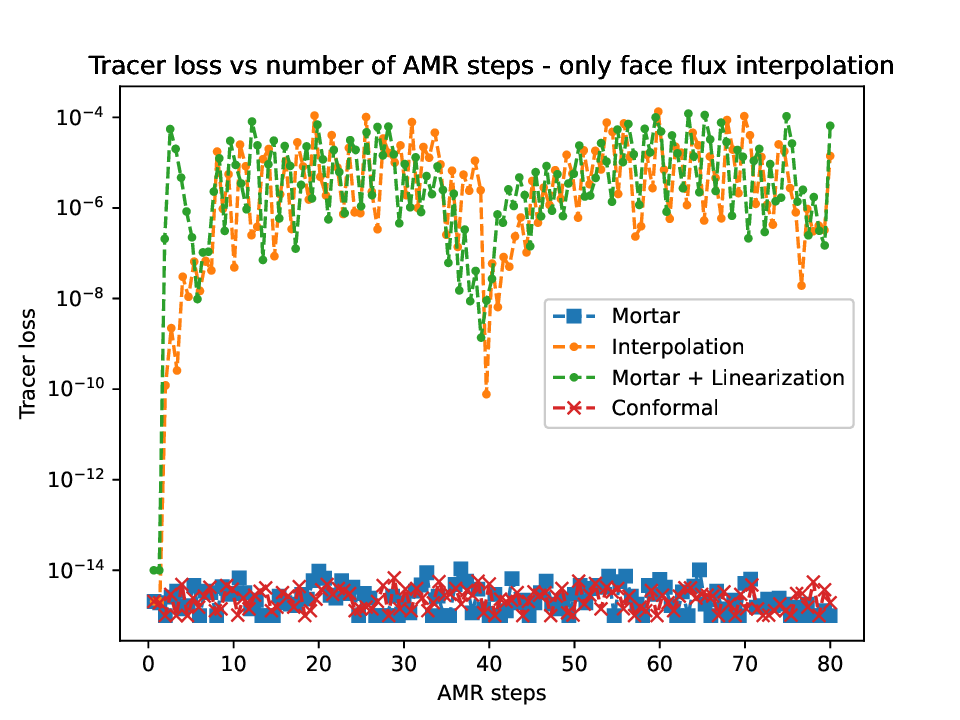}
\end{subfigure}
\caption{Tracer loss plot for the \citet{leveque1996} test case, using tree-based AMR with 80 regridding steps conducted over one full period. On the left, the mortar method demonstrates significantly better tracer conservation compared to the interpolation method. However, the mortar method applied to the linearized form  of the equations, a common practice in CFD, leads to the loss of its conservative property. On the right, tracer loss plot considering only face interpolations, i.e., disregarding error introduced from solution transfer (refinement + coarsening). In this scenario, the mortar method performs as well as the conformal grid. }
\label{leveque-cons}
\end{figure}

We have previously highlighted the importance of the interpolation method for maintaining global conservation in dGSEM. In Fig.\ \ref{leveque-cons}, we compare the mortar method with the interpolation method in terms of conserving the total amount of tracer for the \citet{leveque1996} test case. It is evident that the mortar method outperforms the interpolation method. The tracer loss is close to machine precision, a level of accuracy achieved by the conformal grid simulation. However, if we disregard solution transfer error, which is the error associated with transferring fields from parent to children elements and vice-versa, the mortar method's performance is virtually indistinguishable from that of the conformal grid.  It is important to note that when the mortar method is employed with a non-conservative formulation or an inaccurate linearization method, it loses its conservative properties, as expected. In CFD it is common to linearize the governing equations by considering $\mathbf U$ from the previous iteration and subsequently forming a system of linear equations to be solved implicitly. When the mortar method is applied to the linearized form of the equations, its conservative properties are compromised, and it performs similarly to the interpolation method. 
See Appendix \ref{appendix-linear} for more details. 

\subsection{Rising thermal bubble}
Until now, we have focused on test cases involving the passive advection of a tracer with prescribed wind fields and the isentropic vortex problem. While these are sufficient for validating the AMR implementations, atmospheric processes consist of complex phenomena including convection and wave motions of different scales.

The test case we consider here, the Robert Rising Thermal Bubble (RRTB) \citep{robert1993}, is frequently used to evaluate non-hydrostatic atmospheric dynamical cores.  It depicts a warm bubble evolving within a neutrally stratified atmosphere featuring a constant potential temperature $\theta_0$. The initial potential temperature within the bubble follows the cosine-bell shape as defined by Equation \ref{cosine-bell}. The initial conditions ensure hydrostatic balance, where pressure decreases with height, given by:

\[
p=p_0 \left( 1-\frac{gz}{c_p\theta_0} \right)^{c_p/R}.
\]

The simulation runs for 600 seconds, with the regridding step applied every 25 seconds based on potential temperature. 
The computational domain spans $[0m,1000m]^2$, with the bubble's center at $(x_c,z_c)=(500m,350m)$, 
a radius of $250m$, and $\theta_c=0.5K$, $\theta_0=300$K. 

\begin{figure}
\centering
\begin{subfigure}{0.24\textwidth}
\includegraphics[width=\linewidth]{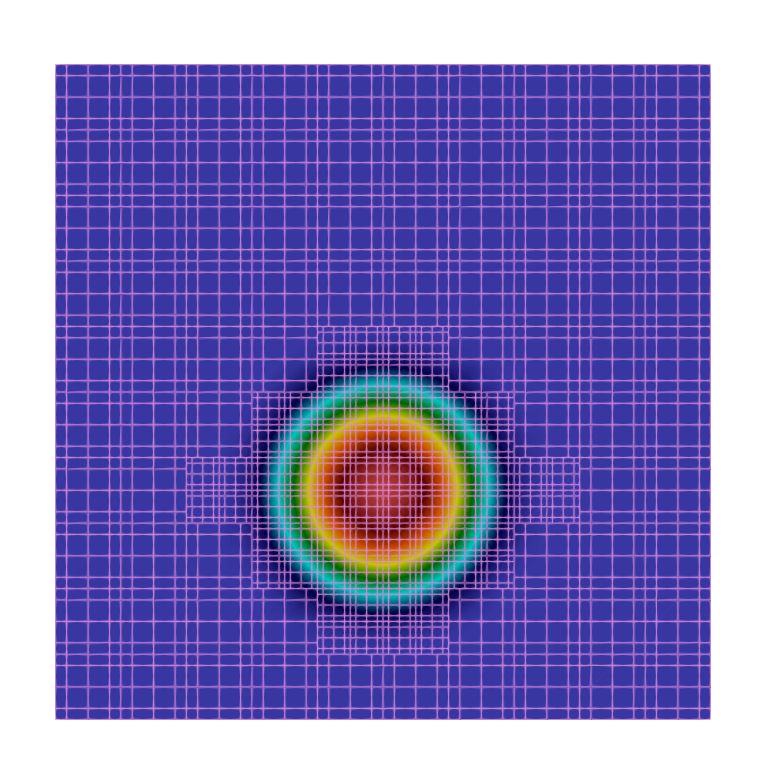}
\caption{t=0}
\end{subfigure}
\begin{subfigure}{0.24\textwidth}
\includegraphics[width=\linewidth]{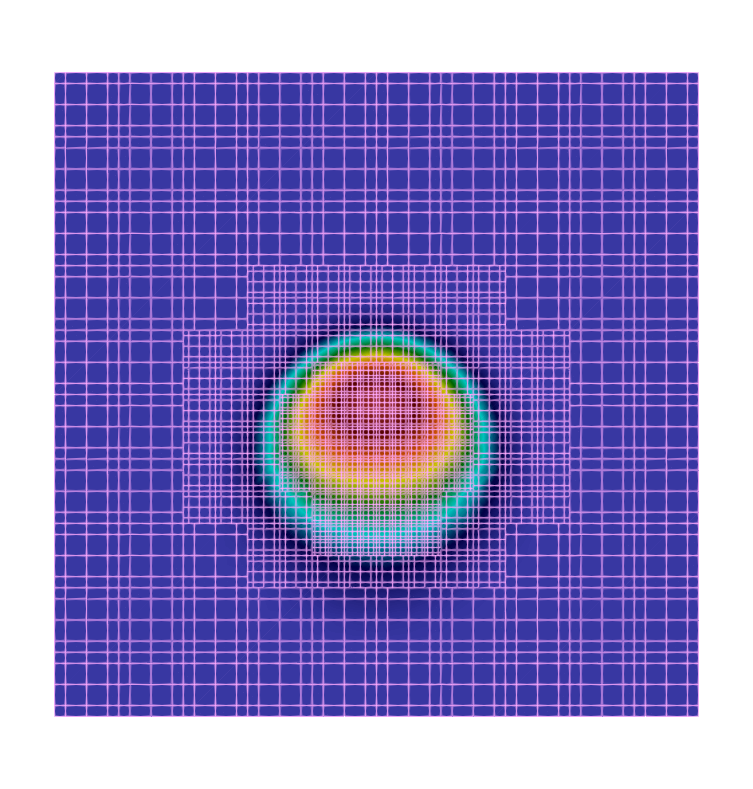}
\caption{t=200s}
\end{subfigure}
\begin{subfigure}{0.24\textwidth}
\includegraphics[width=\linewidth]{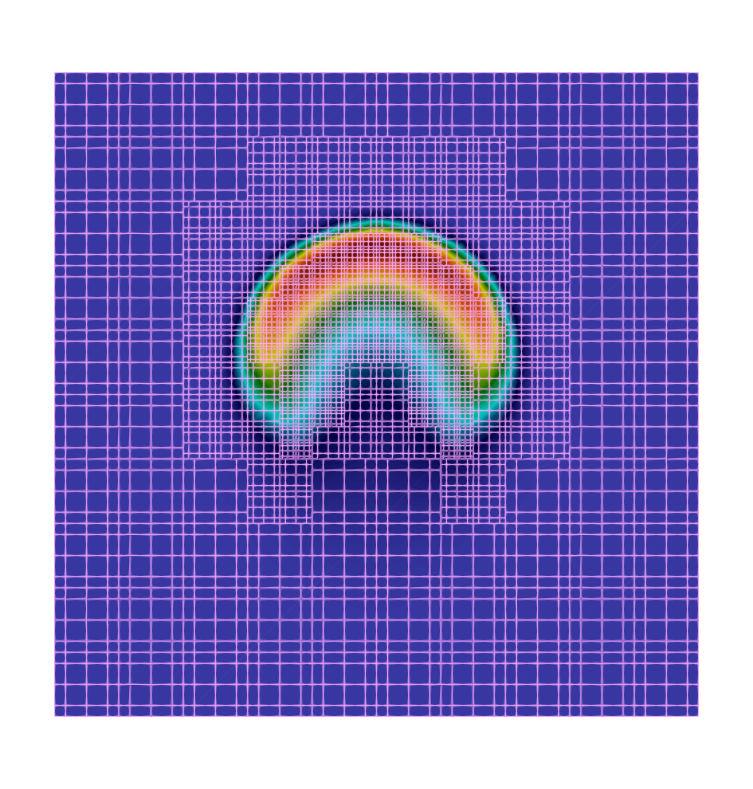}
\caption{t=400s}
\end{subfigure}
\begin{subfigure}{0.24\textwidth}
\includegraphics[width=\linewidth]{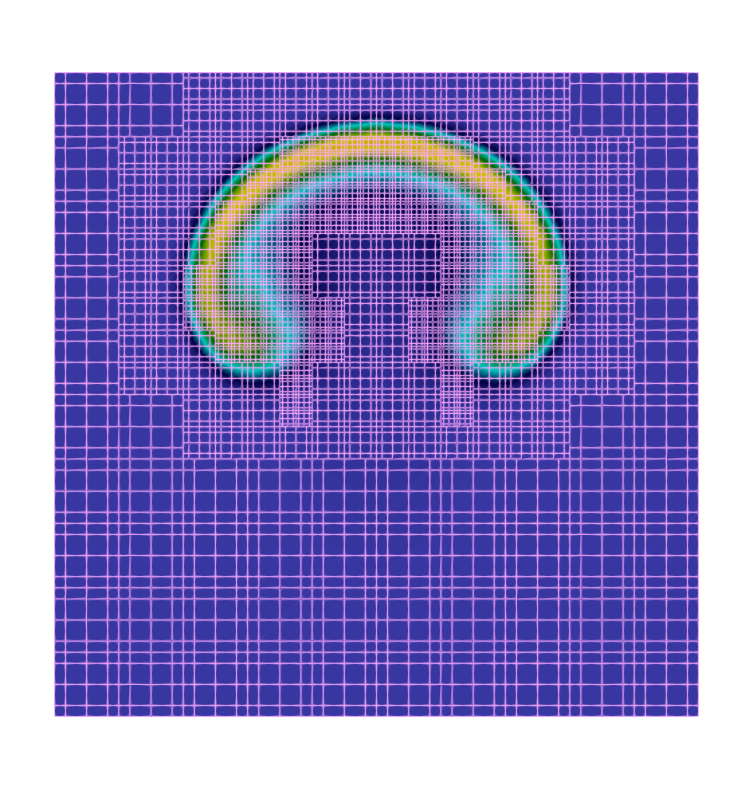}
\caption{t=600s}
\end{subfigure}
\begin{subfigure}{0.24\textwidth}
\includegraphics[width=\linewidth]{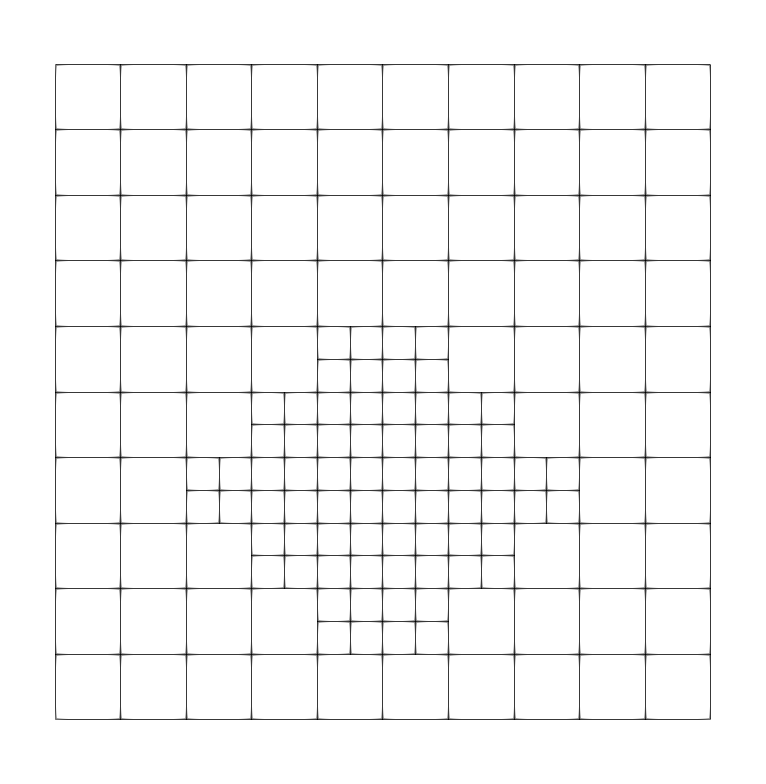}
\end{subfigure}
\begin{subfigure}{0.24\textwidth}
\includegraphics[width=\linewidth]{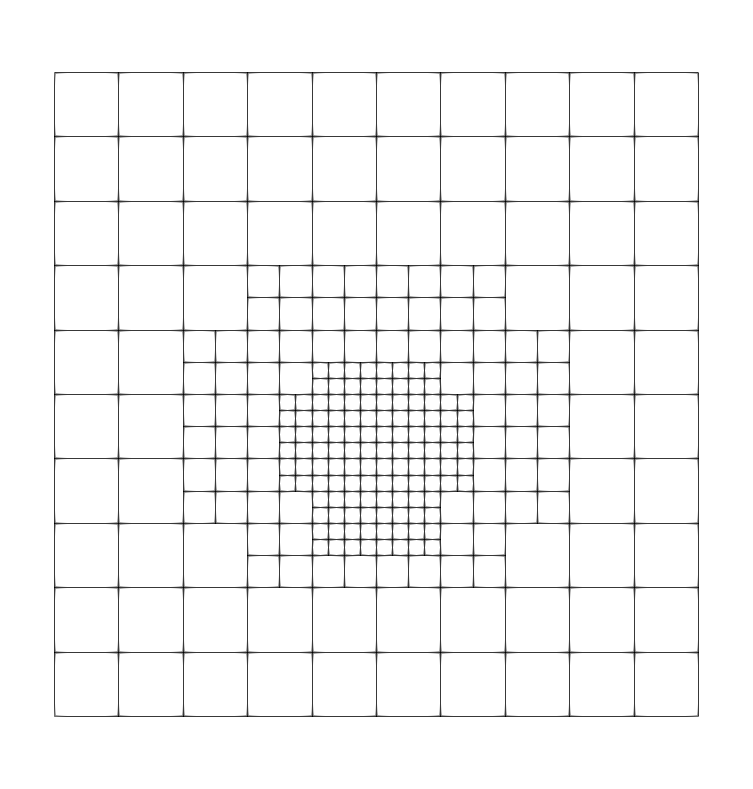}
\end{subfigure}
\begin{subfigure}{0.24\textwidth}
\includegraphics[width=\linewidth]{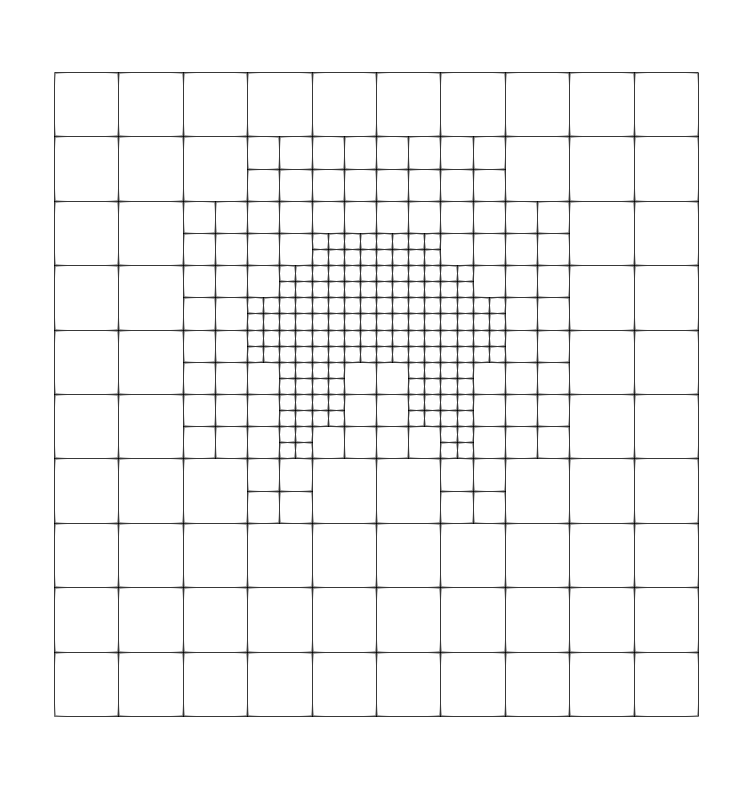}
\end{subfigure}
\begin{subfigure}{0.24\textwidth}
\includegraphics[width=\linewidth]{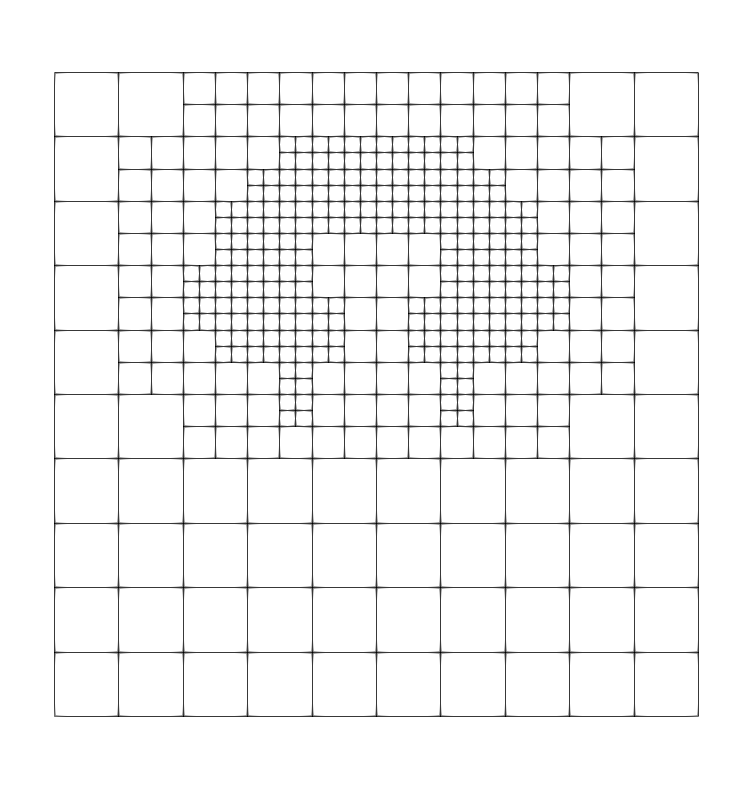}
\end{subfigure}
\begin{subfigure}{0.24\textwidth}
\includegraphics[width=\linewidth]{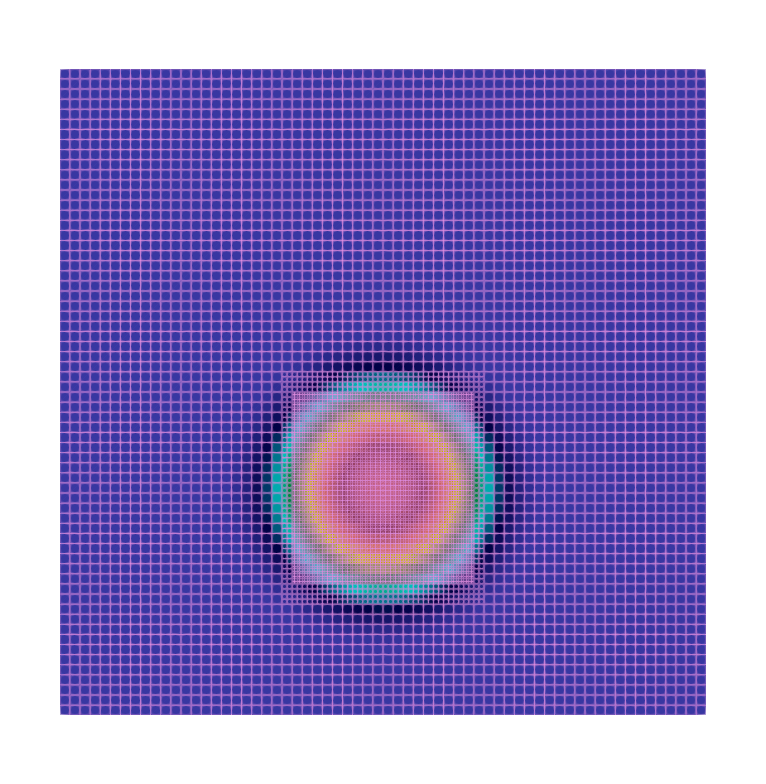}
\caption{t=0}
\end{subfigure}
\begin{subfigure}{0.24\textwidth}
\includegraphics[width=\linewidth]{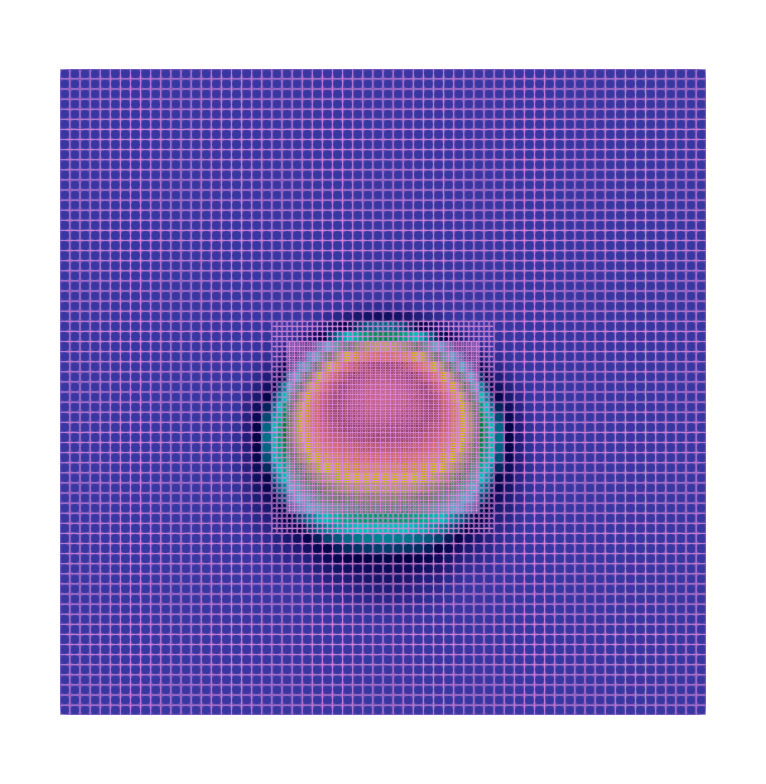}
\caption{t=200s}
\end{subfigure}
\begin{subfigure}{0.24\textwidth}
\includegraphics[width=\linewidth]{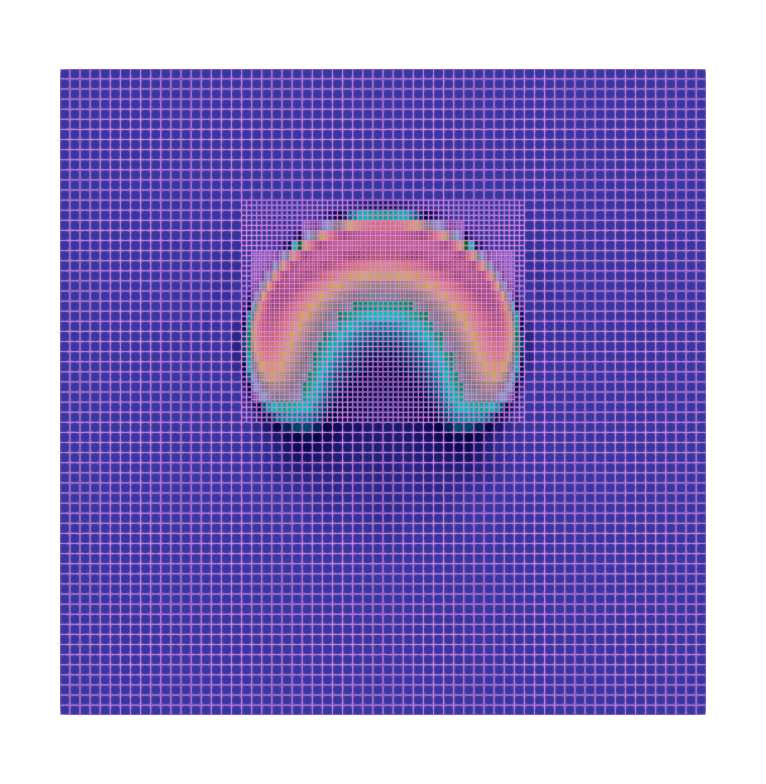}
\caption{t=400s}
\end{subfigure}
\begin{subfigure}{0.24\textwidth}
\includegraphics[width=\linewidth]{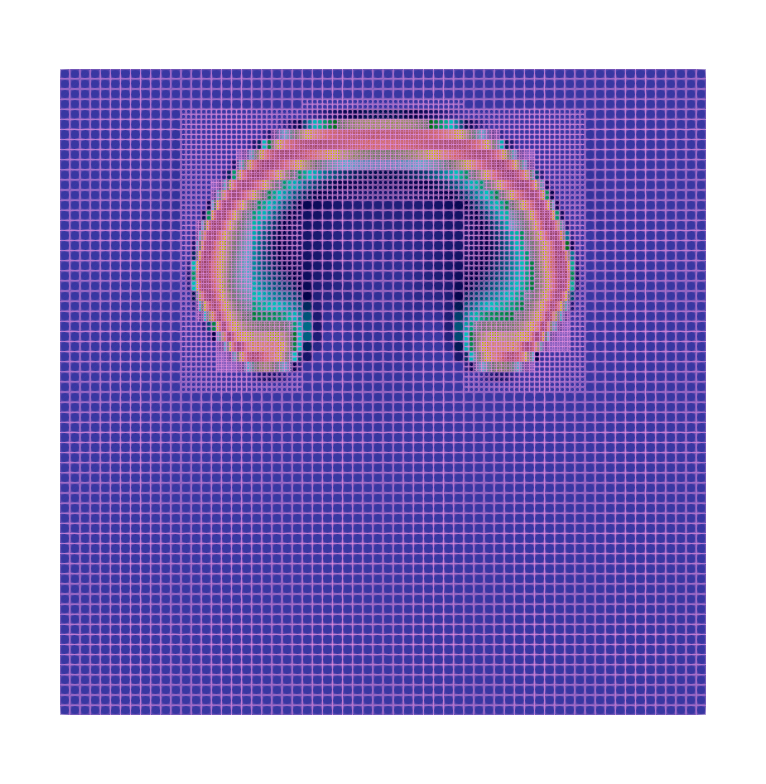}
\caption{t=600s}
\end{subfigure}
\begin{subfigure}{0.24\textwidth}
\includegraphics[width=\linewidth]{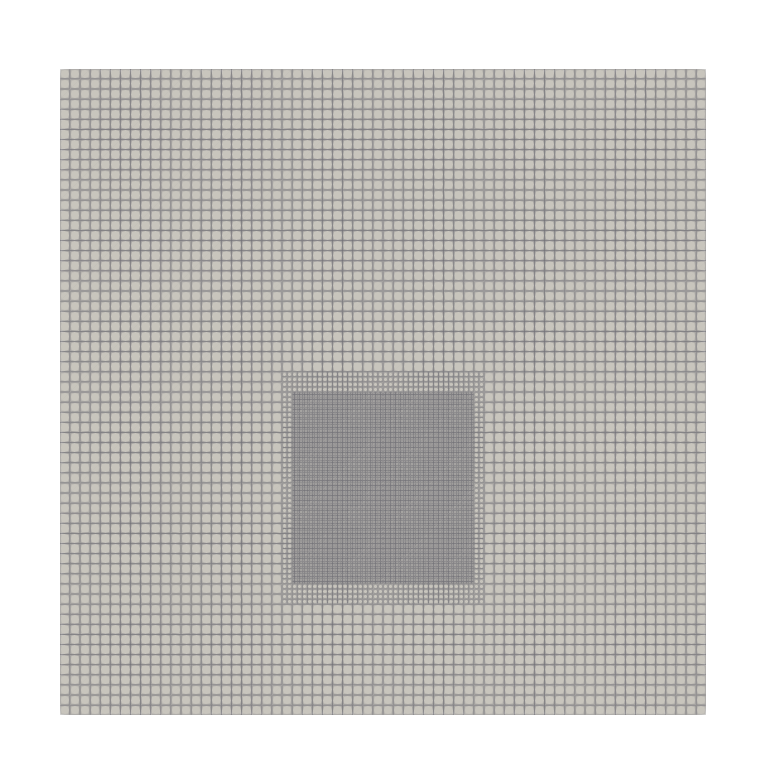}
\end{subfigure}
\begin{subfigure}{0.24\textwidth}
\includegraphics[width=\linewidth]{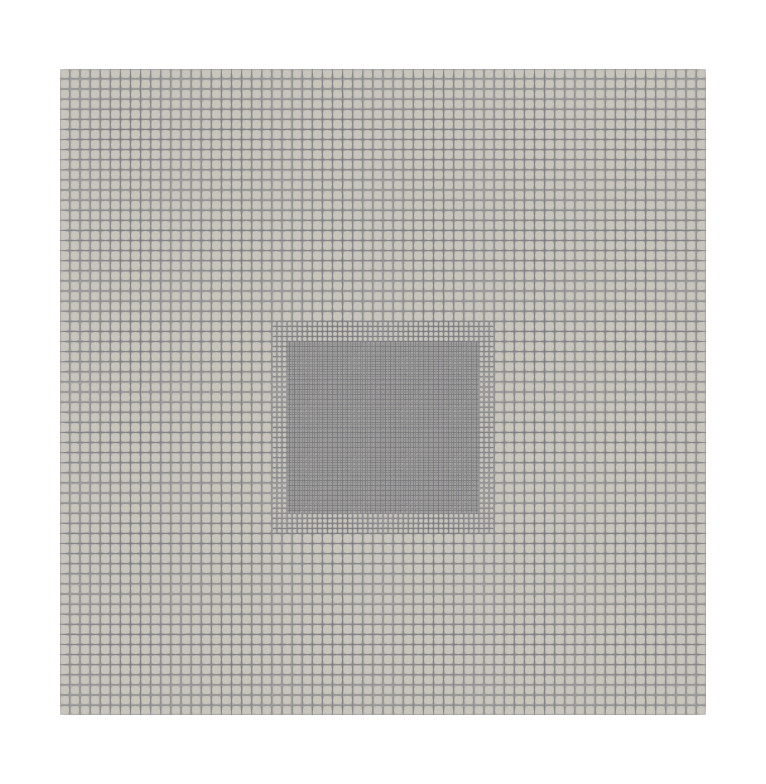}
\end{subfigure}
\begin{subfigure}{0.24\textwidth}
\includegraphics[width=\linewidth]{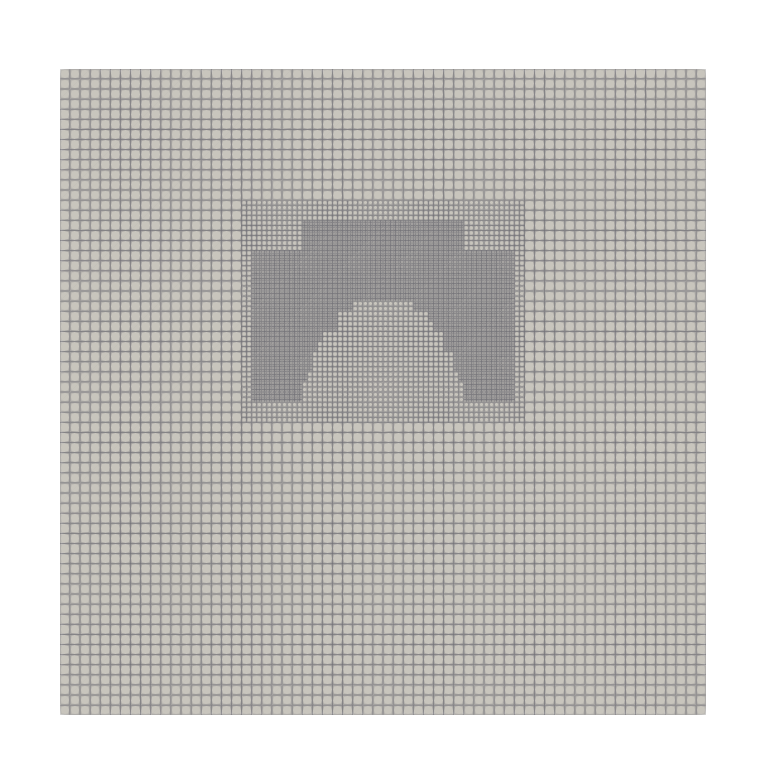}
\end{subfigure}
\begin{subfigure}{0.24\textwidth}
\includegraphics[width=\linewidth]{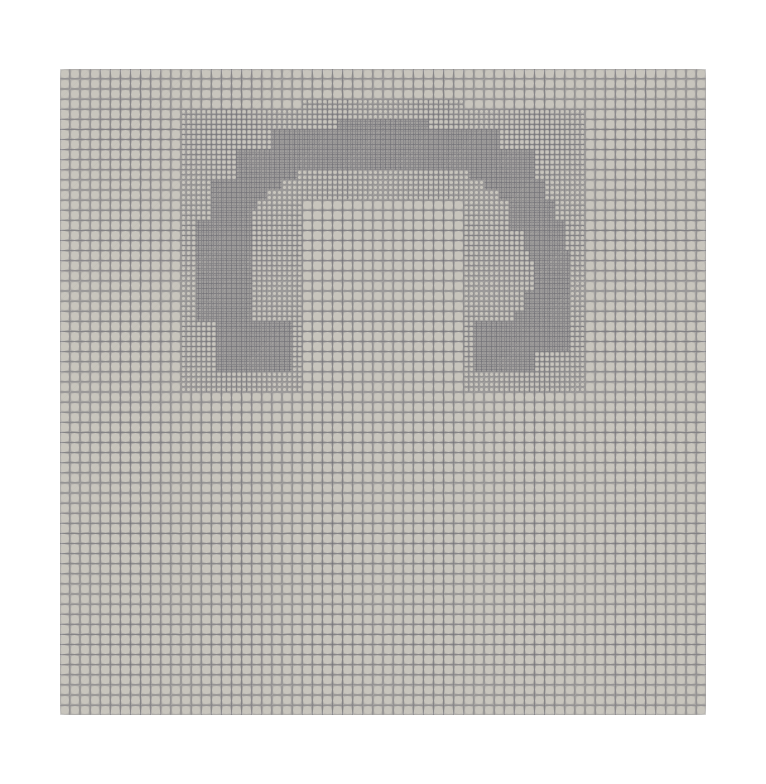}
\end{subfigure}
\caption{AMR for the Robert rising thermal bubble problem \citep{robert1993} with two levels of factor 2 refinement, using tree-based AMR (top two rows),  and level-based AMR (bottom two rows). The refinement criteria for both is based on potential temperature. The figures depict evolution of potential temperature, with AMR applied using a refinement criteria based on potential temperature.}
\label{srtb-amr}
\end{figure}

With the tree-based method the domain is divided into $10\times10$ elements with polynomial order $N=4$. 
The effective model resolution is around 25m. No-flux boundary conditions are enforced on all boundaries. 
To ensure stability of the simulation, artificial viscosity of $\mu=1.5m^2/s$ is applied.
With the level-based method we span the domain with $32\times4\times32$ cells.  
Both methods use two levels of factor 2 refinement, but we note that ERF refines only in the lateral directions.
(While AMReX itself allows refinement in all coordinate directions, ERF, like most other NWP models, 
does not allow refinement in the vertical direction.)
The results of both schemes using dynamic AMR are depicted in Figure \ref{srtb-amr}.

\begin{figure}
\centering
\includegraphics[width=0.3\linewidth]{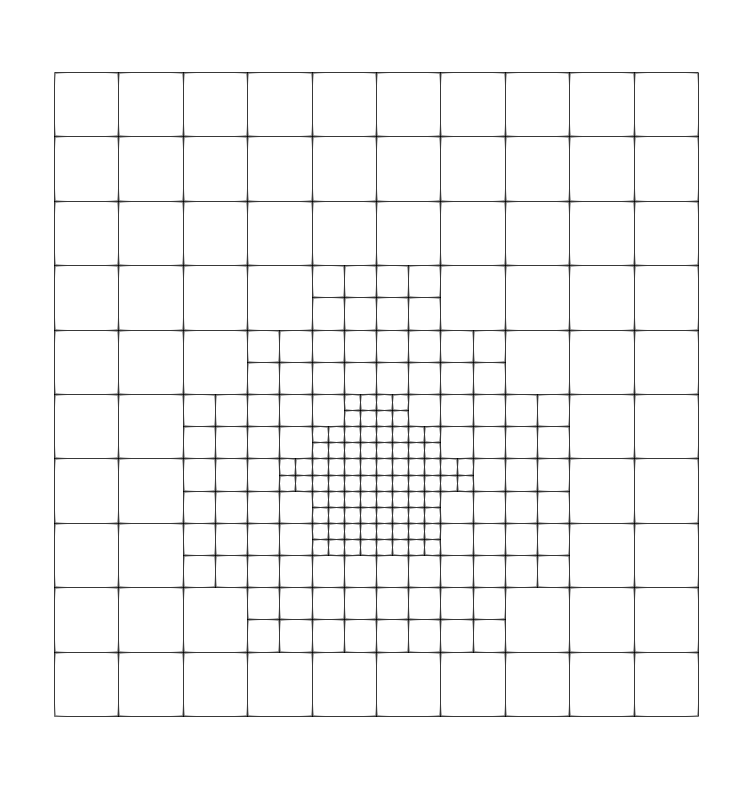}
\includegraphics[width=0.38\linewidth]{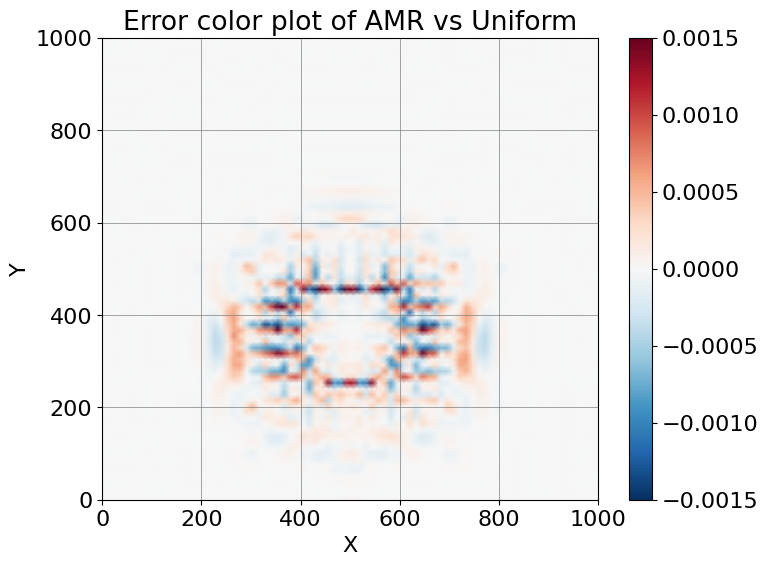}
\includegraphics[width=0.3\linewidth]{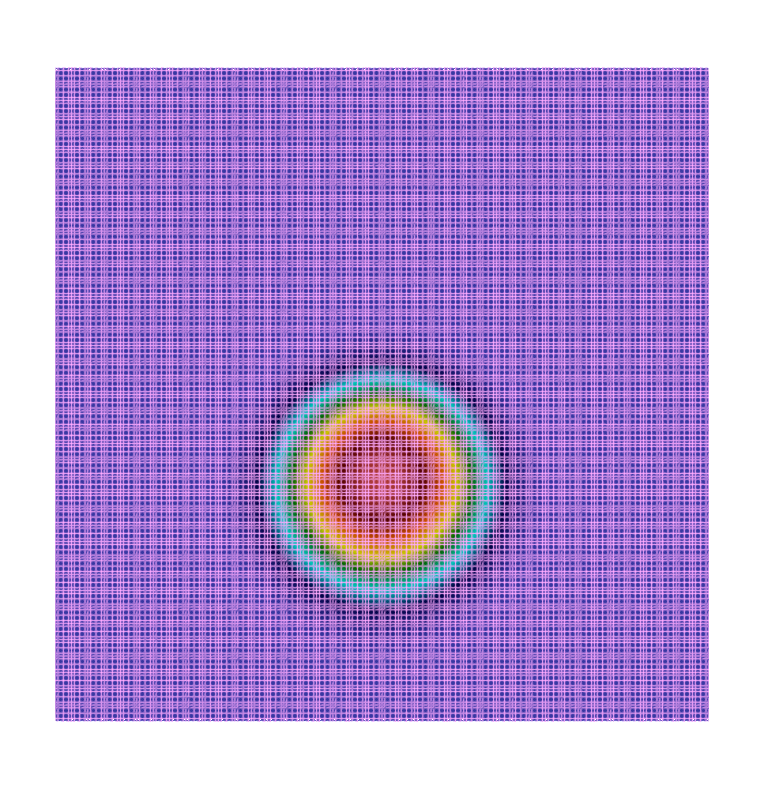}
\caption{Middle) potential temperature error at t=75 sec between the AMR simulation and the uniform grid simulation. Left) AMR grid. Right) $40\times40$ elements uniform grid simulation. Error is minimal outside of the bubble, where no other features are present, and in the inner core of the bubble where the same grid resolution is used as the uniform grid.}
\label{srtb-diff}
\end{figure}

Both the tree-based and level-based AMR schemes effectively track the bubble as it rises.
For the tree-based AMR, the 2:1 balance required for numerical stability and the mortar element method is maintained across all levels. 
While this problem lacks an analytical solution unlike previous test cases, a comparison against a uniform grid at the maximum AMR level (2 levels) provides insight into the accuracy and time-saving benefits of AMR. Figure \ref{srtb-diff} illustrates the difference between a $ 40 \times 40$ uniform grid simulation and the AMR simulation at 75 seconds. An average potential temperature difference of $10^{-4}$K is observed, with minimal errors in the inner core where the same resolution is used and slight differences in the outer ring of the bubble due to differing levels (1 vs. 2 levels). To minimize this error, 2-level refinement for the entire bubble region could be applied. The time-saving achieved by AMR is approximately tenfold, even without the use of multi-rate techniques, which would further enhance the AMR simulation. Although level-based AMR methods lend themselves more readily to a multi-rate implementation, it is also possible, albeit more complex, to apply this approach to non-overlapping tree-based AMR methods, as exemplified in \citet{seny2014,mugg2021} for multi-rate strategies in tree-based AMR coupled with dGSEM discretization.

\begin{figure}
\centering
\begin{subfigure}{0.49\textwidth}
\centering
\includegraphics[width=\linewidth]{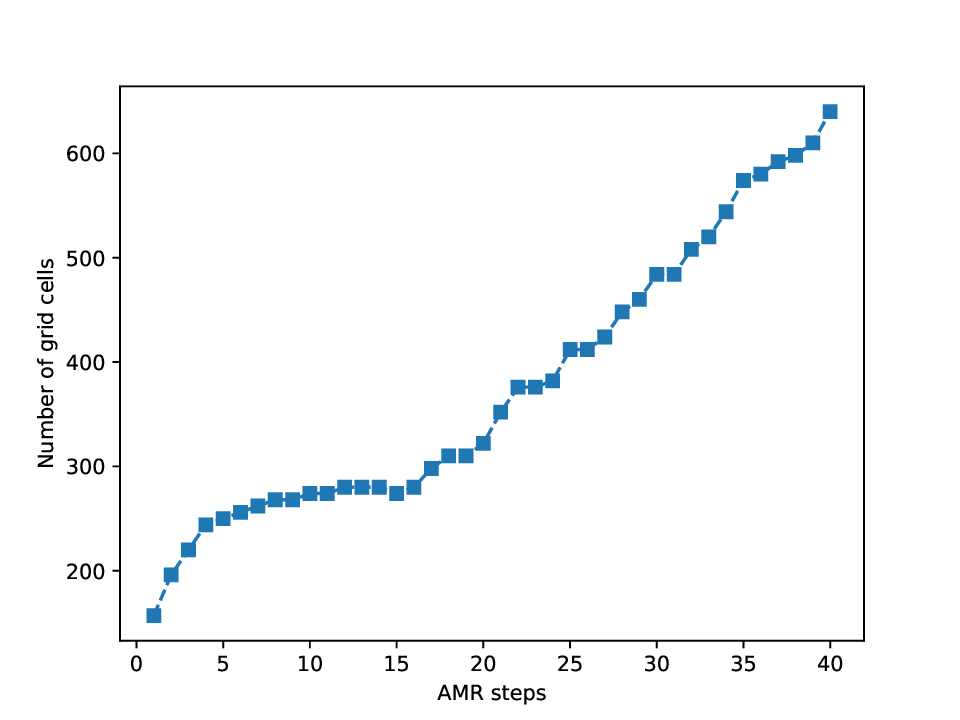}
\caption{Number of grid cells}
\label{srtb-cells}
\end{subfigure}
\begin{subfigure}{0.49\textwidth}
\centering
\includegraphics[width=\linewidth]{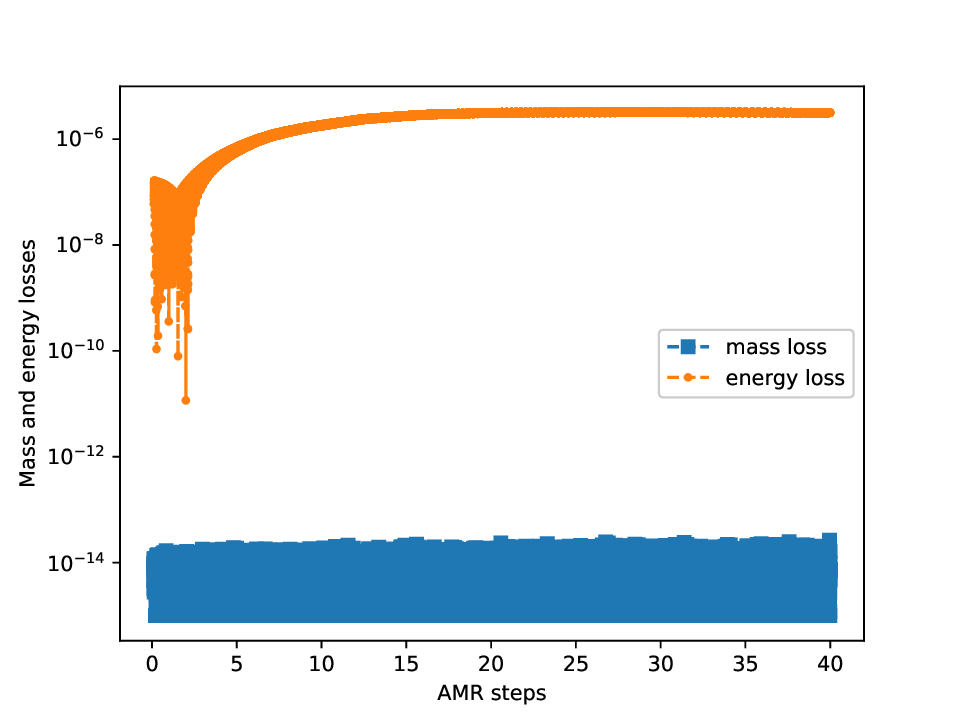}
\caption{Mass and energy loss}
\label{srtb-loss}
\end{subfigure}
\caption{(a) Number of grid cells for tree-based AMR simulation of the Robert rising thermal bubble problem. (b) Mass and energy losses  using the mortar method. The total mass loss is nearly at machine precision, but energy conservation is notably lower due to the governing equations being formulated on potential temperature.}
\label{loss-srtb}
\end{figure}

Figure \ref{loss-srtb} shows the number of grid cells used by the tree-based AMR simulation and the errors in mass and energy conservation. The  grid cell count increases over time as the bubble advects and expands under the influence of buoyancy. Mass is preserved up to nearly machine precision $10^{-15}$ after 1000 seconds, while energy conservation stands at $10^{-6}$. As expected, the total energy in the system is not conserved primarily due to the use of potential temperature as the prognostic variable.

\subsection{Flow on a sphere}
Spherical geometry introduces additional challenges as previously discussed in Sec. \ref{spherical-geom}. To validate the implementation of AMR on the sphere, two test cases are considered: a scalar advection problem with a reversing flow field and an acoustic wave propagation around the globe. The simulations in this section utilize only the tree-based AMR because the level-based ERF is not designed to function as a global circulation model. It is important to note that this distinction does not reflect negatively on level-based methods. Indeed, it is entirely feasible to utilize level-based AMR methods, such as AMReX on the sphere, as demonstrated, for instance, in \citet{ferguson2016} using the Chombo library. Both test cases are 2D shallow atmosphere simulations on a curved surface and require modifications if they are to be solved with the 3D equation sets (see Appendix \ref{appendix-2d} for details). 
\subsubsection{Advection on the sphere}
 To validate AMR for transport on the sphere, we employ a test case with a time-reversing background wind field in a similar fashion to the \citet{leveque1996} test case used in the preceding section. Specifically, we investigate the evolution of two symmetrically positioned cosine-bell tracers under a non-divergent wind field. The tracer's shape is defined by Equation \ref{cosine-bell}, with the radius now determined by the geodesic distance

\[
r = r_e \cos^{-1}[\sin\phi_i \sin\phi + \cos\phi_i\cos\phi\cos(\lambda - \lambda_i)].
\]
Here, $(\lambda_i,\phi_i)$ denote the longitude and latitude of the origins of the two cosine-bells, $r_e$ represents Earth's radius. The time-dependent wind field components are given by:

\[
\begin{aligned}
u = \frac{10 r_e}{T} sin^2(\lambda')sin(2\phi)g(t)+\frac{2\pi r_e}{T} cos(\phi)\\
v = \frac{10 r_e}{T} sin(2\lambda')cos(\phi) g(t)
\end{aligned}
\]
where $g(t)$ is defined as in Equation \ref{cosine-time}, $\lambda'=\lambda-2\pi t/T$.

The computational domain takes the form of  the traditional cubed-sphere, consisting of six tiles. Cubed-sphere grids avoid polar singularities incurred by latitude-longitude grids and also prove more convenient for adaptive mesh refinement. Each tile is initially subdivided into $30\times30$ elements with polynomial order $N = 4$. Each element represents an isoparametric curvilinear quadrilateral, meaning that the curved surface is approximated by the same basis function used to approximate fields inside each element. Given the two-dimensional nature of the problem, no boundary conditions are required. The time-integration extends for a duration of one full period T = 12 days. 

\begin{figure}
\centering
\begin{subfigure}{0.32\textwidth}
\centering
\includegraphics[width=\linewidth]{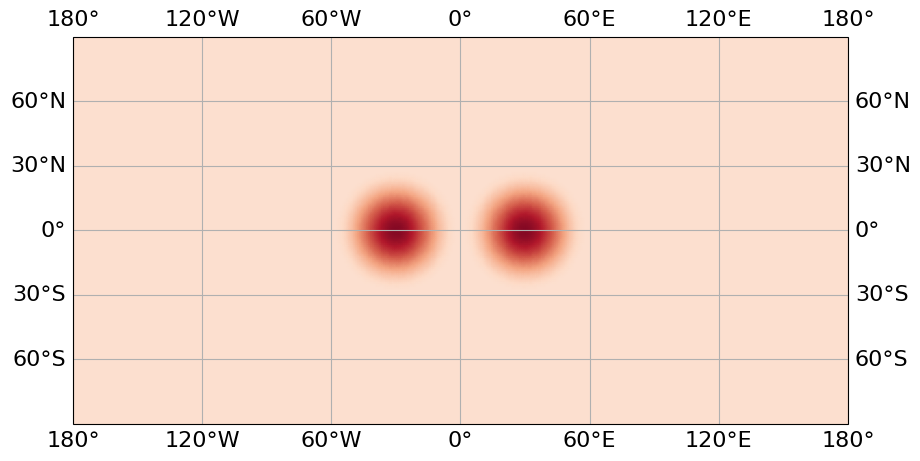}
\caption{Conformal t=0}
\end{subfigure}
\begin{subfigure}{0.32\textwidth}
\centering
\includegraphics[width=\linewidth]{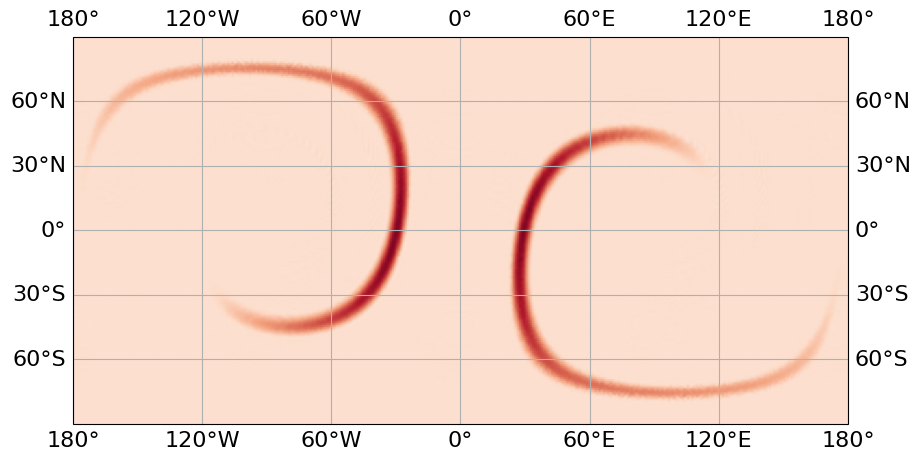}
\caption{Conformal t=T/2}
\end{subfigure}
\begin{subfigure}{0.32\textwidth}
\centering
\includegraphics[width=\linewidth]{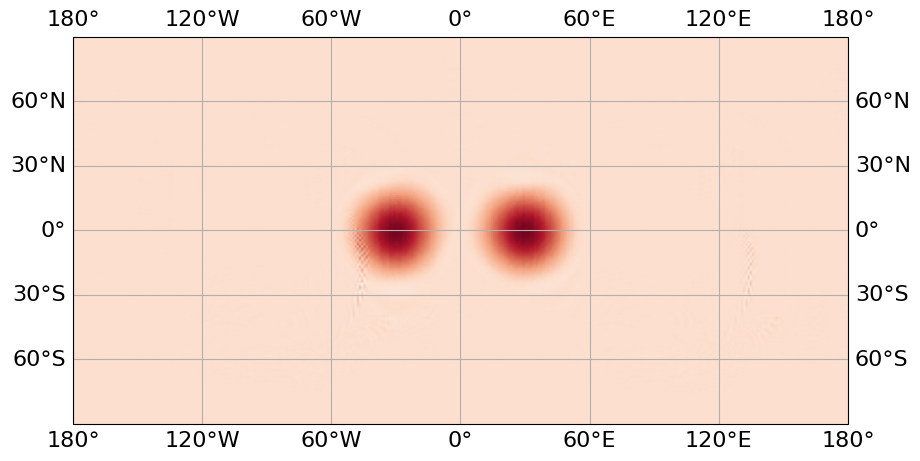}
\caption{Conformal t=T}
\end{subfigure}
\begin{subfigure}{0.32\textwidth}
\centering
\includegraphics[width=0.8\linewidth]{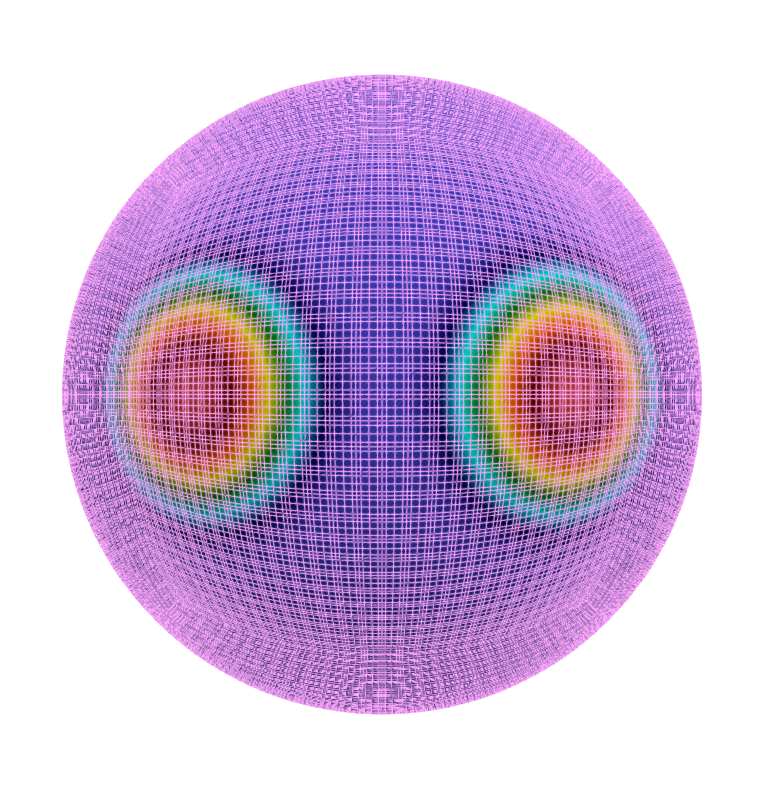}
\end{subfigure}
\begin{subfigure}{0.32\textwidth}
\centering
\includegraphics[width=0.8\linewidth]{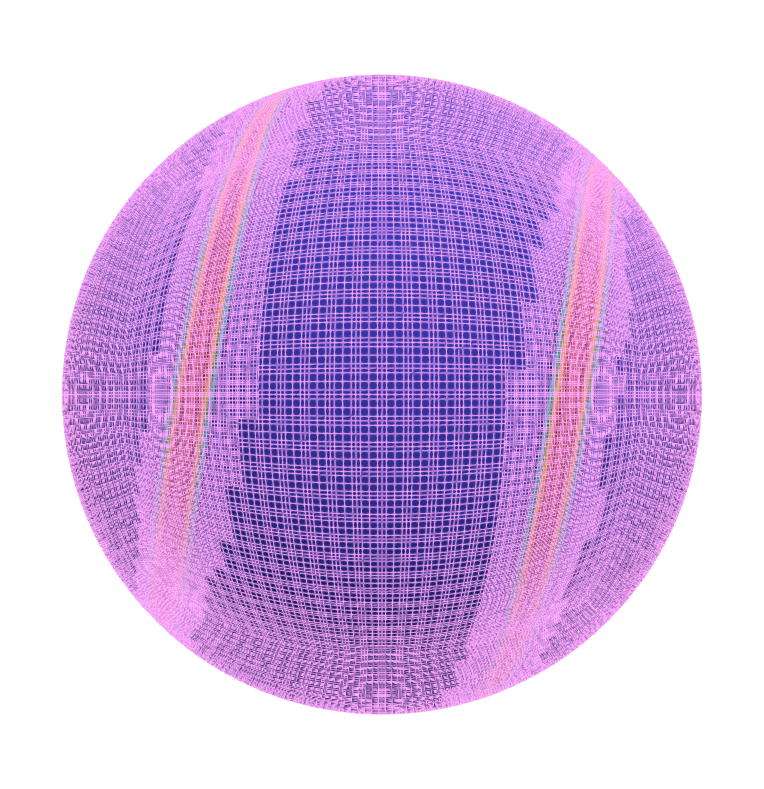}
\end{subfigure}
\begin{subfigure}{0.32\textwidth}
\centering
\includegraphics[width=0.8\linewidth]{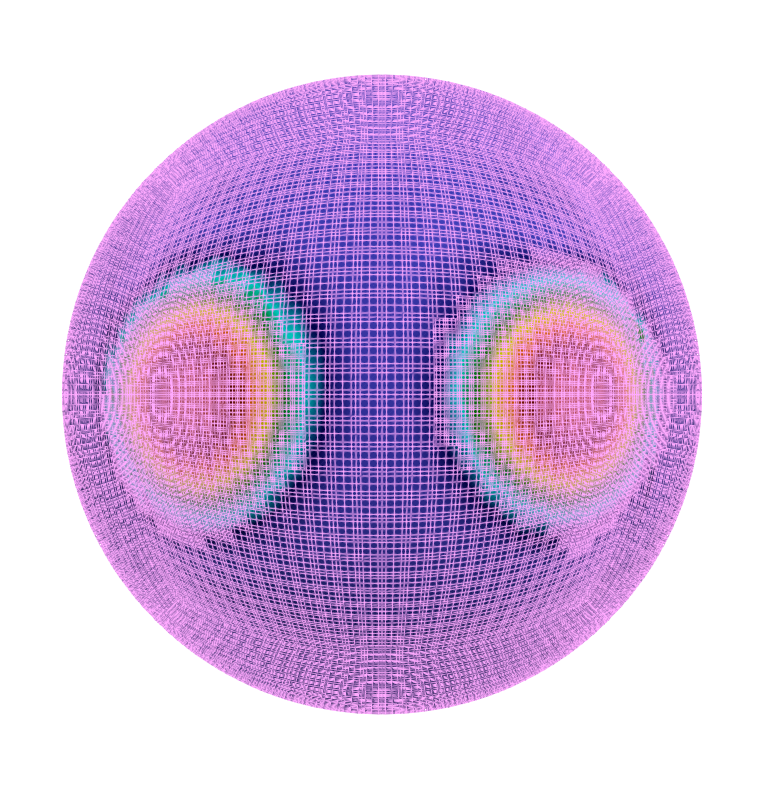}
\end{subfigure}
\begin{subfigure}{0.32\textwidth}
\centering
\includegraphics[width=\linewidth]{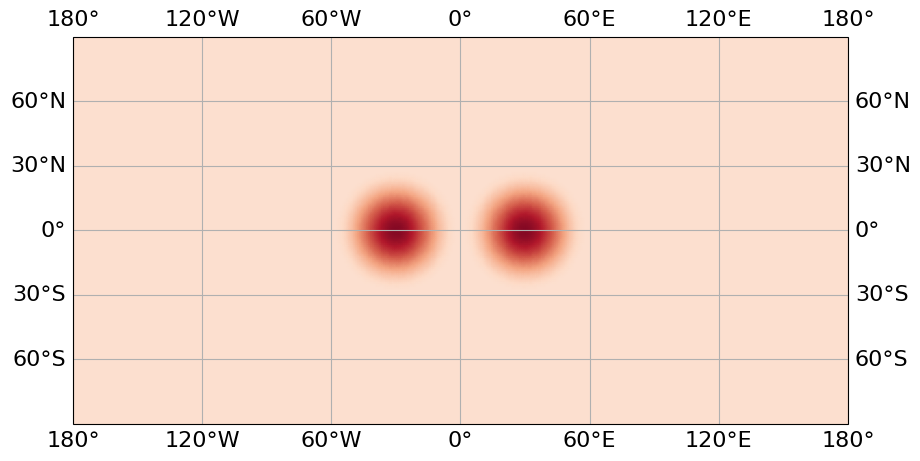}
\caption{AMR t=0}
\end{subfigure}
\begin{subfigure}{0.32\textwidth}
\centering
\includegraphics[width=\linewidth]{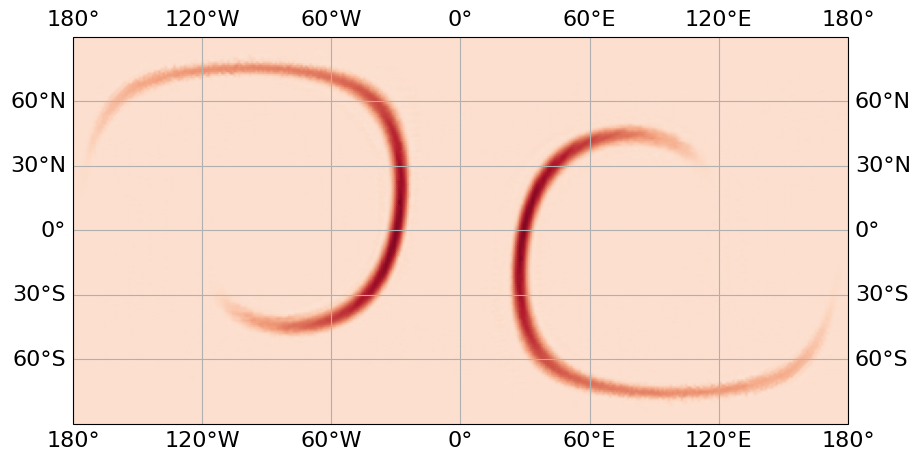}
\caption{AMR t=T/2}
\end{subfigure}
\begin{subfigure}{0.32\textwidth}
\centering
\includegraphics[width=\linewidth]{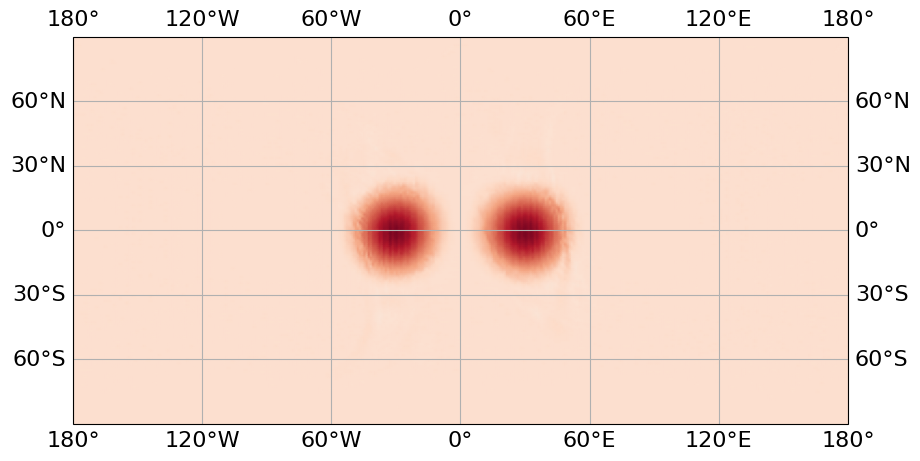}
\caption{AMR t=T}
\end{subfigure}
\begin{subfigure}{0.49\textwidth}
\centering
\includegraphics[width=0.97\linewidth]{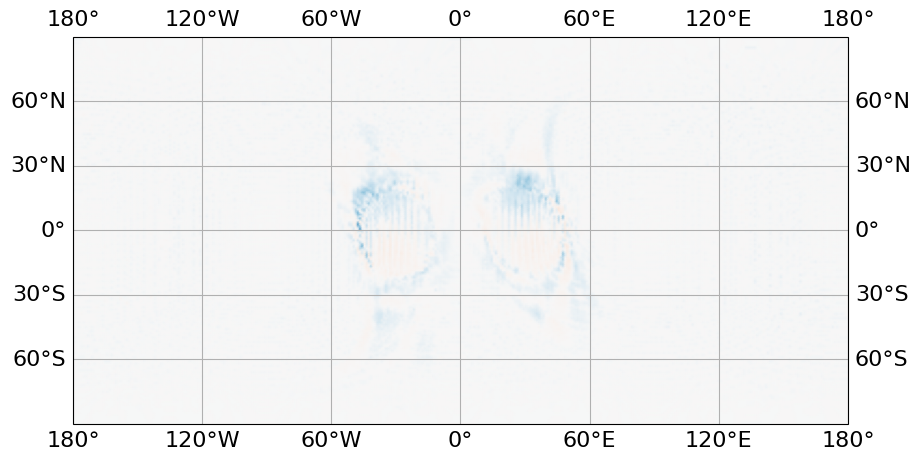}
\caption{Tracer error at t=T for the AMR simulation}
\label{advection-sphere-diff}
\end{subfigure}
\begin{subfigure}{0.49\textwidth}
\centering
\includegraphics[width=0.8\linewidth]{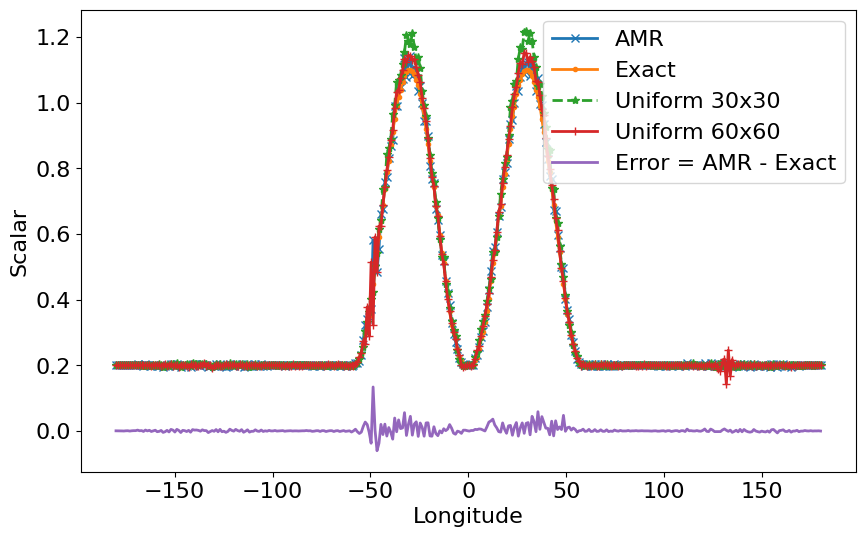}
\caption{Trace plot along the equator at t=T}
\label{advection-sphere-diff-line}
\end{subfigure}
\begin{subfigure}{0.5\textwidth}
\centering
\includegraphics[width=\linewidth]{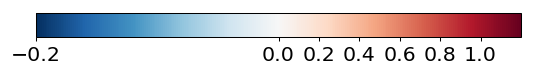}
\end{subfigure}
\caption{Evolution of the tracer advection problem, as studied by \citet{lauritzen2012}, over one complete period using the tree-based AMR method. Figures (a)-(c) depict the outcomes for a uniform grid simulation employing $60\times60$ elements per cubed-sphere face. In contrast, Figures (d)-(f) showcase the corresponding results for the AMR simulation, initialized with a grid of $30\times30$ elements per cubed-sphere face, along with a single level of refinement and a buffer zone of 4 cells. In both cases, a slender filament emerges at the midway point, followed by a reversal in wind flow that causes the filament to deform back to a shape closely resembling its original configuration. The tree-based AMR method effectively captures the evolving shape of the filament. Figure (h) demonstrates a close match between the AMR simulation results and the exact solution. Notably, the uniform grid solution with $30\times30$ elements exhibits a prominent peak compared to the other configurations.}
\label{advection-sphere}
\end{figure}

\begin{figure}
\centering
\begin{subfigure}{0.49\textwidth}
\centering
\includegraphics[width=\linewidth]{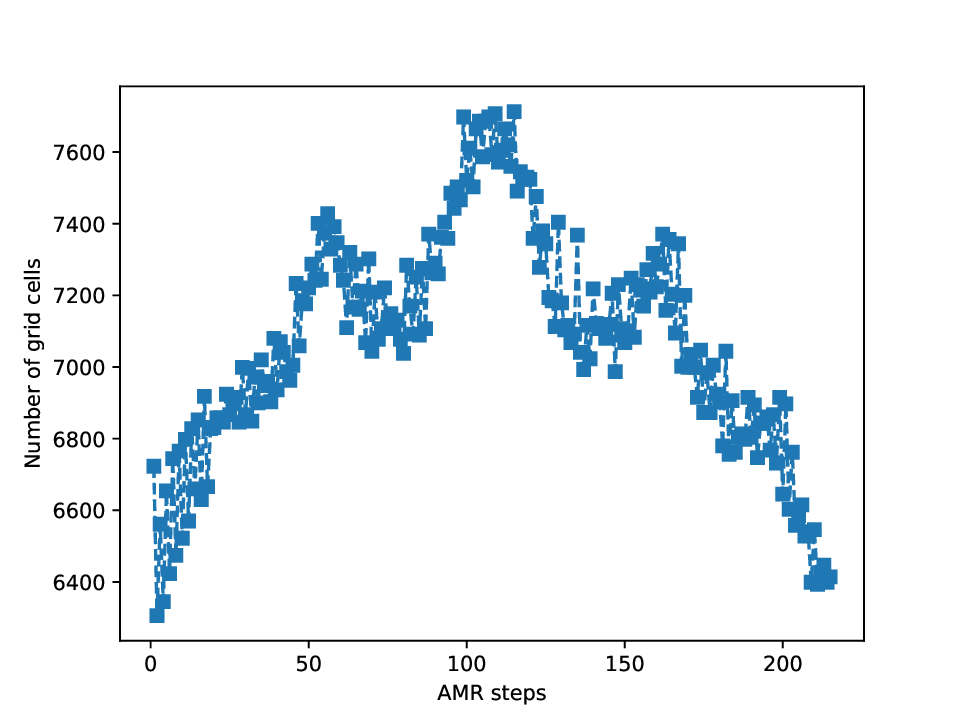}
\caption{Number of grid cells}
\label{sphere-cells}
\end{subfigure}
\begin{subfigure}{0.49\textwidth}
\centering
\includegraphics[width=\linewidth]{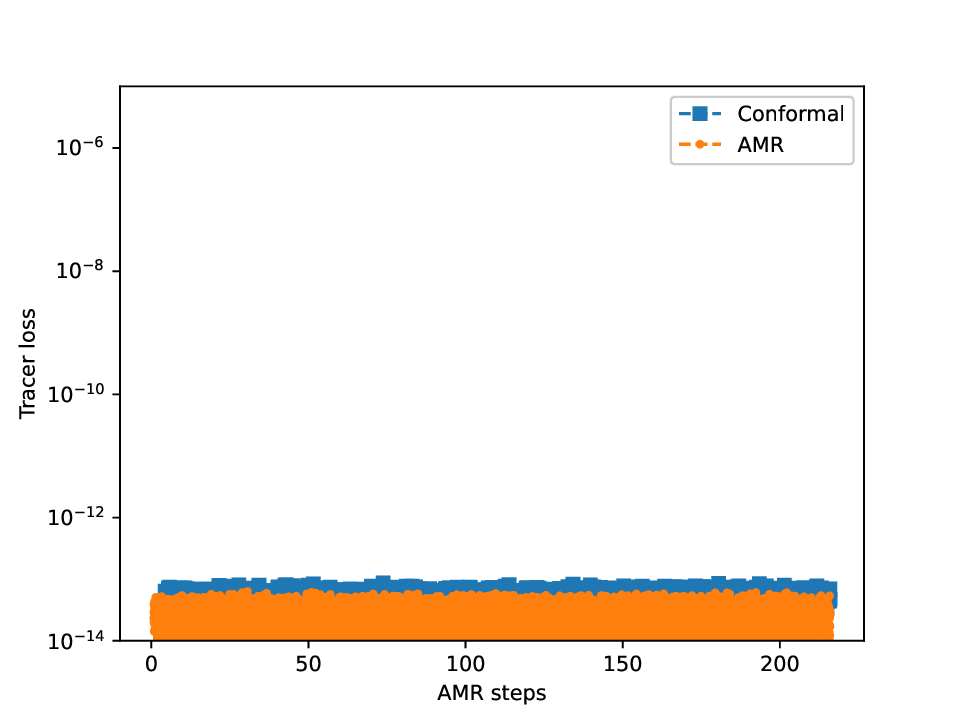}
\caption{Tracer loss}
\label{sphere-loss}
\end{subfigure}
\caption{(a) The number of grid cells utilized by AMR for the linear advection problem of \citet{lauritzen2012}. The grid cell count generally correlates with the length of the filament, peaking at half a period. (b) Tracer loss for conformal and AMR simulations using the mortar method. Tracer loss is in the order of machine precision for both runs.}
\label{loss-sphere}
\end{figure}

The result of the simulation at both half and full periods are presented in Figure \ref{advection-sphere}. The initially cosine-bell-shaped tracer undergoes significant stretching, forming long filaments at $t=T/2$, closely reproducing the result of \citet{lauritzen2012}. In both the conformal and AMR simulations, the tracers are essentially restored to their original shape at the full period. It is also evident that the tree-based AMR simulation, featuring a buffer zone of 4 cells, adeptly tracks the tracer at all stages. The number of grid cells increases as the tracer gets stretched and deformed, and is roughly proportional to its length, as depicted in Figure \ref{sphere-cells}. The tracer error plot in Figure \ref{advection-sphere-diff} for the AMR simulation reveals minimal errors, especially outside of the bubble regions. 

Furthermore, we run two uniform grid simulations at resolutions of $30\times30$ and $60\times60$ elements per cubed-sphere face. As observed in Figure \ref{advection-sphere-diff-line}, the AMR result aligns more closely with the uniform grid simulation using $60\times60$ elements than the coarser uniform grid simulation. The time saving achieved through AMR is roughly $4\times$ less compared to the fine grid simulation using $60\times60$ elements.

In Figure \ref{sphere-loss}, we observe that the tracer loss is close to machine precision at $10^{-14}$ for both the uniform grid and AMR simulations. The latter was possible as a result of  the specific solution strategy discussed in Section \ref{sec-tree-amr} that ensure conservation for simulations on curved surfaces.  This methodology has allowed us to attain similar level of tracer conservation comparable to that observed for advection problems within a box.

\subsubsection{Acoustic wave on the sphere}
For the validation of AMR on the sphere to solve the Euler equations, we utilize a test case from \citet{tomita2004} featuring an acoustic wave traveling around the globe. The initial state is hydrostatically balanced with an isothermal background potential temperature of $\theta_0$=300K. A perturbation pressure $p'$

\[
p'=f(\lambda,\phi) g(r)
\]
is superimposed on the reference pressure where $f$ is defined by the cosine-bell function in  Equation \ref{cosine-bell}, with $h_{max}$ = 100 Pa. The function $g$ is defined as

\[
g(r)=\sin\left(\frac{n_v \pi r}{r_T}\right)
\]
where $n_v=1$, $r_c$ is one third of the radius of the earth $r_e=$6371km, and a model altitude (top) of $r_T$=10km. The simulation runs for 24 hours, with regridding applied every hour based on the pressure perturbation. The grid is cubed-sphere with an initial resolution of $30\times30$ elements per cube face and a polynomial order $N = 4$. A no-flux boundary condition is used at the top and bottom surfaces.

\begin{figure}
\centering
\begin{subfigure}{0.32\textwidth}
\centering
\includegraphics[width=\linewidth]{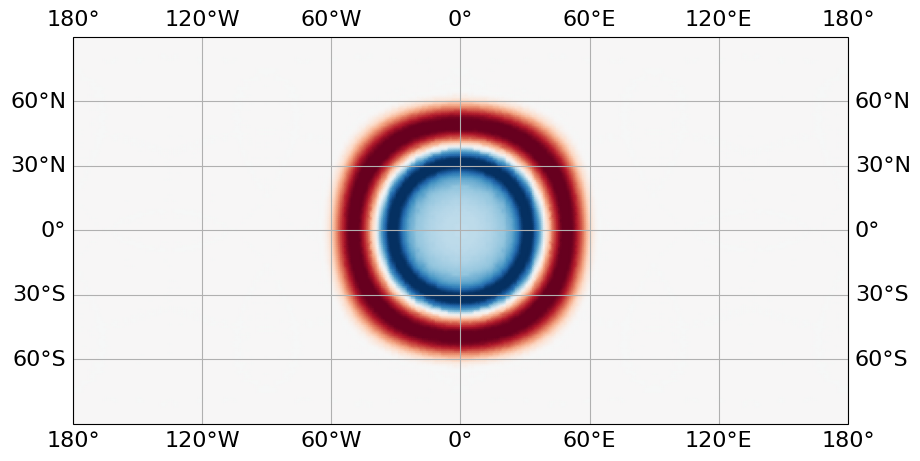}
\caption{t=4 hrs}
\end{subfigure}
\begin{subfigure}{0.32\textwidth}
\centering
\includegraphics[width=\linewidth]{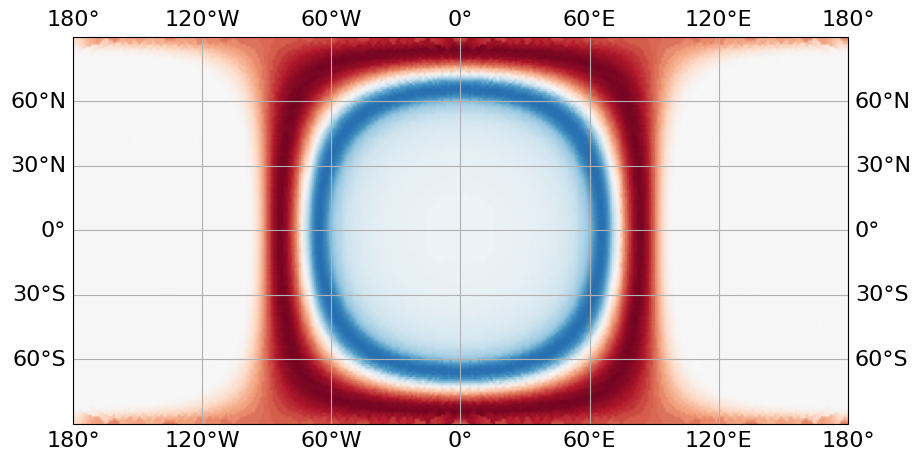}
\caption{t=7 hrs}
\end{subfigure}
\begin{subfigure}{0.32\textwidth}
\centering
\includegraphics[width=\linewidth]{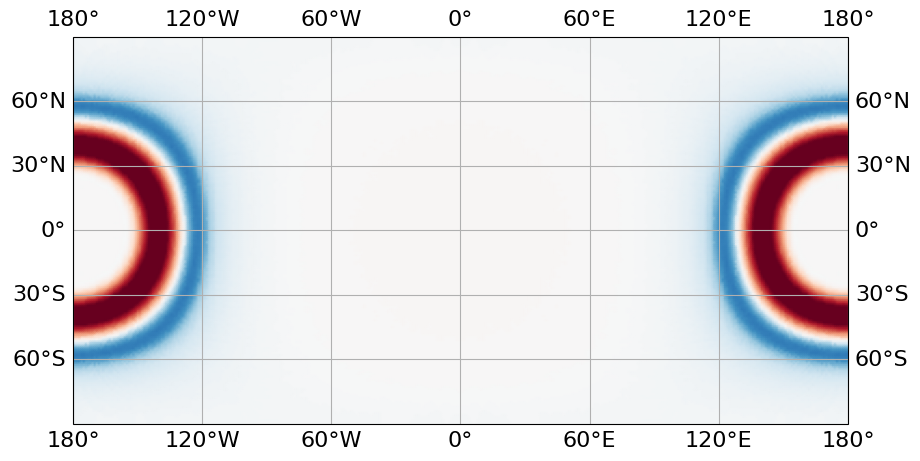}
\caption{t=12 hrs}
\end{subfigure}
\begin{subfigure}{0.32\textwidth}
\centering
\includegraphics[width=0.8\linewidth]{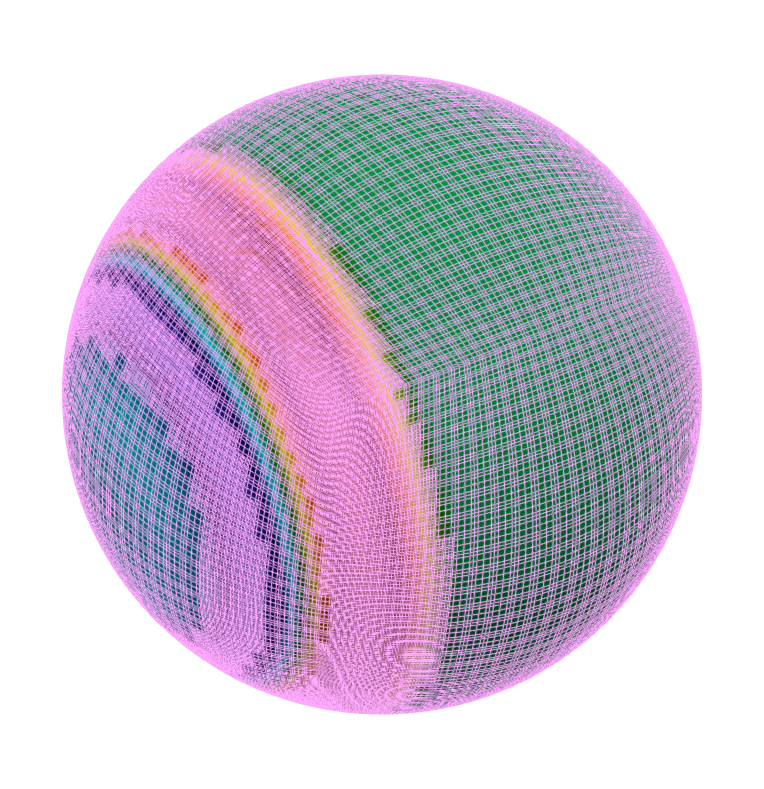}
\end{subfigure}
\begin{subfigure}{0.32\textwidth}
\centering
\includegraphics[width=0.8\linewidth]{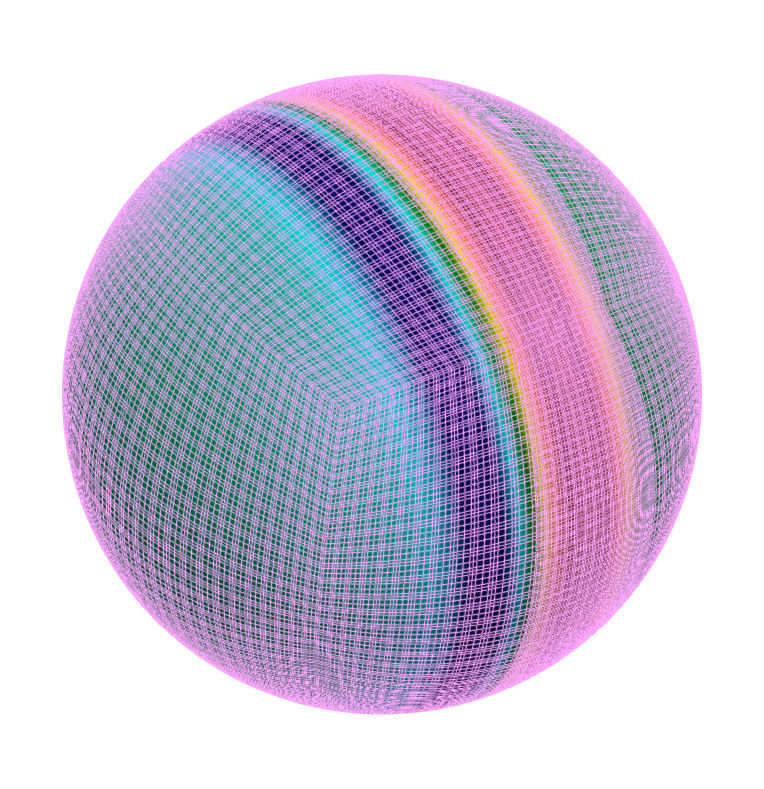}
\end{subfigure}
\begin{subfigure}{0.32\textwidth}
\centering
\includegraphics[width=0.8\linewidth]{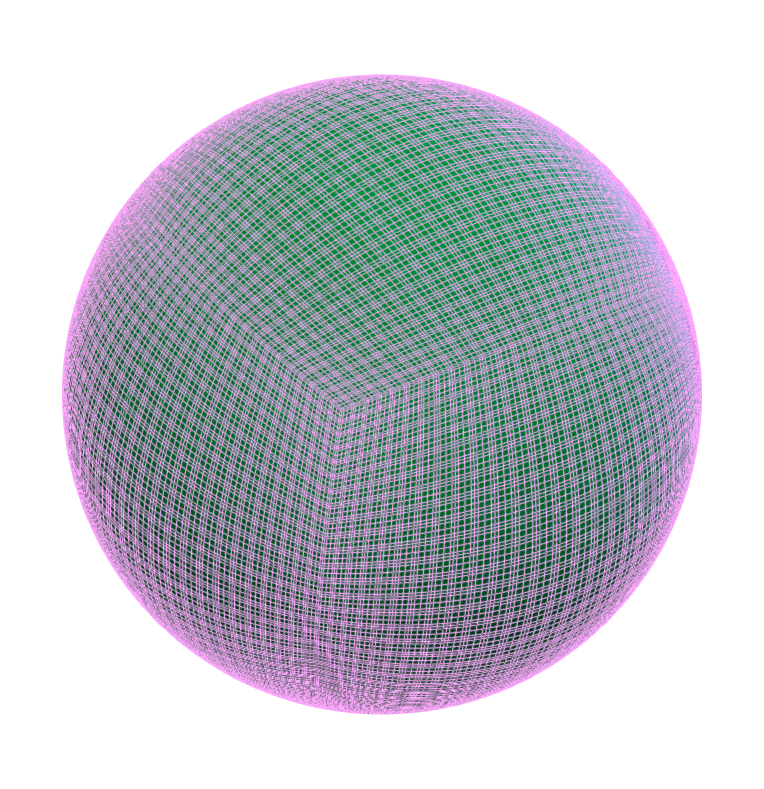}
\end{subfigure}
\begin{subfigure}{0.49\textwidth}
\centering
\includegraphics[width=0.96\linewidth]{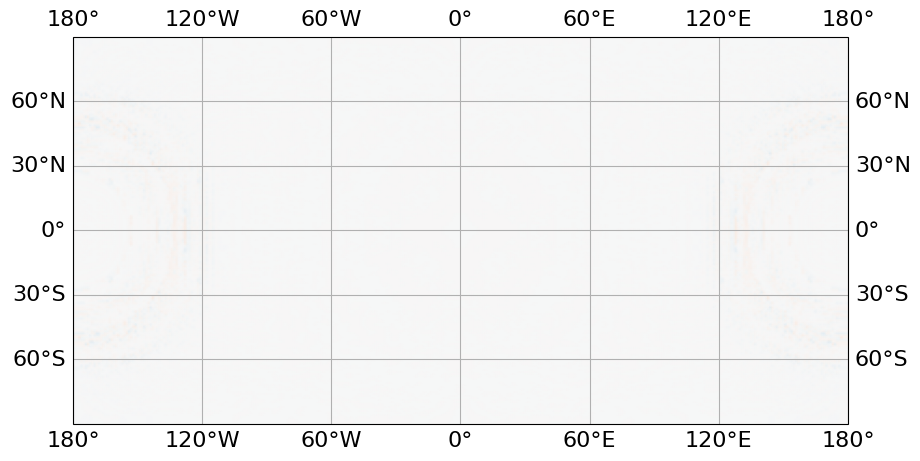}
\caption{Pressure error b/n AMR and uniform grid at t=12 hrs}
\label{acoustic-diff}
\end{subfigure}
\begin{subfigure}{0.49\textwidth}
\centering
\includegraphics[width=0.8\linewidth]{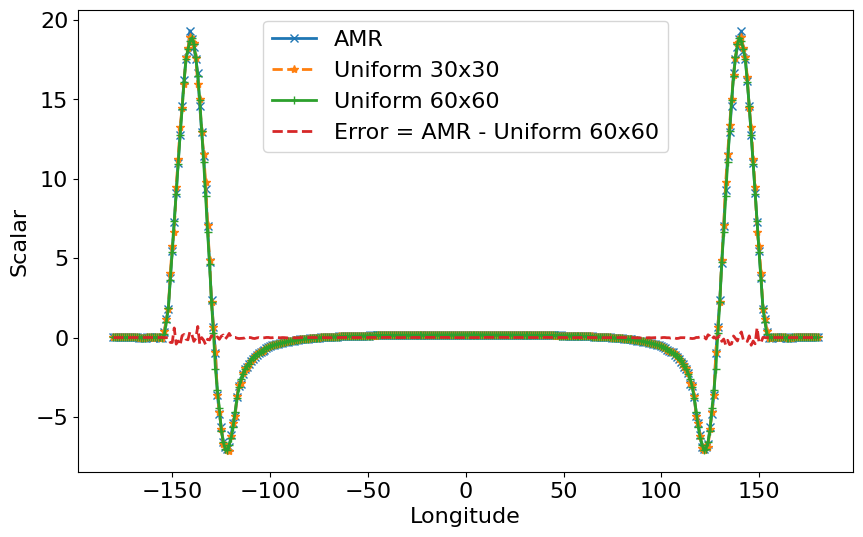}
\caption{Pressure plots along the equator at t=12 hrs}
\label{acoustic-diff-line}
\end{subfigure}
\begin{subfigure}{\textwidth}
\centering
\includegraphics[width=0.4\linewidth]{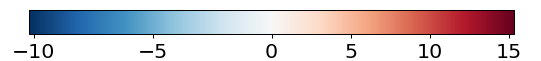}
\end{subfigure}
\caption{Propagation of an acoustic wave on the sphere solved using tree-based AMR and one level of refinement. The pressure perturbation after 4 hours, 7 hours and 12 hours is displayed. The simulation uses a cubed sphere grid with $30\times30$ elements and polynomial order of 4. Bottom figures (d and e) show difference in pressure between the AMR simulation with with $30\times30$ elements and the conformal grid simulation with $60\times60$ elements.}
\label{acoustic-sphere}
\end{figure}

\begin{figure}
\centering
\begin{subfigure}{0.48\textwidth}
\centering
\includegraphics[width=\linewidth]{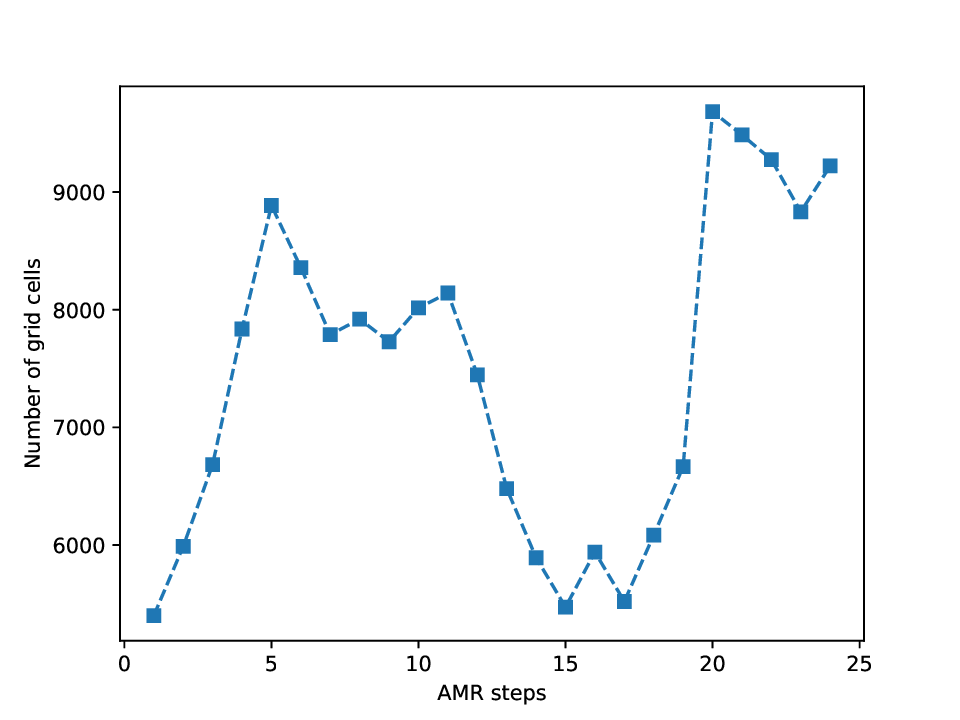}
\caption{Number of grid cells}
\label{acoustic-cells}
\end{subfigure}
\begin{subfigure}{0.48\textwidth}
\centering
\includegraphics[width=\linewidth]{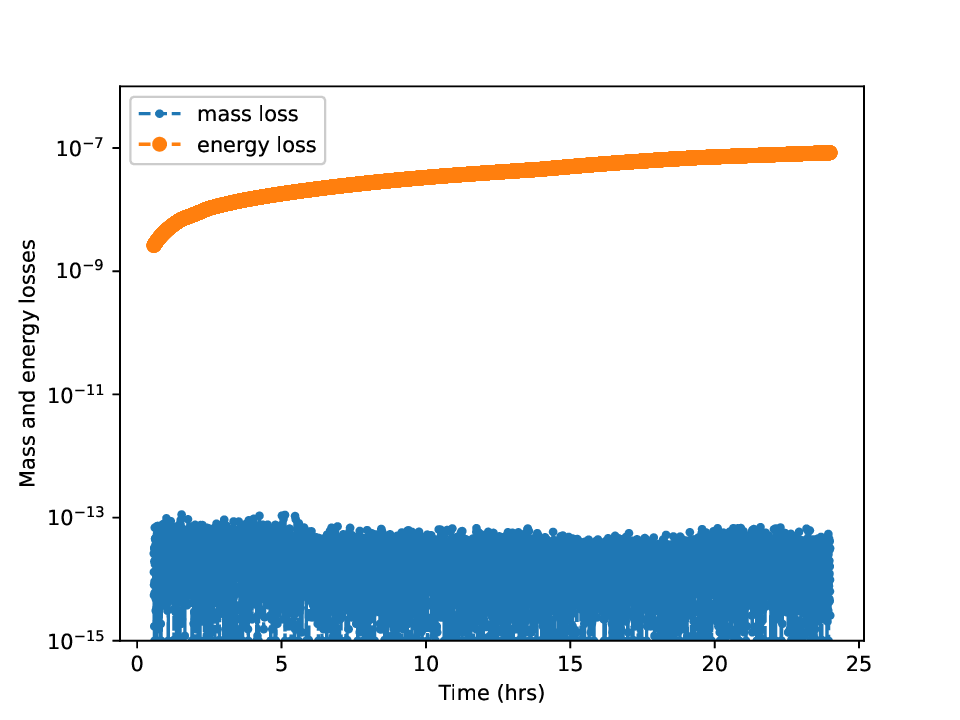}
\caption{Mass and energy loss}
\label{acoustic-loss}
\end{subfigure}
\caption{(a) Number of grid cells for the AMR simulation of the acoustic wave propagation problem from \citet{tomita2004}. (b) Mass and energy losses for the AMR simulations using the mortar method.}
\label{loss-acoustic}
\end{figure}

Figure \ref{acoustic-sphere} illustrates the simulation results, showcasing the propagation of the wave from the initial perturbation location to the antipode. Adaptive mesh refinement, using pressure perturbation as the criterion, effectively tracks the wave rings at all hours. To validate the results, a visual comparison was made with plots from \citet{tomita2004} at different hours, and the time it would take for the wave to reach the antipode was computed. The speed of sound was calculated based on the initial condition: $a = \sqrt{\gamma p / \rho}$ = 347.3 m/s. The observation confirmed that the wave took approximately 16 hours to converge at the antipode, consistent with the AMR simulation. In Figure \ref{acoustic-cells}, we can also see that the AMR simulation introduced only few refined cells between hours 15 and 17. This was followed by a sudden surge in refined cells, driven by the subsequent expanding wave, even surpassing the peak at hour 5. Regarding mass and energy conservation, once again the mass loss is close to machine precision at $10^{-15}$ and the energy loss stands at $10^{-7}$. 


Additionally, we conducted conformal grid simulations with resolutions of $30\times30$ and $60\times60$ elements per cubed-sphere face, and compared them against the AMR simulation. Figures \ref{acoustic-diff} and \ref{acoustic-diff-line} illustrate that the results are nearly indistinguishable from each other. In contrast to the previous test case on the sphere, the coarser uniform grid resolution did not produce noticeable differences with the other simulations. Once again, AMR resulted in significant time savings, approximately $3.7\times$ faster than the finer conformal grid simulation in this case.




\section{Choice of AMR framework for NWP}
\label{sec-choice}
Here, we summarize our insights on various aspects of selecting and integrating an AMR framework into existing NWP models, including both GCMs and LAMs.

The level-based AMR technique has served as the de facto standard for mesh refinement in NWP since the seminal work of \citet{berger1984,skamarock1989}. Its principal advantage lies in the refinement and clustering algorithm, which generates a hierarchy of logically rectangular boxes, preserving the structured nature of the grids. This facilitates the reuse of existing code for nested grids, enabling separate time advancement with subcycling. Furthermore, it permits the use of different physics options, providing flexibility in the development process of regional models. While conceptually straightforward, its implementation is not necessarily simple especially when considering scalability and performance portability.

Numerous implementations of mesh refinement in NWP models involve building the infrastructure from scratch, as exemplified by \citet{skamarock2019,harris2013,zangl2022}. Some instances are justified due to the absence of an existing library to handle the grid type and AMR infrastructure, such as \citet{zangl2022} utilizing an unstructured icosahedral grid. Implementing AMR is often a time-consuming process, as echoed in our experience implementing AMR in NebulaSEM from scratch. The resultant code may not scale well and may not leverage fine-grained computing architectures, such as GPUs. Mesh management code, in particular, is typically not GPU-friendly. Additionally, devising parallel communication strategies for a scalable AMR implementation is intricate in both approaches as outlined in \citet{p4est} for cell-based and \citet{amrex2021} for level-based methods. 

Furthermore, several from-scratch implementations are solely designed for static mesh refinement, rendering them unsuitable for tracking dynamic features like tropical cyclones. To enhance cyclone tracking, the specification of the nested grid placement is often done manually through input files. An adaptive mesh refinement technique can automate this process, utilizing criteria such as relative vorticity for refinement \citep{ferguson2016}.

Both AMR strategies have exhibited scalability to the entirety of exascale supercomputers. Several AMReX based applications have showcased exascale level scalability. For instance, WarpX, a 2022 ACM Gordon Bell Prize winner, has demonstrated a weak scaling efficiency of about 75\% on supercomputers including Frontier, Summit, Fugaku, and Perlmutter \citep{fedeli2022}. Another AMReX based application, PeleC, exhibited weak scaling efficiency of 65\% on Summit \citep{frahan2023}. We anticipate a similar level of scalability for NWP applications that utilize AMReX like  ERF.

The choice of one-way vs two-way nesting is an interesting one, and both approaches are available in ERF.
In many scientific domains two-way nesting is the default, but not necessarily so in NWP.  
The question of whether nudging is needed and if so, how it should be done, is also an open question. 
The answer, as in so many of these questions, is that the answer depends on the numerics used in the single-level time evolution. 
In  \citet{skamarock2019,zangl2022}, nudging techniques are employed to enhance one-way nesting by aligning prognostic variables in the nest 
with corresponding parent values. However, this issue may not be a practical concern in most cases. 
Some argue that two-way nesting does not significantly improve results; 
for instance, \citet{soriano2002} reported worse outcomes with two-way coupling. 
Errors often arise in transferring boundary data from parent to nested grid and 
in integrating the fine-grid solution into the parent grid. \citet{zangl2022} 
contend that using high-order Radial Basis Function (RBF) interpolation minimizes disturbances 
induced at nest boundaries to the point of being negligible for real applications. 
\citet{harris2013} opt to skip mass updates from the nest to the parent grid to 
ensure mass conservation in the two-way nested FV3 global circulation model.

Both level-based and tree-based AMR can ensure conservation of quantities whose evolution is written in conservation form.
Mass conservation in particular is an important property for long-term forecasts. In ERF, for example, if
two-way nesting is selected, then the refluxing operation that occurs after data on the coarse and fine levels reach the same time (e.g. after one coarse step and two fine steps with subcycling and a refinement ratio of 2) ensures conservation by exactly fixing the flux mismatch that occurred during the time evolution.   In one-way nesting, refluxing is not used, thus conservation is not ensured.

In the tree-based method,
this required careful design of interpolation methods, such as the ``mortar method" of \citet{kopriva1996}, as a straightforward point-to-point interpolation fell short in satisfying conservation. The latter approach can also induce instability at non-conformal element boundaries, similar to the challenges faced by the level-based methods at nest boundaries. Solution transfer between parent and child elements also requires careful consideration particularly on curved surfaces like the sphere. 


\section{Conclusions}
\label{sec-conclusions}
In this study, we have investigated two distinct adaptive mesh refinement methods within the context of the multi-scale challenges posed by numerical weather prediction. Our focus was on implementing and validating these methods through a series of benchmark problems, with a specific emphasis on their applicability to the dynamic refinement required for tracking localized features, such as tropical cyclones.

The first AMR method explored is the level-based approach, as implemented in the Energy Research Framework (ERF) code. The ERF code, which utilizes the massively parallel block-structured AMReX framework, provides functionality similar to the well-established Weather Research and Forecasting model. Leveraging the finite-difference method on the Arakawa C-grid, the level-based AMR in AMReX streamlines the integration of AMR infrastructure into existing NWP models. The framework's adaptability to exascale supercomputers and performance portability to GPUs and other accelerators further enhances its practicality. Notably, the familiarity of this level-based AMR method in static form within NWP models makes it a compelling choice for transitioning from static to dynamic mesh refinement.

The second method, a tree-based AMR, is implemented in NebulaSEM, a dynamical core employing the discontinuous Galerkin spectral element method. NebulaSEM, designed for simulating flow on the sphere, employs a non-overlapping, tree-based AMR with a forest of quad-/oct-trees managing refined and coarsened cells. While effective in accurately tracking features, this method presents challenges in terms of code reuse within NWP models accustomed to structured grids. Additionally, the application of multi-rate time-stepping methods to the tree-based AMR can be more intricate.

Verification of both AMR methods was conducted using established benchmark problems in NWP. The assessment of their efficiency relative to a uniform grid resolution solution demonstrated that both methods excel in accurately tracking localized and moving features. Moreover, they outperformed uniform grid simulations in terms of time-to-solution metric  for the same number of grid cells. These findings affirm the viability of both level-based and tree-based AMR methods for enhancing the predictive capabilities of NWP models. The choice between these methods may depend on factors such as familiarity of the approach, ease of integration, specific requirements of the model and existing code structure.

In addition, we have explored the issue of mass conservation encountered when utilizing AMR on non-conformal and high-order dGSEM discretization. The solution transfer between parent and children cells during refinement and coarsening, as well as the computation of fluxes at non-conformal faces, demand careful consideration to uphold the global conservation properties inherent to dGSEM. We have illustrated, through various examples, that simplistic interpolation techniques and inaccurate linearization methods (for the purpose of implicit-explicit treatment of terms) can compromise this conservation property. In contrast, the conservative projection method presented in this study, the ``mortar method" have demonstrated the ability  to conserve mass to close to machine precision. Additionally, we have devised a novel solution transfer strategy for simulations on the sphere that is easy to implement, memory-efficient and most importantly conservative. In this method, the projection matrices are still constructed only on the reference element and then used by all curvilinear elements on the sphere with corrections applied during refinement and coarsening.

In the realm of AMR for GCMs that operate on spherical grids, additional challenges emerge. Among the various grid types used in GCMs, the cubed-sphere grid stands out as particularly well-suited for AMR implementation. It seamlessly integrates with existing AMR libraries designed for quadrilateral/hexahedral grid refinement. We have uncovered that maintaining mass conservation of the dGSEM method on the sphere with AMR is more challenging than that in a box. This is primarily due to the fact that use of isoparametric curvilinear elements to represent the curved surface means that discretized geometry representation changes after refinement/coarsening. Through the application of two widely recognized benchmarks, we have showcased how the tree-based AMR method effectively tracks the transport of a tracer on the sphere and the propagation of an acoustic wave. We should note that while similar results can be anticipated with level-based methods, the fact that ERF is designed to be a limited-area model prevents us from conducting global simulations.   

\section*{Acknowledgment}
The first and fourth authors gratefully acknowledge funding from the Software Engineering for Novel Architectures (SENA) effort at GSL. The work of the second author at LBNL was supported by the U.S. Department of Energy under contract No. DE-AC02-05CH11231. Funding for the development of ERF was provided by the U.S. Department of Energy
Office of Energy Efficiency and Renewable Energy Wind Energy Technologies Office. The third author was funded by the Office of Naval Research under grant \# N0001419WX00721 and the National Science Foundation under grant AGS-1835881. 


\begin{appendices}

\section{Linearization error in advection equation}
\label{appendix-linear}
Here, we examine the linearization of the advection equation given by

\[
\frac{\partial \qb}{\partial t} + \nabla \cdot \left( \frac{\qb \mathbf U}{\rho} \right) = 0\\
\]
discretized into the weak form

\[
\left(v, \frac{\partial \qb}{\partial t}\right)_{\Omega_e} + \left\langle v \mathbf n, \Fm^* (\qb)\right\rangle_{\partial\Omega_e} - \left(\nabla v, \Fm(\qb)\right)_{\Omega_e} = 0.
\]

For the sake of linearization, we need to separate $\qb$ from $\mathbf U$ in $\Fm(\qb) = \qb\mathbf U /\rho = \Fm_m \qb$ where $\Fm_m=\mathbf U / \rho$ represents the mass flux.  To preserve the tracer conservation characteristics of the ``mortar method", both surface and volume integral terms should be treated consistently,  either explicitly or implicitly. Treating only one term implicitly results in a timestep-lagged value  of $\qb$ being used in the other term, leading to a loss of conservation.

The volume integral term can be correctly linearized due to the property 

\[
\left(\nabla v, \Fm(\qb)\right)_{\Omega_e} = \left(\nabla v, \Fm_m\right)_{\Omega_e}\qb.
\]
However, linearization is not feasible for the surface integral term because 

\[
\left\langle v \mathbf n, \Fm^*(\qb) \right\rangle_{\partial\Omega_e} \ne \left\langle v \mathbf n, \Fm_m^* *\qb^*\right\rangle_{\partial\Omega_e} = \left\langle v \mathbf n, \Fm_m^* \right\rangle_{\partial\Omega_e} \qb^*.
\]
However, this is exactly what is employed in finite-volume CFD introducing linearization error. The Rusanov flux used in this study has an average term $\{\Fm(q)\}$ that does not adhere to the product rule. Additionally, the ``mortar methods" gather and scatter operation complicates matters, necessitating the explicit treatment of the flux term. It's worth noting that the volume integral is zero for the finite volume method, rendering the linearization of the term inconsequential.

Consequently, if the surface integral term is treated explicitly without linearization, and the volume integral term is treated explicitly with or without linearization, tracer conservation is maintained by the ``mortar method" AMR. Conversely, if the volume integral term is implicitly treated with linearization, tracer conservation is compromised, even on a conformal grid, due to the use of lagged $\mathbf{q}$ in the explicitly treated surface integral term. 

\section{Solving a 2D problem on a curved surface}
\label{appendix-2d}
In order to solve a 2D problem, such as shallow atmosphere simulations, on a curved surface like the sphere, we modify the 3D Jacobian matrix as follows. Given physical coordinates $\mathbf x  = (x,y,z)$ and reference coordinates $\mathbf{\hat{x}} = (\xi,\eta,\zeta)$, the Jacobian matrix and its inverse

\begin{equation}
J = \begin{bmatrix} 
\frac{\partial x}{\partial \xi} &  \frac{\partial x}{\partial \eta} &  \frac{\partial x}{\partial \zeta}  \\
\frac{\partial y}{\partial \xi} &  \frac{\partial y}{\partial \eta} &  \frac{\partial y}{\partial \zeta}  \\
\frac{\partial z}{\partial \xi} &  \frac{\partial z}{\partial \eta} &  \frac{\partial z}{\partial \zeta}  \\
\end{bmatrix}
\quad and \quad
J^{-1} = \begin{bmatrix} 
\frac{\partial \xi}{\partial x} &  \frac{\partial \xi}{\partial y} &  \frac{\partial \xi}{\partial z}  \\
\frac{\partial \eta}{\partial x} &  \frac{\partial \eta}{\partial y} &  \frac{\partial \eta}{\partial z}  \\
\frac{\partial \zeta}{\partial x} &  \frac{\partial \zeta}{\partial y} &  \frac{\partial \zeta}{\partial z}  \\
\end{bmatrix}.
\end{equation}
can be used for coordinate transformations:
\[
\begin{aligned}
\mathbf{x}=J\mathbf{\hat{x}}\\
\mathbf{\hat{x}}=J^{-1}\mathbf{x}\\
\end{aligned}
\]
and computing derivatives of field $\phi$:
\[
\begin{aligned}
\nabla \phi={J^{-1}}^{\top}\hat{\nabla} \phi\\
\hat{\nabla} \phi=J^{\top}\nabla \phi.\\
\end{aligned}
\]
For 2D problems, the 3D Jacobian matrix becomes singular because $\partial x/\partial \zeta = \partial y/\partial \zeta = \partial z/\partial \zeta= 0$. Hence, we compute the the pseudo-inverse (also known as the Moore-Penrose inverse) $J^+$ to replace $J^{-1}$

\[
J^+={(J^{\top}J)}^{-1}J^{\top}
\]

\end{appendices}
\bibliographystyle{elsarticle-harv}
\bibliography{references,mypublications}

\end{document}